\documentclass[final]{siamltex}
\usepackage{xcolor}
\usepackage{amsmath}
\usepackage{amssymb}
\usepackage{mathrsfs}
\usepackage{graphicx}
\usepackage{bm}
\usepackage{ucs}
\usepackage[utf8x]{inputenc}
\usepackage{wrapfig}
\usepackage{color}
\usepackage{float}
\usepackage{epstopdf}
\usepackage{algorithm,algorithmic}
\usepackage[breaklinks,colorlinks=true,linkcolor=blue,citecolor=red,
  backref=page]{hyperref}

\DeclareMathAlphabet{\itbf}{OML}{cmm}{b}{it}
 \DeclareMathAlphabet\mathbfcal{OMS}{cmsy}{b}{n}
 
\def\EE{\mathbb{E}}
\def\GG{\boldsymbol{\mathcal{G}}}

\def\DD{{{\bf D}}}
\def\JJ{{\bf J}}
\def\FF{\boldsymbol{\mathscr{F}}}
\def\QQ{\boldsymbol{\mathcal{Q}}}
\def\UU{{\bf U}}
\def\tUU{\widetilde{\UU}}
\def\VV{{\bf V}}
\def\WW{\boldsymbol{\mathcal{W}}}
\def\HH{{\bf H}}

\def\bR{\boldsymbol{\mathcal{R}}}
\def\bSigma{\boldsymbol{\Sigma}}
\def\bV{{\bf V}}
\def\bx{{{\itbf x}}}
\def\bz{{{\itbf z}}}
\def\by{{{\itbf y}}}
\def\be{{{\itbf e}}}
\def\bd{{{\bf d}}}
\def\bu{{{\itbf u}}}
\def\bh{{{\itbf h}}}

\def\bW{{\bf W}}
\def\bT{{\bf T}}
\newcommand{\la}{\lambda}
\newcommand{\brho}{\boldsymbol{\rho}}
\newcommand{\ep}{\epsilon}
\newcommand{\de}{\delta}
\newcommand{\om}{\omega}
\newcommand{\bE}{{\itbf E}}
\newcommand{\vcE}{\vec{\boldsymbol{\mathcal{E}}}}
\newcommand{\vE}{\vec{\bE}}

\newcommand{\vx}{\vec{\bx}}

\newcommand{\vz}{\vec{\bz}}
\newcommand{\vy}{\vec{\by}}
\newcommand{\ve}{\vec{\be}}
\newcommand{\bG}{{\bf G}}

\newcommand{\bS}{{\bf S}}
\newcommand{\tDD}{\widetilde{\DD}}

\newcommand{\cA}{\mathcal{A}}
\newcommand{\vJ}{\vec{\itbf j}}
\newcommand{\ts}{\widetilde{\sigma}}
\newcommand{\tbu}{\widetilde{\itbf u}}
\newcommand{\bB}{{\bf B}}
\newcommand{\tbB}{\widetilde{\bB}}
\newcommand{\bDe}{\boldsymbol{\Delta}}
\newcommand{\bGamma}{\boldsymbol{\Gamma}}

\begin{document}

\title{Robust imaging with electromagnetic waves in  noisy environments} 
\author{Liliana Borcea\footnotemark[1] \and Josselin
  Garnier\footnotemark[2]}
\maketitle


\renewcommand{\thefootnote}{\fnsymbol{footnote}}
\footnotetext[1]{Department of Mathematics, University of Michigan,
  Ann Arbor, MI 48109. {\tt borcea@umich.edu}}
\footnotetext[2]{Laboratoire de Probabilit\'es et Mod\`eles
  Al\'eatoires \& Laboratoire Jacques-Louis Lions, Universit{\'e}
  Paris Diderot, 75205 Paris Cedex 13, France.  {\tt
    garnier@math.univ-paris-diderot.fr}}
\markboth{L. BORCEA, J. GARNIER}{ARRAY IMAGING}

\begin{abstract}
We study imaging with an array of sensors that probes a medium with
single frequency electromagnetic waves and records the scattered
electric field. The medium is known and homogenous except for some
small and penetrable inclusions.  The goal of inversion is to locate
and characterize these inclusions from the data collected by the
array, which are corrupted by additive noise.  We use results from
random matrix theory to obtain a robust inversion method.  We assess
its performance with numerical simulations and quantify the benefit of
measuring more than one component of the scattered electric field.
\end{abstract}

\section{Introduction}
We study an inverse problem for Maxwell's equations, where unknown
scatterers in a medium are to be determined from data collected with
an array of sources and receivers. The sources emit signals which
propagate in the medium and the receivers record the backscattered
waves. The recordings, indexed by the source and receiver pair, are
the entries in the response matrix, the data, which are contaminated with additive noise.

Sensor array imaging is an important technology in nondestructive
testing, medical imaging, radar, seismic exploration, and
elsewhere. It operates over application specific ranges of frequencies
and fixed or synthetic apertures. We consider imaging
with single frequency waves and a fixed array of $N$ sensors that play
the dual role of sources and receivers. The unknown scatterers lie in
a homogeneous medium. They are $P$ penetrable inclusions of small size
with respect to the wavelength. As shown in
\cite{ammari2007music} they can be described by their center location
$\vy_p$ and reflectivity tensor $\brho_p$ which depends on the volume
and shape of their support, and their electric permittivity $\ep_p$,
for $p = 1, \ldots, P$.

There are many algorithms for locating scatterers from the response
matrix. Reverse time or Kirchhoff migration \cite{biondi20063d} and
the similar filtered backprojection \cite{cheney2009fundamentals} are
common in seismic and radar imaging with broadband data. Matched field
processing \cite{baggeroer1993overview} is popular in underwater
acoustics. Qualitative methods such as linear sampling
\cite{cakoni2011linear}, the factorization method
\cite{kirsch2008factorization}, and the MUSIC (Multiple Signal
Classification) algorithm \cite{schmidt1986multiple} are useful for
imaging with single frequency waves, and give high resolution
results if the noise  is  weak.  

A study of noise effects on the MUSIC method can be found in
\cite{fannjiang2011music}. It determines upper bounds on the noise
level that guarantee the exact recovery of point scatterers in a sparse
scene.  Here we consider stronger noise and use  results from random matrix theory
\cite{benaych2011eigenvalues,capitaine2012central,baik2005phase,
  marchenko1967distribution} to obtain a robust MUSIC type localization method. Such results have led to scatterer
detection tests in \cite{ammari2011imaging,ammari2012statistical}, and to filters 
of backscattered waves in random
media in \cite{aubry2009detection,alonso2011detection}.  They are also the foundation of  a robust sonar 
array imaging method of well separated  point scatterers in \cite{garnier2014applications}.
 We
extend the results in \cite{garnier2014applications} to
electromagnetic waves, for locating one or more inclusions and estimating
their reflectivity tensor. We also quantify the  benefits of measuring more than one 
component of the scattered electric field.

The paper is organized as follows: We begin in section \ref{sect:form}
with the mathematical formulation of the inverse problem, and then
give in section \ref{sect:data} the data model. The singular value
decomposition analysis of the response matrix is described in section
\ref{sect:inv}. The inversion method uses this analysis and is presented in section
\ref{sect:inv.alg}. Numerical
simulation results are shown in section \ref{sect:num}. We end with a summary in
section \ref{sect:sum}.
\section{Formulation of the problem}
\label{sect:form}
Consider a homogeneous isotropic medium with electric permittivity
$\ep_o$ and magnetic permeability $\mu_o$, which contains $P$
penetrable inclusions supported in the disjoint simply connected
domains $\Omega_p$ with smooth boundaries, centered at locations
$\vy_p$, for $p = 1,\ldots,P$. The inclusions are modelled by the
piecewise constant electric permittivity
\begin{equation}
\ep(\vx) = \sum_{p=1}^P \ep_p 1_{\Omega_p}(\vx) + \ep_o
1_{\mathbb{R}^3 \setminus \Omega}(\vx), \qquad \Omega =
\bigcup_{p=1}^P \Omega_p,
\label{eq:f1}
\end{equation}
where $1_\Omega(\vx)$ is the indicator function of the domain $\Omega$,
and $\ep_p$ are scalar real valued for $p = 1, \ldots, P$.

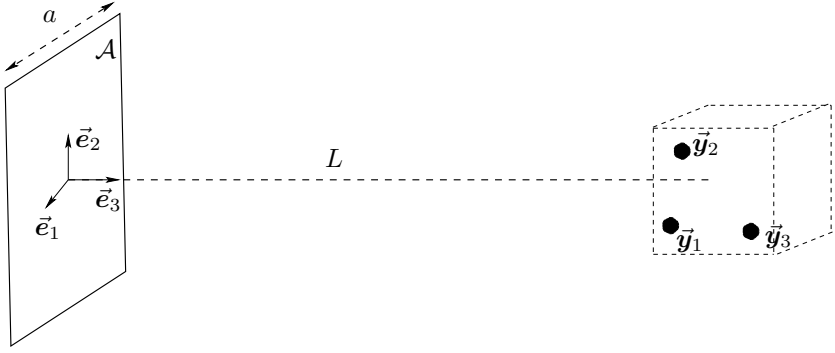
\begin{figure}[t]
\centering
  \begin{picture}(0,100)%
\hspace{-2in}\includegraphics[width = 0.85\textwidth]{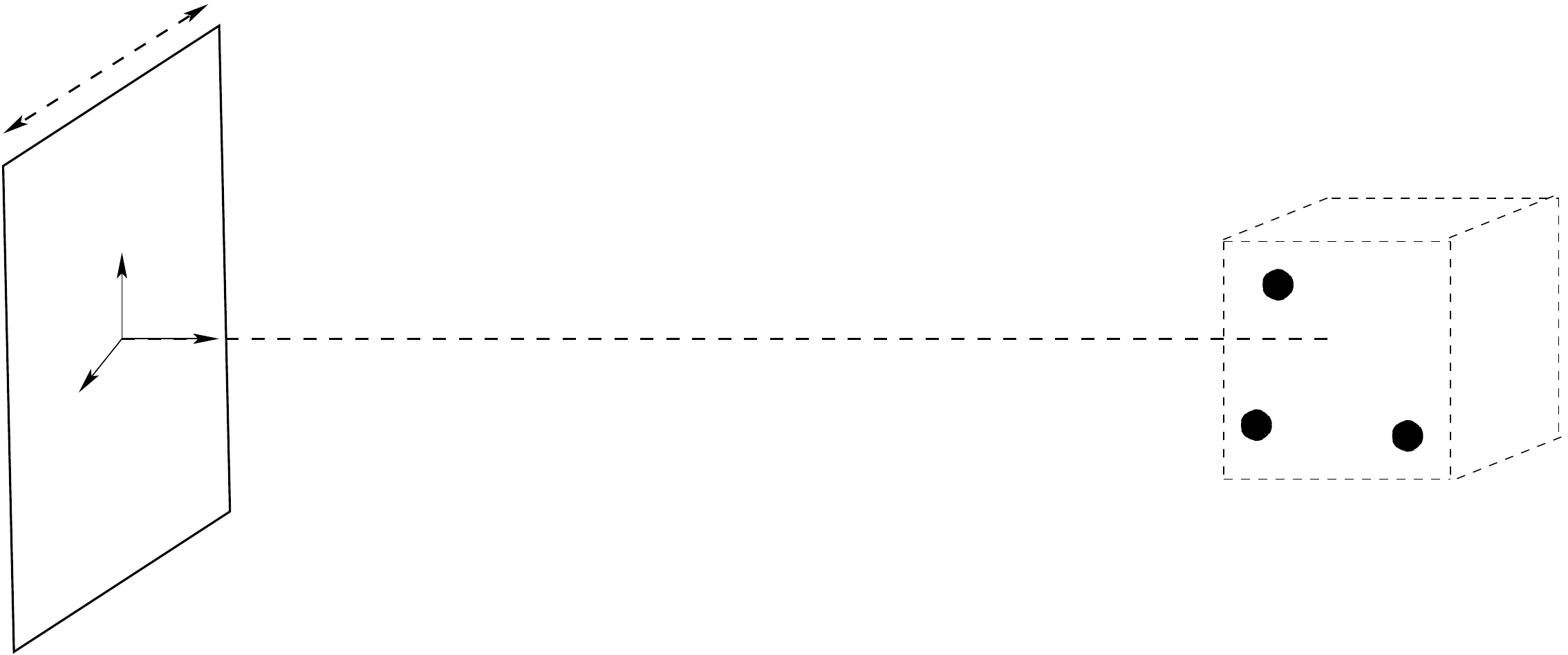}%
\end{picture}%
\setlength{\unitlength}{100sp}%
\begingroup\makeatletter\ifx\SetFigFont\undefined%
\gdef\SetFigFont#1#2#3#4#5{%
  \reset@font\fontsize{#1}{#2pt}%
  \fontfamily{#3}\fontseries{#4}\fontshape{#5}%
  \selectfont}%
\fi\endgroup%
\begin{picture}(0,100)(0,100)
\put(-72000,72000){\makebox(0,0)[lb]{\smash{{\SetFigFont{7}{8.4}{\familydefault}{\mddefault}{\updefault}{\color[rgb]{0,0,0}{\normalsize $\cA$}}%
}}}}
\put(-85000,81000){\makebox(0,0)[lb]{\smash{{\SetFigFont{7}{8.4}{\familydefault}{\mddefault}{\updefault}{\color[rgb]{0,0,0}{\normalsize $a$}}%
}}}}
\put(-15000,45000){\makebox(0,0)[lb]{\smash{{\SetFigFont{7}{8.4}{\familydefault}{\mddefault}{\updefault}{\color[rgb]{0,0,0}{\normalsize $L$}}%
}}}}
\put(-77000,50000){\makebox(0,0)[lb]{\smash{{\SetFigFont{7}{8.4}{\familydefault}{\mddefault}{\updefault}{\color[rgb]{0,0,0}{\normalsize $\ve_2$}}%
}}}}
\put(-87000,28000){\makebox(0,0)[lb]{\smash{{\SetFigFont{7}{8.4}{\familydefault}{\mddefault}{\updefault}{\color[rgb]{0,0,0}{\normalsize $\ve_1$}}%
}}}}
\put(-72000,35000){\makebox(0,0)[lb]{\smash{{\SetFigFont{7}{8.4}{\familydefault}{\mddefault}{\updefault}{\color[rgb]{0,0,0}{\normalsize $\ve_3$}}%
}}}}
\put(72000,25000){\makebox(0,0)[lb]{\smash{{\SetFigFont{7}{8.4}{\familydefault}{\mddefault}{\updefault}{\color[rgb]{0,0,0}{\normalsize $\vy_1$}}%
}}}}
\put(76000,49000){\makebox(0,0)[lb]{\smash{{\SetFigFont{7}{8.4}{\familydefault}{\mddefault}{\updefault}{\color[rgb]{0,0,0}{\normalsize $\vy_2$}}%
}}}}
\put(94000,26000){\makebox(0,0)[lb]{\smash{{\SetFigFont{7}{8.4}{\familydefault}{\mddefault}{\updefault}{\color[rgb]{0,0,0}{\normalsize $\vy_3$}}%
}}}}
\end{picture}%
\caption{Illustration of the inversion setup. The array $\cA$ is
  square planar with side $a$, in the cross-range plane spanned by
  $\ve_1$ and $\ve_2$. The imaging region is at range scale $L$ from
  the array. In the illustration there are three inclusions at
  locations $\vy_p$, for $p = 1, 2, 3.$}
\label{fig:setup}
\end{figure}

The medium is probed with $N$ sensors that are both sources and
receivers. They are closely spaced at locations $\{\vx_r\}_{1\le r \le
  N} $ in a set $\cA$, so they behave like a collective entity, the
array. For convenience we let $\cA$ be planar and square, of side $a$,
the array aperture. This allows us to introduce a system of
coordinates with origin at the center of the array and orthonormal
basis $\{\ve_1,\ve_2,\ve_3\}$. The vectors $\ve_1, \ve_2$ span the plane
of the array, named  the cross-range plane, and the line along
$\ve_3$ is the range direction.  The distance between the array and the
inclusions is of order $L$, the range scale. See Figure
\ref{fig:setup} for an illustration.

We index consistently a source sensor by $s$ and a receiver sensor by
$r$, altough both $s$ and $r$ take values in the same set $\{1,
\ldots, N\}$. The sources are point dipoles with current density
$\vJ_s \de(\vx-\vx_s)$, and the resulting electric field $\vE$
satisfies
\begin{align}
\vec{\nabla} \times \vec{\nabla} \times \vE(\vx;\vx_s) - k^2
  \frac{\ep(\vx)}{\ep_o} \vE(\vx;\vx_s) = i k
\sqrt{\frac{\mu_o}{\ep_o}}\,\vJ_s \de (\vx-\vx_s),
\label{eq:f3}
\end{align}
and the radiation condition
\begin{equation}
\lim_{|\vx| \to \infty} |\vx| \left[ \vec{\nabla} \times
  \vE(\vx;\vx_s) - i k \frac{\vx}{|\vx|} \times \vE(\vx;\vx_s) \right]
= 0.
\label{eq:f5}
\end{equation}
Equation \eqref{eq:f5} is derived from Maxwell's equations at
frequency $\om$, and
\[
k = \frac{\om}{c_o} = \frac{2 \pi}{\la}
\]
is the wavenumber, where $c_o = 1/\sqrt{\mu_o \ep_o}$ is the wave
speed and $\la$ is the wavelength.

We decompose $\vE$ in the direct field $\vE_o$ and the scattered one
denoted by $\vcE$,
\begin{equation}
\vcE(\vx;\vx_s) = \vE(\vx;\vx_s) - \vE_o(\vx;\vx_s).
\end{equation}
The direct field solves the analogue of \eqref{eq:f3}--\eqref{eq:f5}
in the homogeneous medium, and can be written explicitly as
\begin{equation}
 \vE_o(\vx;\vx_s) = i k
\sqrt{\frac{\mu_o}{\ep_o}} \,\bG(\vx,\vx_s) \vJ_s , \label{eq:f3o}
\end{equation}
using  the dyadic Green's tensor
\begin{equation}
\bG(\vx,\vz) = \left( {\bf I}_3 + \frac{\nabla \nabla^T}{k^2} \right)
\frac{e^{i k |\vx-\vz|}}{4 \pi |\vx-\vz|}.
\label{eq:GTens}
\end{equation}
The scattered field $\vcE$ is recorded at the $N$ receivers, and the
experiment may be repeated by emitting waves sequentially from all the
sources.  Each source may use up to three illuminations, with linearly
independent $\vJ_s \in \mathbb{R}^3$, and the receivers may record one
or all three components of $\vcE$. We call the data for all possible
illuminations and recordings \emph{complete}. All other cases, with
one or two illuminations from each source and one or two components of
$\vcE$ measured at the receivers are called \emph{incomplete}. The
goal is to locate and characterize the inclusions using the complete
or incomplete measurements.

We simplify the problem by assuming that the inclusions are small with
respect to the wavelength, so we can represent them as in
\cite{ammari2007music} by their center locations $\vy_p$ and $3\times
3$ reflectivity tensors $\brho_p$, for $p = 1, \ldots, P$.  These
tensors depend on the shape of the domains $\Omega_p$ and the relative
electric permittivity $\ep_p/\ep_o$. The inverse problem is to
determine the locations $\vy_p$ and reflectivity tensors $\brho_p$ from the array
measurements.

\section{Data model}
\label{sect:data}
Let $\alpha$ be a small and positive dimensionless parameter which
scales the size of the inclusions
\begin{equation}
\frac{\left[{\rm Vol}(\Omega_p)\right]^{1/3}}{\la} \sim \alpha \ll 1,
\quad p = 1, \ldots, P,
\label{eq:f4}
\end{equation}
and denote by $\vcE_q(\vx;\vx_s)$ the scattered  field due to
the source at $\vx_s$, with $\vJ_s = \ve_q$. It has the following
asymptotic expansion in the limit $\alpha \to 0$, as shown
in \cite{ammari2007music},
\begin{equation}
\vcE_q(\vx;\vx_s) = i k^3 \sqrt{\frac{\mu_o}{\ep_o}}\sum_{p=1}^P
\bG(\vx,\vy_p) \brho_p \bG(\vy_p,\vx_s) \ve_q 
+O(\alpha^4) \, .
\label{eq:dm1}
\end{equation}
The reflectivity tensor of the $m-$th
inclusion is defined by
\begin{equation}
\brho_p = \alpha^3 \left(\frac{\ep_p}{\ep_o} - 1 \right)
{\bf M}_p,
\label{eq:dm2}
\end{equation}
using the $\alpha$-independent polarization tensor
${\bf M}_p$. We refer to
\cite{ammari2007music} for details on
${\bf M}_p$, and to section
\ref{sect:num} for its calculation in our numerical simulations. 
It suffices to say that ${\bf M}_p$
is a symmetric, real valued $3 \times 3$ matrix which depends on the
shape of $\Omega_p$ and the relative electric permittivity
$\ep_p/\ep_o$.

\subsection{The noiseless data model}
\label{sect:Data.1}
Since $k$, $\mu_o$ and $\ep_o$ are known, let us model the noiseless
data gathered with the receiver and source pair $(r,s)$ by
\begin{equation}
\bd(\vx_r,\vx_s) = \frac{\sqrt{\ep_o/\mu_o}}{ik^3}
\left(\vcE_1,\vcE_2,\vcE_3\right), \quad r, s = 1, \ldots, N.
\label{eq:dm3}
\end{equation} 
This is a $3 \times 3$ complex matrix, the $(r,s)$ block of the $3 N
\times 3N$ \emph{complete data matrix} $\DD$ collected by the array,
\begin{equation}
\DD = \left(  \bd(\vx_r,\vx_s) \right)_{r,s=1,\ldots, N}.
\label{eq:dm4}
\end{equation}
To acquire $\DD$ directly, the array would use the sources
sequentially, with current densities $\ve_q \delta(\vx - \vx_s)$, for
$q =1, 2, 3$ and $s = 1, \ldots, N$, and would measure  the scattered fields $\vcE_q(\vx;\vx_s)$ at the
receivers. This acquisition method is not optimal in noisy environments, 
and may be improved as explained in the next section.

When the sensors emit and receive only along one or two directions, we
have incomplete measurements. Let us denote by $\mathscr{S} \subset
\{1,2,3\}$ the set of indices of the measured components of the electric field, and
define the $3 \times |\mathscr{S}|$ \emph{sensing matrix}
\begin{equation}
\bS = \left(\ve_q \right)_{q \in \mathscr{S}},
\label{eq:dm5}
\end{equation}
with columns given by the unit vectors $\ve_q$ with $q \in
\mathscr{S}$.  We assume for convenience that the source excitation
currents are also limited to these directions, and denote the
incomplete data matrix by $\DD_{_{\hspace{-0.02in}S}}$. It is a
complex $|\mathscr{S}| N \times |\mathscr{S}| N$ matrix that can be
written in terms of $\DD$ as
\begin{equation}
\DD_{_{\hspace{-0.02in}S}} = \mbox{diag} \left(\bS^T, \ldots,
\bS^T\right) \DD \mbox{diag} \left(\bS, \ldots, \bS\right).
\label{eq:dm6}
\end{equation}

The mathematical model of $\DD$ follows from \eqref{eq:dm1} and
\eqref{eq:dm3}.  We write it using the forward map $\FF$, which takes
the unknown locations and reflectivity tensors of the inclusions to the
space $\mathbb{C}^{3N \times 3N}$ of matrices where $\DD$ lies,
\begin{align}
\DD &\approx \FF(\vy_{1}, \ldots, \vy_{P}, \brho_{1}, \ldots, \brho_{P})
= \left( \sum_{p=1}^P \bG(\vx_r,\vy_p) \brho_p
\bG(\vy_p,\vx_s) \right)_{r,s = 1, \ldots, N.}
\label{eq:dm7}
\end{align}
The forward map  for the incomplete
measurements  is deduced from \eqref{eq:dm6}--\eqref{eq:dm7}. 
We assume henceforth that $\FF$ gives a good approximation of $\DD$, as
is the case for sufficiently small $\alpha$, and treat the
approximation in \eqref{eq:dm7} as an equality.

\subsection{Noisy data and the Hadamard acquisition scheme}
\label{sect:Data.2}
In any practical setting the measurements are contaminated with noise,
which may be mitigated by the data acquisition scheme, as we now
explain.

Suppose that the $N$ sources do not emit waves sequentially, but
at the same time, for $3N$ experiments indexed by $(n,q)$, with $n =
1, \ldots, N$ and $q = 1, 2,3.$ In the $(n,q)$ experiment the
excitation from the source at $\vx_s$ is $\vJ_{s,(n,q)}$, and
gathering the measurements for $q = 1, 2, 3$, we obtain the following
model of the recordings at the receiver location $\vx_r$,
\begin{equation}
{\bf d}^{\rm H}_{r,n} = \sum_{s = 1}^N
\bd(\vx_r,\vx_s) {\bf j}_{s,n} +
{\bf w}_{r,n}.
\label{eq:dm9p}
\end{equation}
This is a $3\times 3$ complex matrix, determined by $\bd(\vx_r,\vx_s)$
defined in \eqref{eq:dm3}, the excitation matrices ${\bf j}_{s,n}$ with columns $\vJ_{s,(n,q)}$ for $q = 1, 2,3$, and the
noise matrix ${\bf w}_{r,n}$.  The complete $3N \times 3N$ measurement
matrix ${\bf D}^{\rm H}$  has the block structure
\[
{\bf D}^{\rm H} = \left( {\bf d}^{\rm H}_{r,n} \right)_{r,n=1,\ldots, N},
\]
and it is modeled by 
\begin{equation}
{\bf D}^{\rm H} = \DD {{\JJ}} + \WW,
\label{eq:dm9}
\end{equation}
where ${{\JJ}}$ is the  excitation matrix with blocks
${\bf j}_{s,n}$ for $s,n = 1, \ldots, N$, and $\WW$ is the noise matrix with blocks
${\bf w}_{r,n}$ for $r,n = 1, \ldots, N$.

We assume that the excitation currents are normalized so that the
entries in ${\JJ}$ have unit amplitudes, and take as is usual a noise
matrix $\WW$ with mean zero and  independent complex Gaussian entries
with standard deviation $\sigma$. This gives 
\[
\EE \left[\WW^\dagger \WW\right] = \sigma^2 {\bf  I}_{_{3N}},
\]
where $\dagger$ denotes the complex conjugate and transpose, $\EE$ is
the expectation, and ${\bf I}_{_{3N}}$ is the $3N \times 3N$ identity
matrix.  The acquisition scheme uses an invertible
 $\boldsymbol{{\JJ}}$, to obtain
\begin{equation}
\widetilde \DD = {\bf D}^{\rm H}{\JJ}^{-1} =
\DD + \bW, \qquad \bW = \WW {\JJ}^{-1}.
\label{eq:dm10}
\end{equation}
This  is a corrupted version of $\DD$, but
the noise level in  $\bW$, which has mean zero complex Gaussian
entries, is reduced for a good choice of ${{\JJ}}$.  Such a
choice would give
\begin{equation}
\EE \left[\bW^\dagger \bW\right] =
\big({\JJ}^{-1}\big)^\dagger\EE \left[ \WW^\dagger
  \WW \right]{\JJ}^{-1} = \sigma^2
\big({\JJ}^{-1}\big)^\dagger {\JJ}^{-1} =
\sigma_{_{\bW}}^2{\bf I}_{_{3N}},
\label{eq:dm11}
\end{equation}
and a minimum $\sigma_{_\bW}$.  We conclude that $\JJ$ must be 
unitary, up to a multiplicative constant factor, 
and its determinant should be maximal to get the  smallest
$$
\sigma_{_{\bW}} = \frac{\sigma}{|\det ({{\JJ}})|^{1/(3N)}}.
$$

It is shown in \cite{hadamard1893resolution} that the determinant of
$3N \times 3N$ matrices with entries in the complex unit disk is
bounded above by $(3N)^{3N/2}$, with equality attained by complex
Hadamard matrices. These have entries with modulus one, and are
unitary up to the scaling by their determinant.  Thus, by using a
Hadamard matrix ${\JJ}$, we can mitigate the noise in
$\widetilde{\DD}$ and reduce its standard deviation to
\begin{equation}
\sigma_{_{\bW}} = \sigma/\sqrt{3N}.
\end{equation}

The results extend verbatim to the incomplete measurement setup. The
difference is that in (\ref{eq:dm9p}) we have $\bS^T {\bf
  d}(\vx_r,\vx_s) \bS$ instead of ${\bf d}(\vx_r,\vx_s)$, and the $3
\times 3$ current excitation matrix ${\bf j}_{s,n}$ is replaced by
$\bS^T {\bf j}_{s,n} \bS$, the blocks in the complex $|\mathscr{S}|
N \times |\mathscr{S}| N$ Haddamard matrix. The noise level in
\begin{equation}
\label{eq:dm11S}
\widetilde{\DD}_{_{\hspace{-0.02in} S}} = \DD_{_{\hspace{-0.02in} S}} +
\bW_{_{\hspace{-0.02in} S}}
\end{equation}
is reduced to  $\sigma/\sqrt{|\mathscr{S}|N}.$

\vspace{0.1in} \textbf{Remark:} In case that the noise  at 
nearby receivers is correlated, with known correlation structure, let us denote by $\boldsymbol{\mathcal{C}} = \EE
\left[ \WW^\dagger \WW \right] $  the Hermitian covariance matrix,  and suppose that it is positive definite, with  $O(
(\sigma^{2})^{3N})$ determinant.  Then we can
use its square root to define an invertible matrix $\JJ$ so that
$\boldsymbol{\mathcal{C}}^{1/2} \JJ^{-1}$ is 
unitary up to scaling by its determinant, and   obtain that $\bW$ in \eqref{eq:dm10} has uncorrelated entries
with $\sigma_{_{\bW}} = O\left(\sigma/|\mbox{det}(
\JJ)|^{1/(3N)}\right)$.

\section{Singular value decomposition analysis}
\label{sect:inv}
Our imaging algorithm is based on the Singular Value Decomposition
(SVD) of the response matrix.  We begin in section \ref{sect:inv.1}
with the SVD analysis in the noiseless case.  This uses the rank
calculation of the matrix of Green's functions given in section
\ref{sect:inv.2}. The effect of the additive noise on the SVD is
discussed in section \ref{sect:inv.3}.

\subsection{SVD analysis of the noiseless data matrix}
\label{sect:inv.1}
Let us introduce the $3N \times 3$ matrix $\GG(\vy)$ of Green's
functions, with block structure
\begin{equation}
\GG(\vy) = \begin{pmatrix} \bG(\vx_1,\vy) \\ \vdots
  \\ \bG(\vx_N,\vy)\end{pmatrix},
\label{eq:in2}
\end{equation}
and use the reciprocity relation $ \bG(\vx,\vy) = \bG(\vy,\vx)^T$ to
write the model in \eqref{eq:dm7} as
\begin{equation}
\DD = \sum_{p=1}^P \GG(\vy_p) \brho_p \GG(\vy_p)^T.
\label{eq:in1}
\end{equation}
We show in section \ref{sect:inv.2} that in our
setting $\GG(\vy)$ is full rank, so let
\begin{equation}
\GG(\vy) = \HH(\vy) \bSigma(\vy) \VV(\vy)^\dagger
\label{eq:in3}
\end{equation}
be its SVD, with $\HH \in \mathbb{C}^{3N \times 3}$ satisfying $
\HH^\dagger \HH = {\bf I}_3, $ and unitary $\VV \in \mathbb{C}^{3
  \times 3}$. The singular values of $\GG(\vy)$ are in the diagonal
positive definite $3 \times 3$ matrix $\bSigma(\vy)$. 
Substituting  \eqref{eq:in3} in \eqref{eq:in1} gives
\begin{equation}
\DD = \sum_{p=1}^P \HH(\vy_p) \bR_p \HH(\vy_p)^T,
\label{eq:in4}
\end{equation}
with complex symmetric $3 \times 3$ matrices 
\begin{equation}
\bR_p =\bSigma(\vy_p) \VV(\vy_p)^\dagger \brho_p \overline{\VV(\vy_p)}
\bSigma(\vy_p),
\label{eq:in5}
\end{equation}
where the bar denotes complex conjugate. Moreover, 
if 
\begin{equation}
\bR_p = \boldsymbol{\mathfrak{U}}_p \boldsymbol{\mathfrak{S}}_p
\boldsymbol{\mathcal{V}}_p^\dagger
\label{eq:in6}
\end{equation}
is an SVD of $\bR_p$, we obtain from \eqref{eq:in4} that 
\begin{equation}
\DD = \sum_{p=1}^P \HH(\vy_p) \boldsymbol{\mathfrak U}_p \boldsymbol{\mathfrak{S}}_p
\boldsymbol{\mathcal V}_p^\dagger \HH(\vy_p)^T.
\label{eq:in7}
\end{equation}

Our imaging algorithm  uses the relation
between the column space (range) of $\GG(\vy)$, or equivalently
$\HH(\vy)$, and the subspace spanned by the left singular vectors of
$\DD$. This follows from \eqref{eq:in7} and depends on the number of
inclusions and their separation distance, as we now explain.

\vspace{0.06in} \noindent \textbf{Single inclusion.} When $P =1$ we
see from \eqref{eq:in7} that $\DD$ has an SVD with nonzero singular
values $\sigma_1, \sigma_2, \sigma_3,$ the entries in
$\boldsymbol{\mathfrak{S}}_1$, and a $3N \times 3$ matrix of
corresponing left singular vectors $\HH(\vy_1) \boldsymbol{\mathfrak
  U}_1$. However, even when the singular values are distinct, 
  the singular vectors are defined up to an arbitrary
phase, so we cannot assume that the computed $3N \times 3$ matrix $\UU$ of left
singular vectors of $\DD$ equals $\HH(\vy_1) \boldsymbol{\mathfrak
  U}_1$. We work instead with projection matrices which are uniquely
defined, and satisfy
\begin{equation}
\UU \UU^\dagger = \HH(\vy_1) \boldsymbol{\mathfrak
  U}_1\boldsymbol{\mathfrak U}_1^\dagger \HH(\vy_1)^\dagger =
\HH(\vy_1) \HH(\vy_1)^\dagger,
\label{eq:in8}
\end{equation}
because $\boldsymbol{\mathfrak U}_1$ is unitary.  We conclude that the
column space of $\DD$ is the same as that of $\GG(\vy_1)$, and
moreover, that the $3N \times 3N$ projection on this space can be
calculated as $\UU \UU^\dagger $ using the left singular vectors of
$\DD$, or equivalently as $\HH(\vy_1) \HH(\vy_1)^\dagger$, using the
left singular vectors of $\GG(\vy_1)$.

\vspace{0.06in} \noindent \textbf{Multiple inclusions.} When $1 < P <
N$, the rank of $\DD$ depends on the locations of the inclusions
\cite{devaney2000super}. Generically\footnote{The rank may be smaller
  than $3P$ for very special locations of the inclusions
  \cite{devaney2000super}.}  it equals $3P$, as we assume here to
simplify the presentation. Let $\sigma_1, \ldots, \sigma_{3P}$ be the
nonzero singular values of $\DD$ and
\[
\UU = (\bu_1, \ldots, \bu_{3P})
\] 
the $3N \times 3P$ matrix of its corresponding left
singular vectors. The projections $\HH(\vy_p)\HH(\vy_p)^\dagger$ are
no longer the same as those obtained from $3N \times 3$ blocks in
$\UU$, unless the inclusions are so far apart that the column spaces
of $\GG(\vy_p)$ and $\GG(\vy_{p'})$ are orthogonal for $ p \ne p'.$
Nearby inclusions interact and the singular values and vectors of
$\DD$ are not associated with a single inclusion. However, we conclude
from \eqref{eq:in1} and \eqref{eq:in7} that
\begin{equation}
\mbox{range} \left( \GG(\vy_p) \right) = \mbox{range} \left( \HH
(\vy_p) \right) \subset \mbox{range}(\DD) = \mbox{span} \{ \bu_1, \ldots, \bu_{3P}\}. \label{eq:in8m}
\end{equation}
Moreover, it is shown in \cite[Proposition 4.3]{ammari2007music} that
\begin{equation}
\GG(\vy) \ve \in \mbox{range}(\DD)  \quad \mbox{iff} ~ \vy \in \{\vy_1, \ldots, \vy_P\},
\label{eq:in.8}
\end{equation}
for any $\ve$ such that $\GG(\vy) \ve \ne 0$. In our setting
$\GG(\vy)$ has full rank, so \eqref{eq:in.8} holds for any unit vector
$\ve \in \mathbb{R}^3.$ This implies in particular that the left
singular vectors of $\GG(\vy)$, the columns of $\HH(\vy)$, are in the
range  of $\DD$ if and only if $\vy$
coincides with the location of an inclusion.

The results are very similar for incomplete measurements, with data
modeled by
\begin{equation}
\DD_{_{\hspace{-0.02in} S}} = \sum_{p=1}^P \GG_{_{\hspace{-0.02in}
    S}}(\vy_p) \brho_p \GG_{_{\hspace{-0.02in} S}}(\vy_p)^T, 
\label{eq:in1S}
\end{equation}
in terms of the $|\mathscr{S}|N \times 3$ matrices of Green's functions
\begin{equation}
\GG_{_{\hspace{-0.02in} S}}(\vy) = \mbox{diag}\big(\bS^T,\ldots,
\bS^T\big) \GG(\vy). 
\label{eq:in2S}
\end{equation}
We show in the next section that the rank of $\GG_{_{\hspace{-0.02in}
    S}}$ is still three, as in the complete measurement case, but its
condition number may be much worse. This plays a role in the inversion
with noisy measurements.  We denote with the same symbol $\UU$ the
matrix of left singular vectors of $\DD_{_{\hspace{-0.02in} S}}$ and
$\HH(\vy)$ the matrix of left singular vectors of
$\GG_{_{\hspace{-0.02in} S}}(\vy)$. We conclude as in the complete
measurement case that when $P = 1$, 
\[
\UU \UU^\dagger = \HH(\vy_1)\HH(\vy_1)^\dagger,
\]
and when $1 < P < N$, using \cite[Proposition 4.3]{ammari2007music},
\begin{equation}
\mbox{range} \left( \GG_{_{\hspace{-0.02in} S}}(\vy) \right) = \mbox{range}
\left(\HH(\vy) \right) \subset \mbox{range}
\left(\DD_{_{\hspace{-0.02in} S}}\right) \quad \mbox{iff} ~ \vy \in
\{\vy_1, \ldots, \vy_P\}.
\end{equation}

\subsection{Rank and conditioning of matrix of Green's functions}
\label{sect:inv.2}
As mentioned above and explained in section \ref{sect:inv.alg}, the
rank and condition number of $\GG(\vy)$ play a role in imaging with
noisy measurements. We analyze them here using the Hermitian matrix
\begin{align}
\QQ(\vy) &= \GG(\vy)^\dagger \GG(\vy) = \sum_{r=1}^N
\bG(\vx_r,\vy)^\dagger \bG(\vx_r,\vy) \nonumber \\
&\approx \sum_{r=1}^N
\frac{1}{\left( 4 \pi|\vx_r-\vy| \right)^2} \left[ {\bf I}_3 -
  \frac{(\vx_r-\vy)}{|\vx_r-\vy|} \frac{(\vx_r-\vy)^T}{|\vx_r-\vy|}
  \right],
\label{eq:ra1}
\end{align}
whose eigenvalues equal the square of the singular values of
$\GG(\vy)$, the entries in $\bSigma(\vy)$. The approximation
in \eqref{eq:ra1} uses definition \eqref{eq:GTens} of the Green's
tensor, written as
\begin{equation}
\bG(\vx_r,\vy) = \frac{e^{ik|\vx_r-\vy|}}{4 \pi |\vx_r-\vy|} \left[
  {\bf I}_3 - \frac{(\vx_r-\vy)}{|\vx_r-\vy|}
  \frac{(\vx_r-\vy)^T}{|\vx_r-\vy|} +
  O\hspace{-0.02in}\left(\frac{1}{k L}\right) \right], \label{eq:ra1p}
\end{equation}
and neglects the $O(1/(kL))$ residual which is small when the
inclusions are many wavelengths away from the array.  

Equation \eqref{eq:ra1} shows that $\QQ(\vy)$ is, up to some scalar
factors, the sum of projection matrices on the plane orthogonal to
$\vx_r-\vy$. As $\vx_r$ varies in the array aperture, $\vx_r-\vy$
changes direction in a cone of opening angle of order $a/L$. The
smaller this angle, the worse the conditioning of $\QQ(\vy)$.  
Then let us consider the small aperture regime $a \ll L$, where
we can approximate the right hand side in \eqref{eq:ra1} by
\begin{align}
\QQ(\vy) \approx \frac{N}{(4 \pi L)^2}\hspace{-0.03in} \left[{\bf
  I}_3 - \ve_3 \ve_3^T - \frac{a}{L N} \sum_{r=1}^N\hspace{-0.03in}
 \begin{pmatrix} \frac{az_{r,1}^2}{L} & \frac{az_{r,1}
      z_{r,2}}{L} & -z_{r_1}\vspace{0.03in} \\ \vspace{0.03in}
    \frac{az_{r,1} z_{r,2}}{L} & \frac{az_{r,2}^2}{L} & - z_{r_2} \\ -
    z_{r_1} & - z_{r_2} & - \frac{a|{\itbf z}_r|^2 }{L}\end{pmatrix} \right],
\label{eq:ra2}
\end{align}
using the notation 
\[{\itbf z}_r = (\bx_r - \by)/a =
(z_{r,1},z_{r,2}), 
\] for $\vx_r = (\bx_r,0)$ and $\vy = (\by,L)$ satisfying $|\by| \lesssim a
\ll L.$ Here we neglected an $O(a^3/L^3)$ matrix residual in the sum over
the $N$ sensors.   The leading part of \eqref{eq:ra2},
${N}/{(4 \pi L)^2} [ {\bf I}_3 - \ve_3 \ve_3^T]$, 
is proportional to
the orthogonal projection on the cross-range plane, and determines the
two larger eigenvalues of $\QQ(\vy)$, which are approximately $N/(4
\pi L)^2$. The third eigenvalue is much smaller, by a factor of order
$a^2/L^2$. Thus, $\QQ(\vy)$ is full rank, but it is poorly
conditioned in the small aperture regime. 

For incomplete measurements we can determine the rank of
$\GG_{_{\hspace{-0.02in} S}}(\vy)$ from that of
\begin{align}
\QQ_{_{\hspace{-0.02in} S}}(\vy) = \GG_{_{\hspace{-0.02in}
      S}}(\vy)^\dagger \GG_{_{\hspace{-0.02in} S}}(\vy) = \sum_{r=1}^N
  \bG(\vx_r,\vy)^\dagger \bS \bS^T \bG(\vx_r,\vy).
\label{eq:ra3}
\end{align}
The estimation of the eigenvalues of $\QQ_{_{\hspace{-0.02in}
    S}}(\vy)$ is similar to the above. We obtain that when $a \ll L$, 
and for measurements along a single direction $\ve_1$ or
$\ve_2$, the matrix $\QQ_{_{\hspace{-0.02in}}}(\vy)$ has one large
eigenvalue approximated by $N/(4 \pi L)^2$, a  second
eigenvalue smaller by a factor of order $a^2/L^2$, and an even smaller third
one. When the measurements are made in the longitudinal direction
$\ve_3$, all the eigenvalues are much smaller than $N/(4 \pi
L)^2$. This is the worse measurement setup. Finally, when all the
transversal components of the field are measured, we note that since 
\[
\bS \bS^T = {\bf I}_3 - \ve_3 \ve_3^T, \qquad \bS = (\ve_1,\ve_2),
\]
$\QQ_{_{\hspace{-0.02in}S}}$ is a small perturbation of $\QQ$. It has
almost the same condition number as $\QQ(\vy)$, and this improves as
the ratio $a/L$ grows.

\subsection{Additive noise effects on the SVD}
\label{sect:inv.3}
To unify the discussion for complete and incomplete measurements, let
us consider  the generic problem
\begin{equation}
\tbB = \bB + \bW,
\label{eq:RM0}
\end{equation}
for a low rank $\mathfrak{R}$ matrix $\bB \in \mathbb{C}^{M \times M}$ corrupted
with  additive noise.    In our context $M = 3N$ for complete measurements
and $|\mathscr{S}|N$ for incomplete ones, and $\mathfrak{R}= 3P$. Moreover,
$\tbB$ equals $\tDD$ or $\tDD_{_{\hspace{-0.02in}S}}$, and it is
modeled by \eqref{eq:dm10} and \eqref{eq:dm11S} in the Hadamard data
acquisition scheme described in section \ref{sect:Data.2}.  The noise matrix $\bW$ has mean zero and independent, identically
distributed complex Gaussian entries of standard deviation $\sigma/\sqrt{M}$.

The matrix $\tbB$ has singular values $\ts_j$ and left singular vectors 
$\tbu_j$, for $j = 1, \ldots, M$. Its rank is typically $M$ because of the almost surely  full rank noise matrix, but we are  interested in 
its few singular values that can be distinguished from noise, and the associated singular vectors. 
These are perturbations of the singular values $\sigma_j$ and left singular vectors $\bu_j$ 
of $\bB$, for $j= 1, \ldots, \mathfrak{R}$, and are described below in the asymptotic limit $ M \to \infty$.  

We begin with the asymptotic approximation of the square of the
Frobenius norm
\[
\|\tbB\|_{_F}^2 = \sum_{j,q = 1}^{M} |\widetilde B_{jq}|^2 =
\sum_{j=1}^{M} \ts_j^2,
\]
which satisfies
\cite{marchenko1967distribution,johnstone2001distribution,
  baik2005phase,capitaine2012central}. 
\begin{equation}
M \left( \frac{1}{M} \sum_{j=1}^{M} \ts_j^2 - \sigma^2 \right) \to
\sum_{j=1}^{\mathfrak{R}} \sigma_j^2 +  \sigma^2 \mathcal{Z}_0, \quad
\mbox{as} ~ M \to \infty,
\label{eq:RM1}
\end{equation}
where the convergence is in distribution, and $\mathcal{Z}_0$ follows
a Gaussian distribution with mean zero and variance one. We use this
result in section \ref{sect:inv.alg} to estimate the noise level
$\sigma$ from the measurements. The behavior of the leading singular
values and singular vectors of $\tbB$ described in the next theorem is
used in section \ref{sect:inv.alg} to obtain a robust localization of
the inclusions.

\begin{theorem}
\label{thm}
Let $\bB \in \mathbb{C}^{M \times M}$ be the matrix of fixed rank
$\mathfrak{R}$ in \eqref{eq:RM0}, and $\tbB$ its corrupted version by
the noise matrix $\bW$ with mean zero, independent, identically
distributed complex Gaussian entries of standard deviation $\sigma/\sqrt{M}$.
\\ (i) For $j = 1, \ldots, \mathfrak{R}$, and in the limit $M \to
\infty$, the perturbed singular values satisfy
\begin{equation}
\sqrt{M} \left[ \ts_j - \sigma_j \left(1 + \frac{\sigma^2}{\sigma_j^2}
  \right)\right] \to \frac{\sigma}{\sqrt{2}} \left( 1 -
\frac{\sigma^2}{\sigma_j^2}\right)^{\hspace{-0.02in} 1/2}
\mathcal{Z}_0, \quad \mbox{if} ~ \sigma_j > \sigma,
\label{eq:RM2}
\end{equation}
and 
\begin{equation}
M^{2/3} ( \ts_j - 2 \sigma ) \to \frac{\sigma}{2^{2/3}} \mathcal{Z}_2,
\quad \mbox{if} ~ \sigma_j < \sigma,
\label{eq:RM3}
\end{equation}
where the convergence is in distribution, $\mathcal{Z}_0$ follows a
Gaussian distribution with mean zero and variance one, and
$\mathcal{Z}_2$ follows a type-2 Tracy Widom distribution.  \\ (ii)
Under the additional assumption that the $\mathfrak{R}$ nonzero
singular values of $\bB$ are distinct, and for indices $j$ such that
$\sigma_j > \sigma$, the left singular vectors $\bu_j$ and $\tbu_j$ of
$\bB$ and $\tbB$ satisfy in the limit $M \to \infty$
\begin{equation}
|\tbu_j^\dagger \bu_j|^2 \to 1 - \frac{\sigma^2}{\sigma_j^2}, 
\label{eq:RM4}
\end{equation}
and for $q \ne j \le \mathfrak{R}$, 
\begin{equation} 
|\tbu_q^\dagger \bu_j| \to 0,
\label{eq:RM5}
\end{equation} 
where the convergence is in probability.
\end{theorem}

The type-$2$ Tracy-Widom distribution has the cumulative distribution function $\Phi_{\rm TW2}$ given by
\begin{equation}
\label{eq:cdftw2}
\Phi_{\rm TW2} (z) = \exp \Big( - 
 \int_z^\infty
(x-z) \varphi^2(x) dx \Big)  , 
\end{equation}
with $\varphi(x)$ the solution of the Painlev\'e equation
\begin{equation} 
\label{painv} 
\varphi''(x)=x \varphi(x) +2\varphi(x)^3, \quad 
 \varphi(x)
\simeq {\rm Ai}(x), \, x \to \infty,\end{equation}
 and ${\rm Ai}$ the Airy function.
 The expectation of $\mathcal{Z}_2$  is $ \EE[\mathcal{Z}_2] \simeq -1.771$ and its variance is 
${\rm Var}(\mathcal{Z}_2) \simeq 0.813$.
Details about the Tracy-Widom distributions can be found in \cite{baik08}.

Theorem \ref{thm} was already stated in \cite{garnier2014applications} for the special case $\mathfrak{R}=1$.
The proof of the extension to an arbitrary rank $\mathfrak{R}$ 
can be obtained from the method described in \cite{benaych2011eigenvalues}. 
Note that formula (\ref{eq:RM2}) seems to  predict
that the standard deviation of the perturbed singular value $\ts_j$ cancels
when $\sigma_j \searrow   \sigma $, but this is true only to leading 
order $M^{-1/2}$. In fact the standard deviation  becomes of order $M^{-2/3}$.
Following \cite{baik2005phase},  we can anticipate that there are interpolating distributions which appear when 
$\sigma_j  =    \sigma + {w}{M^{-{1}/{3}}}  $ for some fixed $w$.

\section{Inversion with noisy data}
\label{sect:inv.alg}
We use the singular value decomposition analysis in the previous
section to obtain a robust inversion method for noisy array data.  We
begin in section \ref{sect:noiselev} with the estimation of the noise
level. Then we formulate in section \ref{sect:localize} the method for
localizing the inclusions. The estimation of their reflectivity tensor is
described in section \ref{sect:getref}.

\subsection{Estimation of the noise level}
\label{sect:noiselev}
A detailed analysis of the estimation of $\sigma$, for a variety of
cases, is given in \cite{garnier2014applications}. Here $N \gg
P$ and we assume that  at least one singular value of the noisy
data matrix can be distinguished from the others, so we can image. This simplifies
the estimation of $\sigma$, as we now explain.

Let us begin by rewriting  \eqref{eq:RM1} as
\begin{equation}
\sum_{j=\mathfrak{R} + 1}^M \hspace{-0.05in}\ts_j^2 -
(M-4\mathfrak{R}) \sigma^2 \sim \sum_{j=1}^{\mathfrak{R}}
\left[\sigma_j^2 - \big(\ts_j^2-4 \sigma^2\big) \right] +  
\sigma^2 \mathcal{Z}_0,
\label{eq:IN1}
\end{equation}
where we recall that  $\mathfrak{R} = 3P$, and $M = 3N$ for complete measurements and $|\mathscr{S}|N$
otherwise. The symbol $\sim$ stands for approximate, in the asymptotic
regime $M \gg 1$, and there is no bias in the right hand side in the
absence of the inclusions i.e., when $\sigma_j=0$ and 
$\ts_j \sim 2 \sigma$ for the first indices $j$ by \eqref{eq:RM3}.

The unbiased estimate of the noise level $\sigma$ follows from
\eqref{eq:IN1},
\begin{equation}
\sigma^e = \frac{1}{M - 4 \mathfrak{R}}
\sum_{j=\mathfrak{R}+1}^M \hspace{-0.05in}\ts_j^2,
\label{eq:IN2}
\end{equation}
but it requires prior knowledge of the number of inclusions.  If this is 
unknown, we can use the empirical estimate
\begin{equation}
\sigma^e = \frac{1}{M - 4 \mathfrak{R}^e} \sum_{j=\mathfrak{R}^e+1}^M \ts_j^2,
\label{eq:IN3}
\end{equation}
where $\mathfrak{R}^e$ is the number of singular values of $\tbB$ that
are significantly larger than the others. The estimate \eqref{eq:IN3}
is very close to \eqref{eq:IN2} because $M \gg \mathfrak{R} \ge \mathfrak{R}^e$.

\subsection{Localization of the inclusions}
\label{sect:localize}
We assume henceforth that the noiseless matrix $\bB$, equal to
$\DD$ for complete measurements and $\DD_{_{\hspace{-0.02in} S}}$ for
incomplete ones, has $\mathfrak{R} = 3 P$ distinct singular values
indexed in decreasing order as $\sigma_1, \ldots,
\sigma_{_{\mathfrak{R}}}$. The corresponding left singular vectors are
the columns $\bu_j$ of the $M \times \mathfrak{R}$ matrix 
\begin{equation}
\UU = \left(\bu_1, \ldots, \bu_{\mathfrak{R}}\right).
\end{equation}
The singular values of the noisy matrix $\tbB$ are 
denoted by $\ts_j$, for 
$j = 1, \ldots, M$. We are  interested 
in the first $\widetilde{\mathfrak{R}}$ of them, which can be distinguished from the noise. 

The effective rank 
 $\widetilde{\mathfrak{R}}$  is determined using Theorem \ref{thm},
 \begin{equation}
\widetilde{\mathfrak{R}} = \max\{j =1, \ldots, M ~ ~ \mbox{such that}
~ ~ \ts_j >   \sigma^e r_\theta \}.
\label{eq:IN4}
\end{equation}
where the threshold $r_\theta$ is defined 
\begin{equation}
\label{ralpha}
r_\theta =  2+    \frac{1}{(2M)^{\frac{2}{3}}}   \Phi_{{\rm TW}2} ^{-1} (1-\theta),
\end{equation}
with $\Phi_{{\rm TW}2} $ the cumulative distribution function (\ref{eq:cdftw2}) of
the  Tracy-Widom distribution of type~$2$.  
As shown in \cite{garnier2014applications,ammarigarnier2015}, by the
Neyman-Pearson lemma, the decision rule (\ref{eq:IN4}) maximizes
 the probability of detection  for a given false alarm rate $\theta$. This is a user defined number
 satisfying $0 < \theta \ll 1$. In our case $M$ is large, so for all $\theta$ we  collect the singular values that are  significantly larger than $2 \sigma^e$.

We let
\begin{equation}
\widetilde{\UU} = \left(\tbu_1, \ldots, \tbu_{_{\widetilde{\mathfrak{R}}}} \right)
\label{eq:IN5}
\end{equation}
be the matrix of the leading left singular vectors of $\tbB$.
The MUSIC method \cite{ammari2007music} determines the locations of the inclusions by
projecting $\GG(\vy)\ve$ or $\GG_{_{\hspace{-0.02 in}S}}(\vy) \ve$ on
the range of $\widetilde{\UU}$, for some vector $\ve$ and search points
$\vy$.  With $\ve = \ve_1$ we obtain the MUSIC imaging
function
\begin{equation}
\mathcal{I}_{_{\rm MUSIC}}(\vy) = \left\| \left({\bf I}_{M} - \tUU
\tUU^\dagger\right)\GG(\vy)\ve_1\right\|_{_{F}}^{-1},
\label{eq:IN6}
\end{equation}
in the case of complete measurements, where $M = 3N$. It is the same
for incomplete measurements, except that $\GG(\vy)$ is replaced by
$\GG_{_{\hspace{-0.02 in}S}}(\vy)$ and $M = |\mathscr{S}|N$.  

If noise were negligible, \eqref{eq:IN6} would peak at the locations
of the inclusions due to \eqref{eq:in.8}.  As the noise level 
grows the effective rank $\widetilde{\mathfrak{R}}$ decreases, the
subspace spanned by the columns of $\tUU$ changes, and the MUSIC
images deteriorate. We show next how to improve the 
localization method using the results in Theorem \ref{thm}.

\subsubsection{Localization of one inclusion or well separated inclusions}
\label{sect:localize.1}
The robust localization  is based on 
\eqref{eq:in8} and its equivalent for incomplete measurements, and
Theorem \ref{thm}. It consists of the following steps:

\vspace{0.1in} \noindent \textbf{Step 1.} Estimate the
$\widetilde{\mathfrak{R}}$ singular values of the unknown unperturbed
matrix $\bB$ from equation \eqref{eq:RM2}, rewritten as
\[
\ts_j \approx \sigma_j \left(1 + \frac{\sigma^2}{\sigma_j^2}\right).
\]
Replacing $\sigma$ by its estimate $\sigma^e$, and choosing the root 
that is larger than $\sigma^e$, we obtain 
\begin{equation}
\sigma_j^e = \frac{1}{2} \left[ \ts_j + \sqrt{\ts_j^2 - (2
    \sigma^e)^2}\right], \quad j = 1, \ldots,
\widetilde{\mathfrak{R}}.
\label{eq:IN7}
\end{equation}

\vspace{0.05in} \noindent \textbf{Step 2.} Estimate the changes of
direction of the leading singular vectors, using \eqref{eq:RM4},
\begin{equation}
|\tbu_j^\dagger \bu_j|^2 \approx \cos^2 \theta_j^e = 1 -
\left(\frac{\sigma^e}{\sigma_j^e}\right)^2, \quad j = 1, \ldots,
\widetilde{\mathfrak{R}}.
\label{eq:IN8}
\end{equation}
We also have from \eqref{eq:RM5} that 
\begin{equation}
|\tbu_q^\dagger \bu_j|^2 \approx 0, \qquad j,q = 1, \ldots,
\widetilde{\mathfrak{R}}, ~ ~ j \ne q.
\label{eq:IN9}
\end{equation}

\vspace{0.05in} \noindent \textbf{Step 3.} 
For the search point $\vy$
calculate the matrix $\GG(\vy)$ given by \eqref{eq:in2} for complete
measurements, or $\GG_{_{\hspace{-0.02 in}S}}(\vy)$ given by
\eqref{eq:in2S} for incomplete measurements. Determine the matrix
$\HH(\vy)$ of their left singular vectors, and calculate
\begin{equation}
\bT(\vy) = \widetilde \UU^\dagger \HH(\vy) \HH(\vy)^\dagger \widetilde\UU.
\label{eq:IN10}
\end{equation}
When $\vy = \vy_1$, we conclude from \eqref{eq:in8} that 
\[
\bT(\vy_1) = \widetilde \UU^\dagger \UU \UU^\dagger \widetilde \UU.
\]
Here $\UU$ is unknown, but 
\eqref{eq:IN8}--\eqref{eq:IN9} give that the components of
$\bT(\vy_1)$ satisfy
\begin{equation}
T_{jq}(\vy_1) \approx \delta_{jq} \cos^2 \theta_j^e, \qquad 
j,q = 1, \ldots, \widetilde{\mathfrak{R}}.
\label{eq:IN12}
\end{equation}

\vspace{0.05in} \noindent \textbf{Step 4.} Calculate the imaging
function 
\begin{equation}
\mathcal{I}(\vy) = \left\{\sum_{j,q=1}^{\widetilde{\mathfrak{R}}}
\gamma_j^2 \left[ T_{jq}(\vy)- \delta_{jq} \cos^2 \theta_j^e \right]^2
\right\}^{-1/2},
\label{eq:IN13}
\end{equation}
where $\gamma_j$ are some positive weighting coefficients. We use them
to emphasize the contributions of the terms for large singular values,
and give less weight to those for singular values close to the estimated
noise level $\sigma^e.$ In the numerical simulations 
\begin{equation}
\gamma_j = \mbox{min} \left\{1, \frac{3 (\sigma_j^e -
  \sigma^e)}{\sigma^e} \right\}.
\label{eq:IN14}
\end{equation}
The estimator of $\vy_1$ is the argument of the maximum of $\mathcal{I}(\vy)$.

\vspace{0.06in} This algorithm may be used for localizing multiple
inclusions that are sufficiently far apart, so that the column spaces
of $\HH(\vy_j)$ and $\HH(\vy_q)$, for $j \ne q$, are approximately
orthogonal. To understand what sufficiently far means, we can analyze
the decay of the inner products of the columns of $\GG(\vy_j)$ and
$\GG(\vy_q)$ with the distance $|\vy_q - \vy_j|$. This is similar to
studying the spatial support of
\begin{equation}
\mathcal{F}(\vy_j,\vy_q) = \sum_{r=1}^N \frac{e^{i k
    \left(|\vx_r-\vy_q|-|\vx_r-\vy_j|\right)}}{(4 \pi)^2 |\vx_r-\vy_j||\vx_r-\vy_q|},
\label{eq:IN15}
\end{equation}
the point spread function of the reverse time (Kirchhoff) migration
method for sonar imaging. The known resolution limits of
this method give that the inclusions are well separated when 
\begin{equation}
|\by_q - \by_j| \gg \frac{\la L}{a}, \qquad |y_{q,3}-y_{j,3}| \gg
\frac{\la L^2}{a^2}.
\label{eq:IN16}
\end{equation}
Here we used the notation $\vy_q = (\by_q,y_{q,3})$ and similar for $\vy_j$.

\subsubsection{Localization of multiple inclusions}
\label{sect:localize.2}
Imaging of nearby inclusions is more difficult because they interact,
and multiple inclusions may be associated with one singular value and
vector. To address this case we modify Steps 3-4 of the imaging method
using the results  \eqref{eq:in8m}--\eqref{eq:in.8}, which imply that
\begin{equation}
\bh_q(\vy) = \sum_{l=1}^{\mathfrak{R}} [\bu_l^\dagger \bh_q(\vy)]
\bu_l \quad \mbox{iff} ~  \vy \in\{\vy_1, \ldots, \vy_P\}.
\label{eq:IN17}
\end{equation}
Here $\bh_q(\vy)$ are the left singular vectors of $\GG(\vy)$, the columns of $\HH(\vy)$ for $q = 1, 2, 3$,
and the right hand side is their projection on the span of $\{\bu_1, \ldots, \bu_{\mathfrak{R}}\}$, the range of the noiseless data matrix.
It suffices to work with $q = 1$, which gives the largest
contribution to the data model in equations
\eqref{eq:in4}--\eqref{eq:in5}.  The analysis in section
\ref{sect:inv.2} shows that the third singular value of $\GG(\vy)$ is
very small in the small aperture regime, and for incomplete measurements
along $\ve_1$ or $\ve_2$ the matrix $\GG_{_{\hspace{-0.02in} S}}(\vy)$
has only one large singular value. Thus, the terms with $\bh_2$ and
$\bh_3$ may give small contributions to $\DD$ or
$\DD_{_{\hspace{-0.02in}S}}$ and consequently, $\bh_2$ and $\bh_3$ may
have small projections on the subspace spanned by $\{\bu_j\}_{1 \le j \le \widetilde{\mathfrak{R}}}$ when $j  \le \widetilde{\mathfrak{R}} <
\mathfrak{R}$. This is why we use $\bh_1$ in the localization.

\vspace{0.05in}
We cannot work with \eqref{eq:IN17} directly, because we do not know
the unperturbed singular vectors. However, we can calculate their projection on the 
$\mbox{span}\{\tbu_l\}$ for $l = 1, \ldots,
\widetilde{\mathfrak{R}}$, which satisfies by 
\eqref{eq:IN9} and \eqref{eq:IN17}
\begin{equation}
\tbu_l[\tbu_l^\dagger \bh_1(\vy)] = \tbu_l \sum_{q=1}^{\mathfrak{R}} (\tbu_l^\dagger \bu_q) [\bu_q^\dagger \bh_1(\vy)]
 \approx  \tbu_l (\tbu_l^\dagger \bu_l) 
[\bu_l^\dagger \bh_1(\vy)], \label{eq:IN18}
\end{equation}
for $\vy \in\{\vy_1, \ldots, \vy_P\}.$ 
Taking the Euclidian norm in \eqref{eq:IN18} and using \eqref{eq:IN8} we obtain
\begin{equation}
| \bu_l^\dagger \bh_1(\vy) | \approx \frac{|\tbu_l^\dagger
  \bh_1(\vy)|}{|\cos \theta_l|}, \quad l = 1, \ldots,
\widetilde{\mathfrak{R}},
\end{equation}
and the norm of \eqref{eq:IN17} gives
\begin{equation}
\label{eq:IN18p}
1 = \|\bh_1(\vy)\|^2 = \sum_{l=1}^{\mathfrak{R}} | \bu_l^\dagger
\bh_1(\vy) |^2, \quad ~ \vy \in\{\vy_1, \ldots, \vy_P\}.
\end{equation}
We have only the first $\widetilde{\mathfrak{R}}$ terms in the right
hand side of \eqref{eq:IN18p}, and we use them to obtain the imaging
function
\begin{equation}
\mathcal{I}(\vy) = \left[ 1 -
  \sum_{j=1}^{\widetilde{\mathfrak{R}}} \frac{|\tbu_j^\dagger
  \bh_1(\vy)|^2}{\cos^2 \theta_j^e} \right]^{-1/2}_.
\label{eq:IN19}
\end{equation}
We expect from \eqref{eq:IN17} that when 
$\widetilde{\mathfrak{R}} = \mathfrak{R}$, the imaging function
\eqref{eq:IN19} peaks at $\vy = \vy_j$, for $j = 1,\ldots, P$.  This
happens only for  weak noise. In general $\widetilde{\mathfrak{R}}
< \mathfrak{R}$, and \eqref{eq:IN19} peaks at the
location of the stronger inclusions, which contribute to the
$\widetilde{\mathfrak{R}}$ distinguishable singular values.

\subsection{Estimation of the reflectivity tensor}
\label{sect:getref}
If we knew the locations
$\{\vy_p\}_{1\le p \le P}$ of the inclusions, we could obtain from
\eqref{eq:dm10} and \eqref{eq:in4} that
\begin{align}
\HH(\vy_p)^\dagger \tDD \, \overline{\HH(\vy_p)} = 
\bR_p + \bDe_{\bR_p}, \quad ~ p = 1, \ldots, P.
\label{eq:IN20}
\end{align}
This is for complete measurements,  and the error 
\begin{equation}
\bDe_{\bR_p} = \HH(\vy_p)^\dagger \bW \overline{\HH(\vy_p)} + \sum_{l
  \ne p, l = 1}^P [ \HH(\vy_p)^\dagger \HH(\vy_l)] \, \bR_l\, [
  \HH(\vy_p)^\dagger \HH(\vy_l)]^T
\label{eq:IN21}
\end{equation}
is due to the noise and interaction of the inclusions. The
interaction  is small when the inclusions are well separated. 

We only have estimates $\vy_p^e$ of $\vy_p$, the peaks of the imaging
function defined in the previous section, so we calculate instead
\begin{align}
\widetilde{\bR}_p = \HH(\vy_p^e)^\dagger \tDD \,
\overline{\HH(\vy_p^e)}, \quad ~ p = 1,
\ldots, P.
\label{eq:IN22}
\end{align}
These are modeled by 
\begin{equation}
\widetilde{\bR}_p = \bR_p + \widetilde{\bDe}_{\bR_p},
\label{eq:IN22p}
\end{equation}
where $\widetilde{\bDe}_{\bR_p}$ includes \eqref{eq:IN21} and 
additional terms due to location errors $\vy_p^e - \vy_p$. The estimate
of $\brho_p$ is motivated by definition \eqref{eq:in5} of $\bR_p$,
\begin{align}
\brho_p^e = \bGamma(\vy_p^e) \widetilde \bR_p \bGamma(\vy_p^e)^T, 
\label{eq:IN23}
\end{align}
where 
\begin{equation}
\label{eq:defGamma}
\bGamma(\vy) = \bV(\vy)  \bSigma(\vy)^{-1}.
\end{equation}
Recall from \eqref{eq:in3} that $\bSigma(\vy)$ is the
$3\times 3$ matrix of singular values of $\GG(\vy)$ and $\bV(\vy)$ the 
$3 \times 3$ matrix of its right singular vectors. 

We can model  \eqref{eq:IN23}  by 
\begin{equation}
\brho_p^e = \brho_p + \bDe_{\brho_p},
\label{eq:IN24}
\end{equation}
with error $\bDe_{\brho_p}$ that consists of two parts. The
first part depends on the noise, the interaction of the inclusions and
the error in their estimated locations
\begin{align}
\bDe^{(1)}_{\brho_p} = \bGamma(\vy_p^e)\, \widetilde{\bDe}_{\bR_p}
\bGamma(\vy_p^e)^T.
\end{align}
The second part vanishes when the inclusion is correctly localized
\begin{align}
\bDe^{(2)}_{\brho_p} = \bGamma(\vy_p^e) \bR_p \bGamma(\vy_p^e)^T
-\bGamma(\vy_p) \bR_p \bGamma(\vy_p)^T.
\end{align}

\subsubsection{Discussion}
\label{sect:disc} The numerical results in section \ref{sect:num} 
show that the inclusion localization in cross-range is excellent for
the case of complete measurements. However, the range localization
deteriorates quickly with noise in the small aperture regime.  This can be
understood from the analysis in section \ref{sect:inv.2}, equation
\eqref{eq:ra2} in particular, which shows that $\QQ(\vy) =
\GG(\vy)^\dagger \GG(\vy)$ has two leading orthonormal eigenvectors
${\itbf v}_q$ which are approximately in $\mbox{span}\{\ve_1,
\ve_2\}$, for $ q = 1,2$. They correspond to eigenvalues that are
approximately $N/(4 \pi L^2)$. The third
eigenvector is ${\itbf v}_3 \approx \ve_3$, for a much smaller
eigenvalue by a factor of the order of $a^2/L^2$. These eigenvalues are the
squares of the entries in $\bSigma(\vy)$, and the eigenvectors are the
columns of $\bV(\vy)$. 

We also see from  \eqref{eq:in2} and
\eqref{eq:ra1p} that the  left singular vectors of $\GG(\vy)$
are
\begin{equation}
\bh_q(\vy) \approx N^{-1/2} \begin{pmatrix} 
e^{i k |\vx_1-\vy|} {\itbf v}_q \\ \vdots \\
e^{i k |\vx_N-\vy|} {\itbf v}_q \end{pmatrix}, \qquad q = 1,2, \quad a \ll L,
\end{equation}
with phases
\begin{equation}
\label{eq:phase}
k |\vx_r - \vy| = k L \left( 1 + \frac{y_3-L}{L} +
\frac{|\bx_r-\by|^2}{2 L^2}\right) + O \left(\frac{a^3}{\la L}\right).
\end{equation}
These vectors define the projection matrices $\bh_q(\vy)
\bh_q(\vy)^\dagger$ used in the localization method, for $q = 1,2$. The singular vector  
$\bh_3(\vy)$ is unlikely to play a role when noise is present. Moreover, the
projection matrices are approximately independent of the range
component $y_3$ of $\vy = (\by,y_3)$, due to cancellation in the 
difference of the phase \eqref{eq:phase}. Thus, range estimation becomes impossible in the small aperture regime.

Both terms in the error $\bDe_{\brho_p}$ in the estimation of the reflectivity tensor
\eqref{eq:IN23} have large components when the diagonal matrix
$\bSigma = \mbox{diag}(\Sigma_1,\Sigma_2,\Sigma_3)$ has a small
entry. This is because we take its inverse in
\eqref{eq:defGamma}. Since in the small aperture regime
\[
\Sigma_3(\vy) \ll \Sigma_{1}(\vy) \approx \Sigma_2(\vy),
\]
and the third column of $\bV(\vy)$ is approximately $\ve_3$, we expect
$\left(\bDe_{\brho_p}\right)_{ij}$ to be large for $i$ or $j$ equal to
$3$. The other components of the error should be of the order of the
noise and inclusion interaction. This is indeed demonstrated by the
numerical simulations in section \ref{sect:num}. The condition number
of $\bSigma(\vy)$ improves for larger ratios $a/L$, and so do the
estimates of $\brho_p$.

\subsubsection{Incomplete measurements}
The estimation of the reflectivity tensors from incomplete measurements can be
done formally as above, with $\HH(\vy)$ the matrix of left singular
vectors of $\GG_{_{\hspace{-0.02in} S}}(\vy)$ and $\bSigma(\vy)$ the
matrix of its singular values. The results are expected to be much
worse when only one component of the electric field is measured, that is for 
$\bS = (\ve_q)$ and some $q \in \{ 1, 2 , 3\}$. As shown in section
\ref{sect:inv.2}, matrix $\bSigma(\vy)$ has at most one large entry in this case,
so the effective rank $\widetilde{\mathfrak{R}}$ in the presence of
noise is smaller than in the complete measurement case. This
translates into worse localization of the inclusions. Moreover, more
components of the errors $\bDe_{\brho_p}$ are large, due to the inversion
of $\bSigma(\vy)$ in $\eqref{eq:defGamma}$.

\section{Numerical results}
\label{sect:num}
In this section we present numerical results that assess the
performance of the imaging method described in section
\ref{sect:inv.alg}.  We begin with the setup of the simulations in
section \ref{sect:num.1}. Then we compare in section \ref{sect:num.0}
the asymptotic limits stated in Theorem \ref{thm} to the empirical statistics of the singular values and vectors obtained with
Monte Carlo simulations. The inversion results are in section
\ref{sect:num.2} for one inclusion and in section \ref{sect:num.3} for
multiple inclusions.
\subsection{Description of the numerical simulations}
\label{sect:num.1}

We present results for two scattering regimes. The first is  called a 
\emph{large aperture regime},   because the array aperture is as large as the range of the inclusions $L = a = 10
\la$. The second is a \emph{small aperture regime} with $ L = 100 \la \gg a
= 10 \la$. The sensors are placed on a regular square grid in the
array aperture, with spacing $\la/2$, so $N = 441$. The three-dimensional search domain
is sampled in steps of $\la/2$.  

We model the scattered field using equation \eqref{eq:dm1}, for
inclusions with reflectivity tensor \eqref{eq:dm2}. Since the equation is
linear, we factor out the small scaling factor $\alpha^3$.  The
inclusions are ellipsoids with scaled semiaxes $a_{p,j}$, for $p =1,
\ldots, P$ and $ j = 1,2,3$.  Their scaled volume is \[|\Omega_p| =
\frac{4 \pi}{3} \prod_{q=1}^3 a_{p,q}.\] In a system of coordinates
with axes of the ellipsoid, the polarization tensor in \eqref{eq:dm2}
is diagonal. We use the system of coordinates centered at the array,
with basis $\{\ve_1,\ve_2,\ve_3\}$, which is a rotation of that of the
ellipsoids by some matrix $\boldsymbol{\mathscr{R}}_p$. The polarization tensors are
\cite{ammari2007music}
\begin{equation}
{\bf M}_p = |\Omega_p| \, \boldsymbol{\mathscr{R}}_p
\,\mbox{diag} \left( \frac{1}{1 + (\ep_p/\ep_o -1) \mathscr{D}_{p,q}}, ~
q = 1, 2,3 \right) \boldsymbol{\mathscr{R}}_p^T,
\end{equation}
where $\mathscr{D}_{p,q}$ are the depolarization factors of the
ellipsoids \cite[Section
  3.3]{sihvola1992polarizability}, given by the elliptic integrals 
\begin{equation}
\mathscr{D}_{p,q} = \frac{|\Omega_p|}{2} \int_0^\infty \frac{ds}{(s +
  a_{p,q}^2) \big[\prod_{l=1}^3 (s+a_{p,l})^2\big]^{1/2}}.
\end{equation}
For example, in the case of one inclusion at $\vy_1 =
(1,-1,L)$, we take $a_{1,q} = q$, for $q = 1, 2, 3$, the contrast
$\ep_1/\ep_{o} = 10$, and some rotation $\boldsymbol{\mathscr{R}}_1$ 
to  obtain the reflectivity tensor
\begin{equation}
\alpha^{-3} \brho_1 = \left(\begin{matrix} 55.4 & -7.28 & -13.43 \\
-7.28 & 70.82 & -22.64 \\
-13.43 & -22.64 & 70.75
\end{matrix} \right).
\label{eq:rhoT1}
\end{equation}

The noise matrix $\bW$ is generated
with the MATLAB command \emph{randn}, as in
\begin{equation}
\label{eq:noiseW}
\bW = \frac{\sigma}{\sqrt{2 M}}\left(\mbox{randn}(M) + i \,
\mbox{randn}(M)\right),
\end{equation}
where $M = 3N$ for complete measurements and $|\mathscr{S}|N$ for
incomplete ones. The noise level $\sigma$ is chosen as a percentage of
the largest singular vale $\sigma_1$ of the noiseless data
matrix $\DD$. Thus, when we say $50\%$ noise, we mean that $\sigma = 0.5
\sigma_1$.

\vspace{0.05in}
\textbf{Remark:} 
The typical amplitude of the  entries of the unperturbed matrix 
$\DD$ is $\sigma_1/M$, while the  noise entries in $\bW$ 
 (with the Hadamard acquisition scheme) have standard deviation $\sigma/\sqrt{M}$. Therefore, the  measured  $\widetilde \DD =\DD+\bW$
has the signal-to-noise ratio ${\rm SNR} = \sigma_1/(\sigma \sqrt{M})$.  In the simulations $\sigma/\sigma_1$ varies  between 
$10\%$ and $75\%$. This is very strong noise. For example, when $\sigma/\sigma_1 = 50\%$, the SNR is $2/\sqrt{M} = 0.05$ in the complete measurement case, where $M = 3 N = 1323$. This small SNR means that 
the signal is very weak and buried in noise.

\subsection{Statistics of the singular values and singular vectors}
\label{sect:num.0}
We compare here the asymptotic behavior of the singular values and
singular vectors stated in Theorem \ref{thm} to the empirical
estimates of their statistics obtained with Monte Carlo
simulations. We take a single inclusion with reflectivity tensor
\eqref{eq:rhoT1} and one thousand samples of the noise matrix
\eqref{eq:noiseW}.  

\subsubsection{Complete measurements}
Because there is a single inclusion, the data matrix $\bB = \DD$ has rank 
three. In the large aperture regime its singular values are 
\[
\sigma_1 = 1.703, ~ ~ \sigma_2 = 1.126, ~ ~ \sigma_3 = 0.183,
\]
and in the small aperture  regime they are 
\[
\sigma_1 = 0.021, ~ ~ \sigma_2 = 0.015, ~ ~ \sigma_3 = 3.01 \cdot 10^{-5}.
\]
The difference in the magnitudes of $\sigma_1$ and $\sigma_2$ in the
two regimes is due to the geometrical spreading factor of order $1/(4
\pi L)^2$.  As expected from the discussion in section
\ref{sect:inv.2}, we have $\sigma_3 \ll \sigma_2 \lesssim \sigma_1$ in
the small aperture regime.

\begin{figure}[t]
\begin{center}
\includegraphics[width=0.34\textwidth]{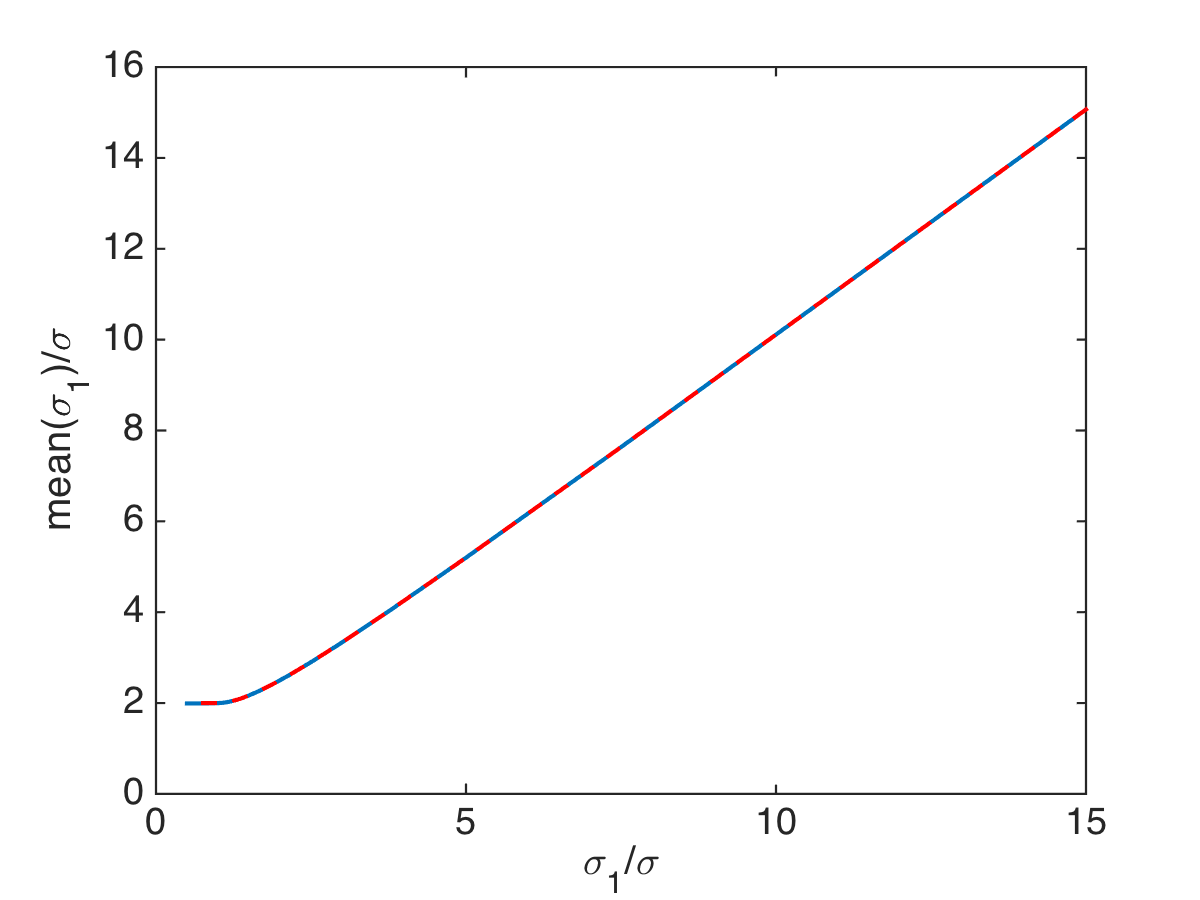}
\hspace{-0.15in}\includegraphics[width=0.34\textwidth]{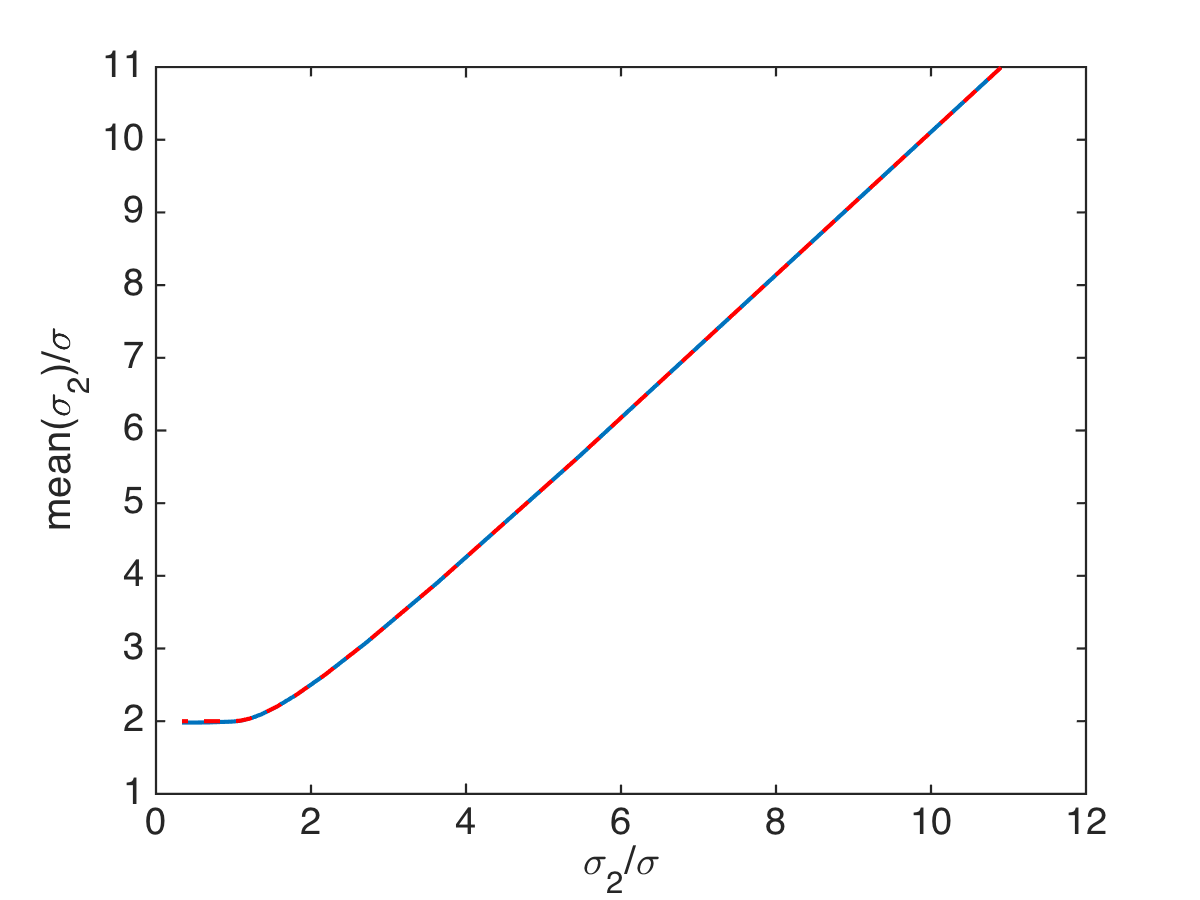}
\hspace{-0.15in}\includegraphics[width=0.34\textwidth]{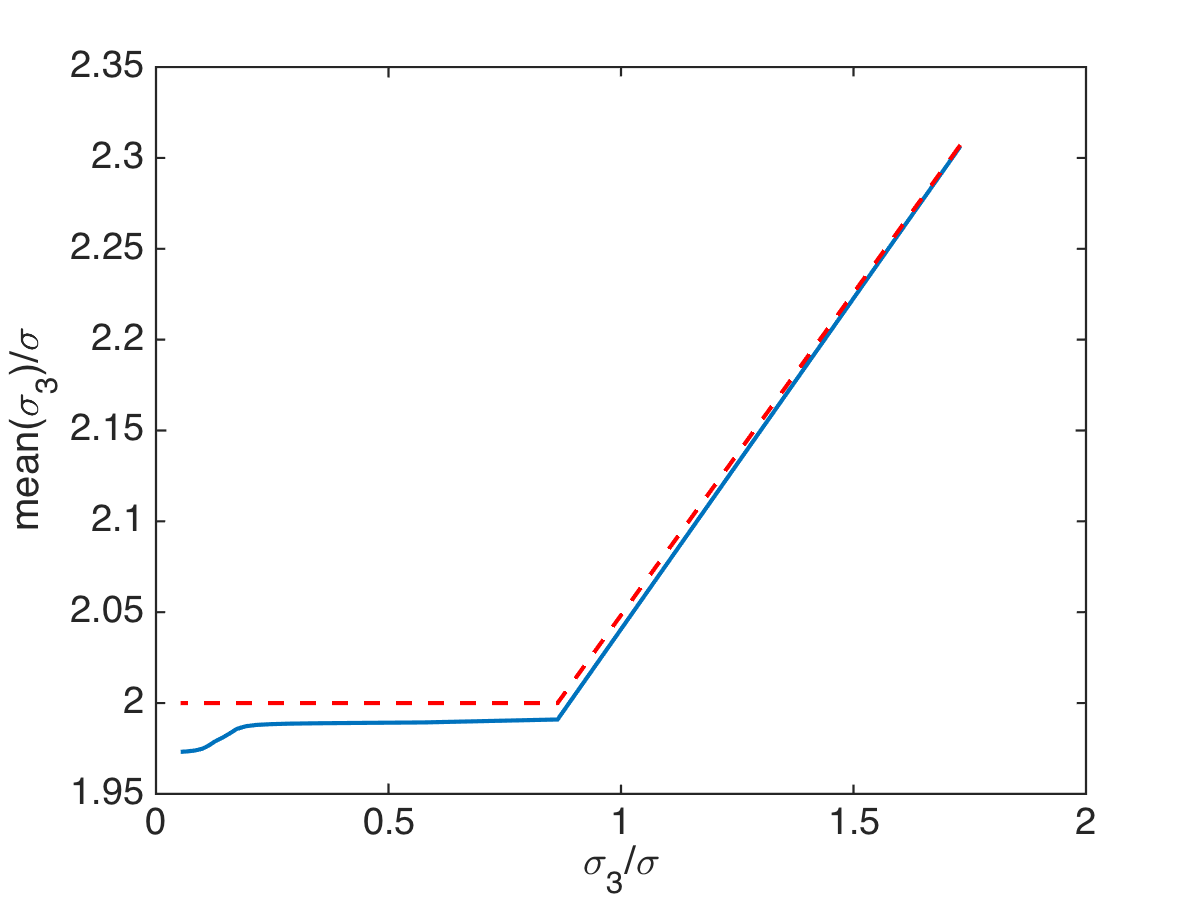}\\
\includegraphics[width=0.34\textwidth]{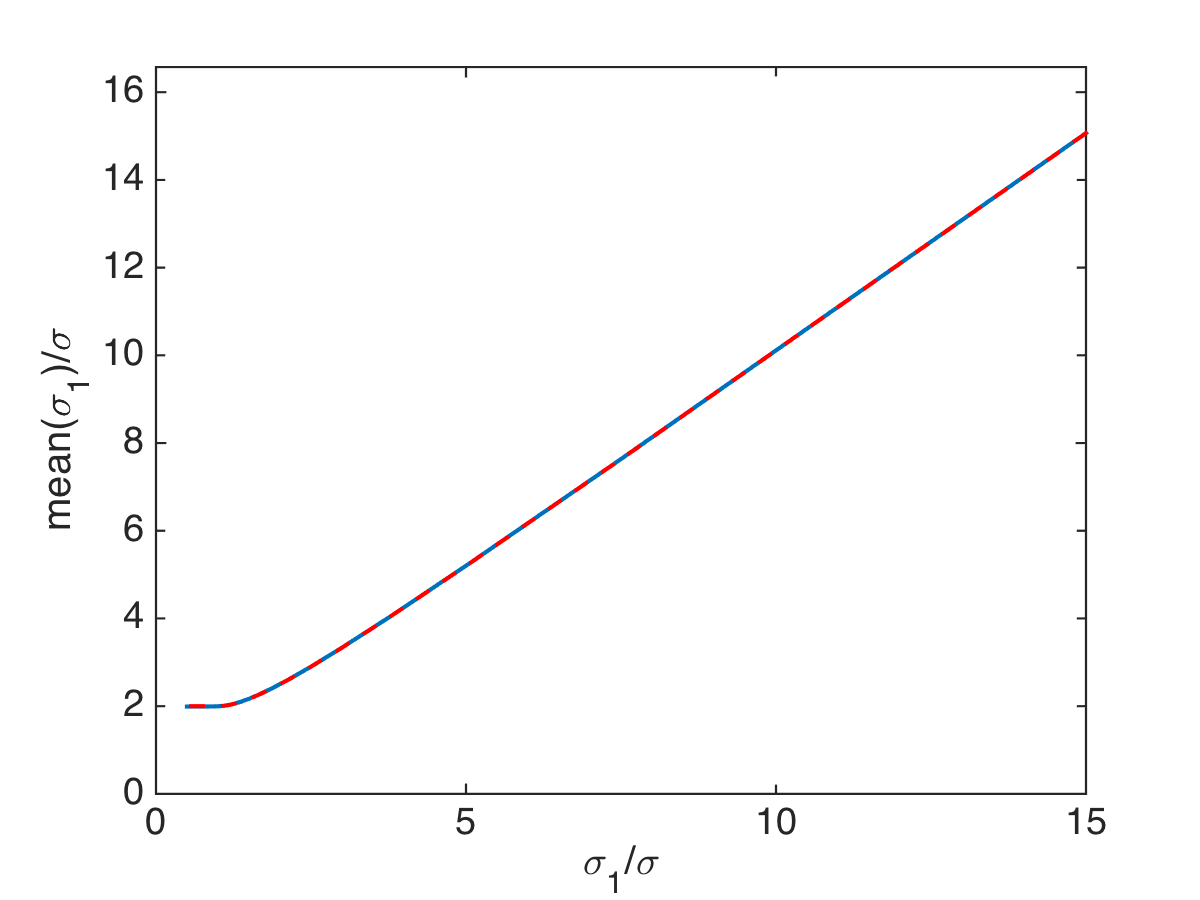}
\hspace{-0.15in}\includegraphics[width=0.34\textwidth]{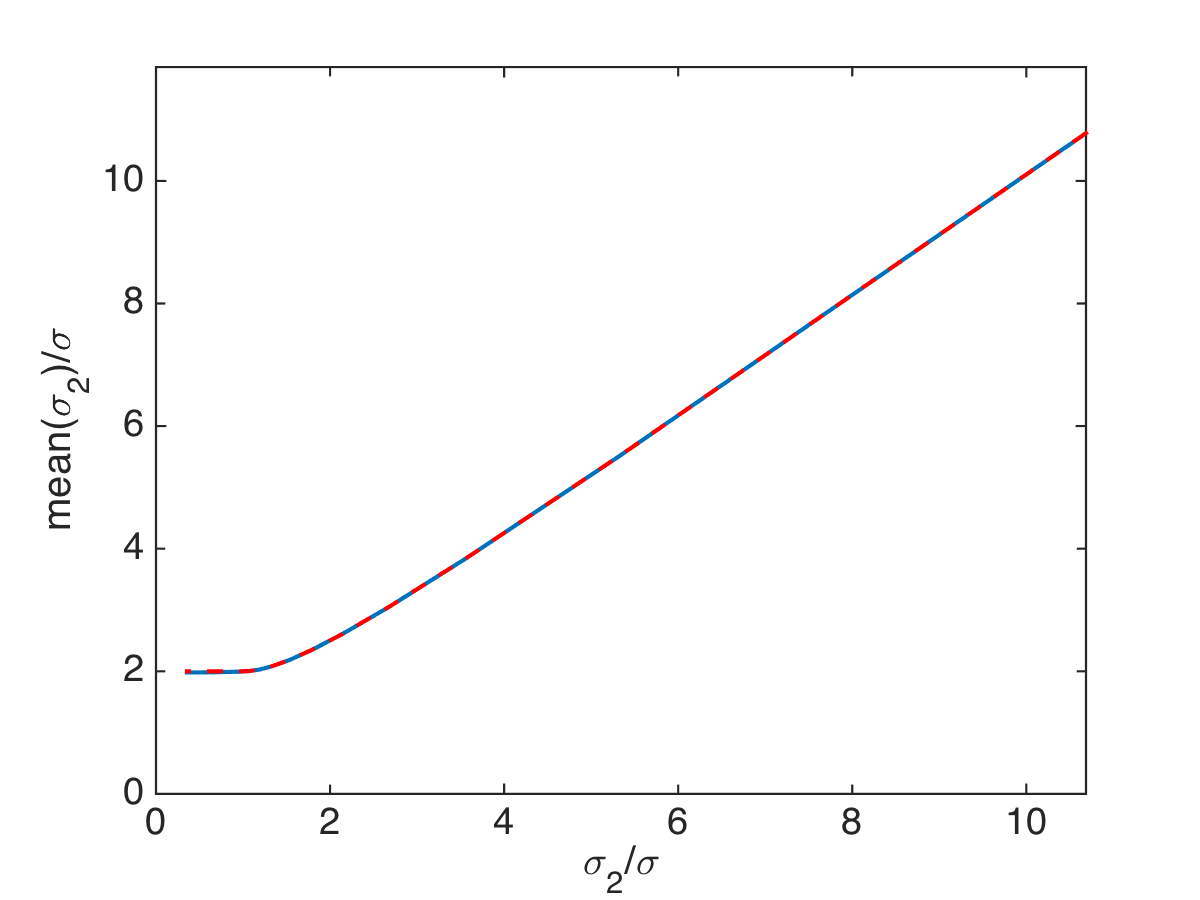}
\hspace{-0.15in}\includegraphics[width=0.34\textwidth]{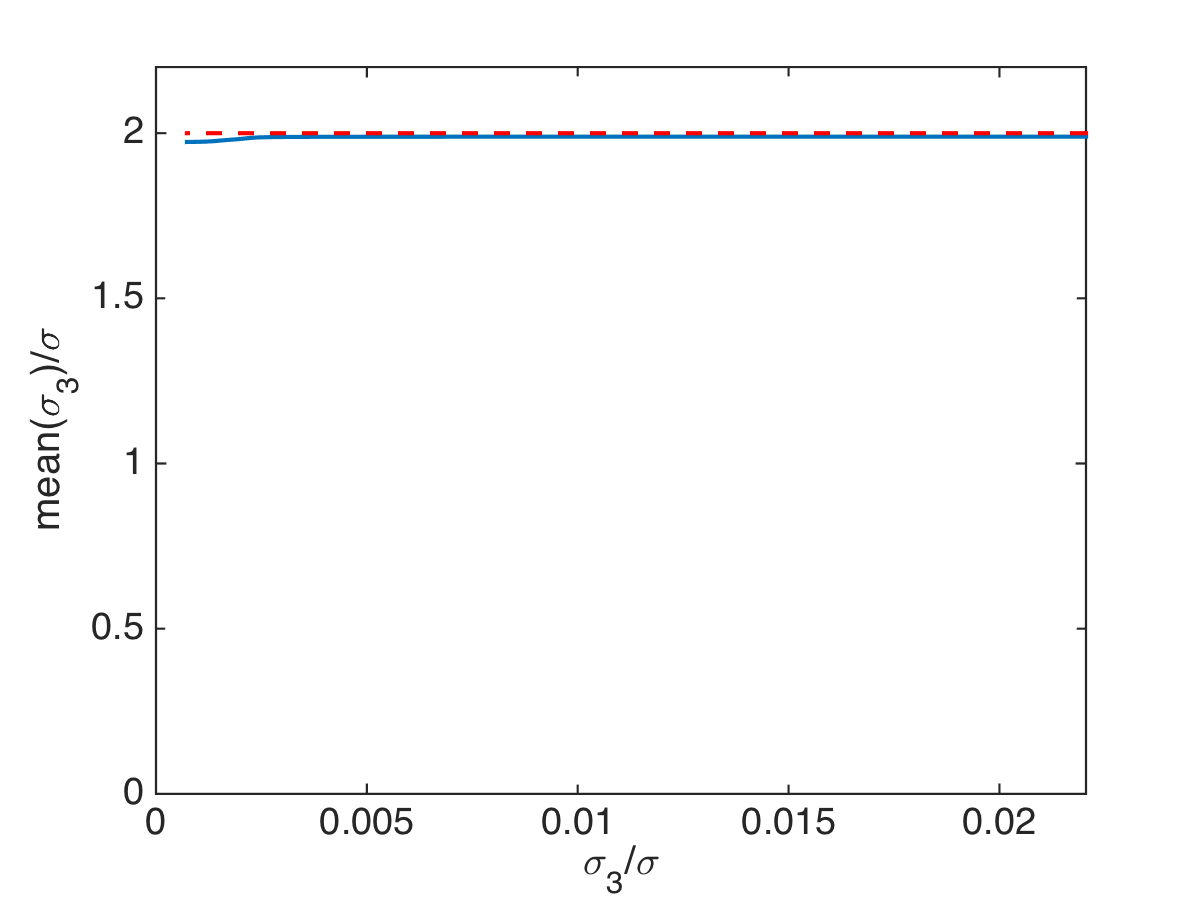}
\end{center}
\vspace{-0.1in}
\caption{Plots of $\mbox{mean}[\ts_j]/\sigma$ vs. $\sigma_j/\sigma$,
  for $ j = 1, 2,3$. The empirical estimates of the means are shown with
  solid blue lines and the asymptotic estimates \eqref{eq:SVDTh} are shown
  with dotted red lines. The top row is for the large aperture regime and
  the bottom row for the small aperture regime.}
\label{fig:SVD1}
\end{figure}

We plot with  solid blue lines in Figure \ref{fig:SVD1} the empirical means
of the singular values $\ts_j$ of $\tDD$, for $j = 1,2,3$ and $\sigma
\in (0, 2 \sigma_1)$. The abscissa is $\sigma_j/\sigma$ and the
ordinate is $\mbox{mean}[\ts_j]/\sigma$. The results are in excellent
agreement with the asymptotic formulas
\begin{equation}
\ts_j \approx \left\{ \begin{array}{ll} \sigma_j \left(1 +
  \sigma^2/\sigma_j^2\right) \quad & \mbox{if} ~ \sigma_j > \sigma,
  \\ 2 \sigma & \mbox{if} ~ \sigma_j < \sigma,
\end{array}
\right. 
\label{eq:SVDTh}
\end{equation}
plotted in the figure with dotted red lines. 
The top row in Figure~\ref{fig:SVD1} is for the large aperture regime and the bottom row for the
small aperture  regime. The plots are similar, except that in the small aperture case
 $\sigma_3$ cannot be distinguished from noise. 
\begin{figure}[h]
\begin{center}
\includegraphics[width=0.34\textwidth]{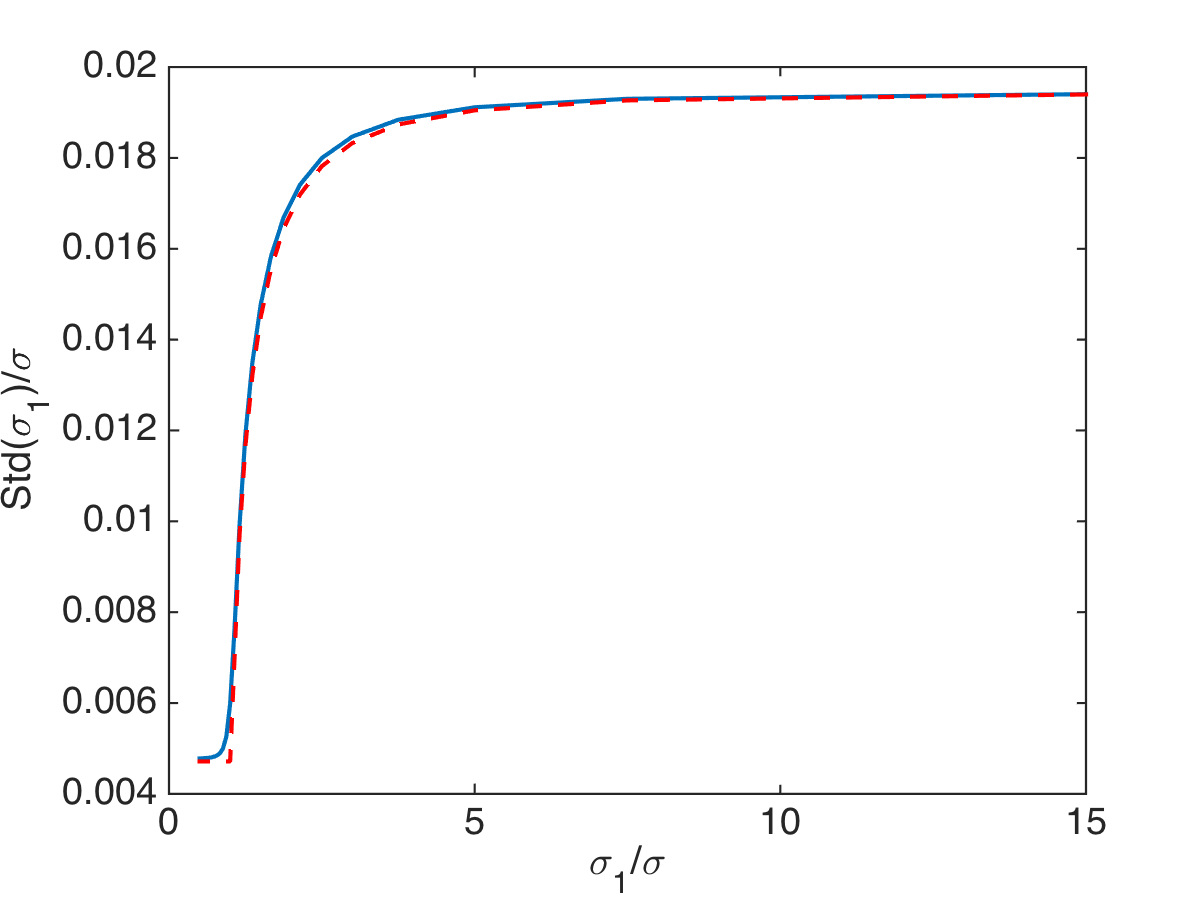}
\hspace{-0.15in}\includegraphics[width=0.34\textwidth]{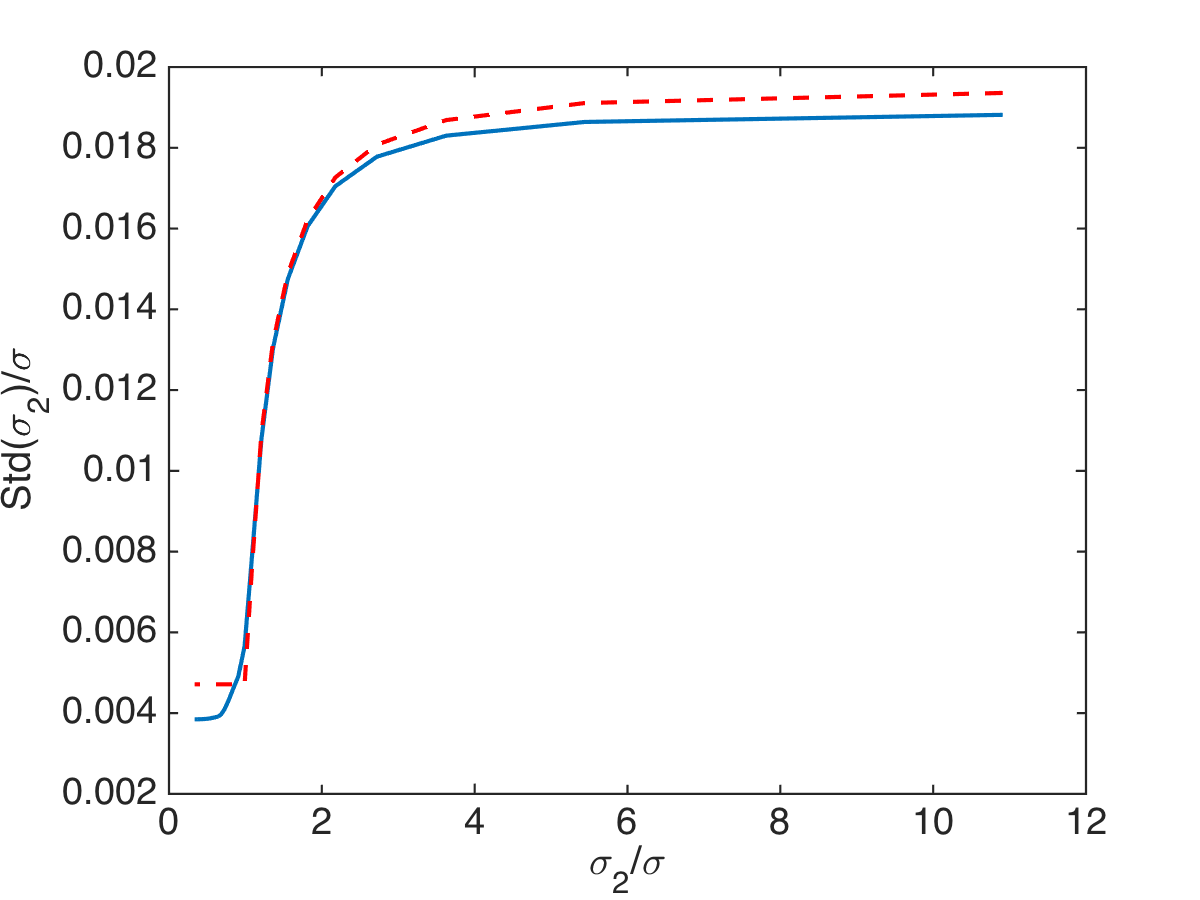}
\hspace{-0.15in}\includegraphics[width=0.34\textwidth]{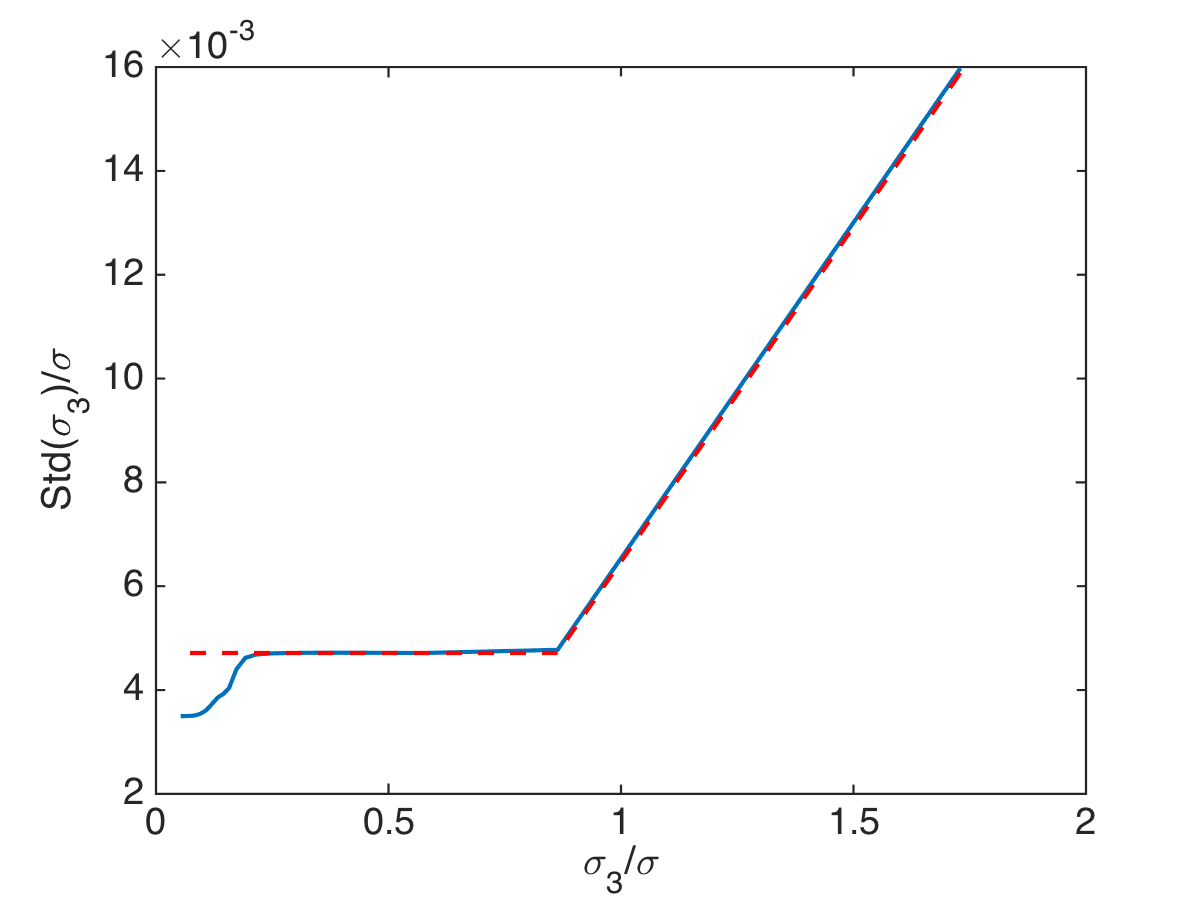}
\end{center}
\vspace{-0.1in}
\caption{Plots of $\mbox{std}[\ts_j]/\sigma$ vs. $\sigma_j/\sigma$,
  for $ j = 1, 2,3$ in the large aperture regime. The empirical estimates are
  shown with solid blue lines and the asymptotic estimates
  \eqref{eq:SVDTh1} with  dotted red lines. }
\label{fig:SVD2}
\end{figure}

In Figure~\ref{fig:SVD2} we plot the standard deviations of $\ts_j$
and note that they are much smaller than the means, 
so $\ts_j$ are approximately deterministic. 
We show the plots of $\mbox{std}(\ts_j)$ in the large aperture regime,
and compare them with the asymptotic ones in Theorem~\ref{thm},
\begin{equation}
\mbox{std}(\ts_j) \approx \left\{ \begin{array}{ll} \frac{\sigma}{(6 N)^{1/2}}
  (1-\sigma^2/\sigma_j^2)^{1/2} \quad &\mbox{if} ~ \sigma_j > \sigma,
  \\ \frac{\sigma}{(6 N)^{2/3}}  \sqrt{0.813} & \mbox{if}
  ~ \sigma_j < \sigma.
\end{array} \right.
\label{eq:SVDTh1}
\end{equation}
Here we used that $\mbox{std}(\mathcal{Z}_2) \approx \sqrt{0.813}.$
The results at small aperture are similar.

\begin{figure}[t]
\begin{center}
\includegraphics[width=0.34\textwidth]{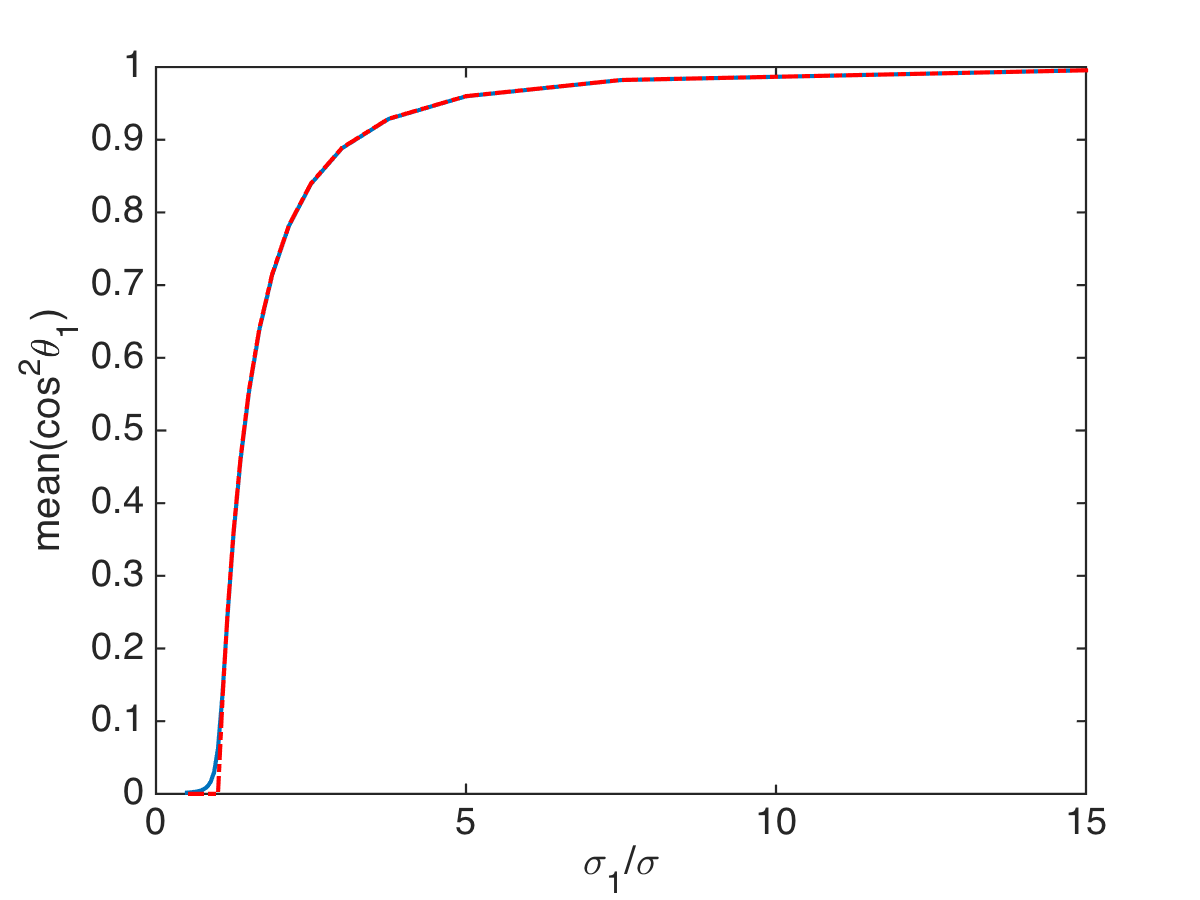}
\includegraphics[width=0.34\textwidth]{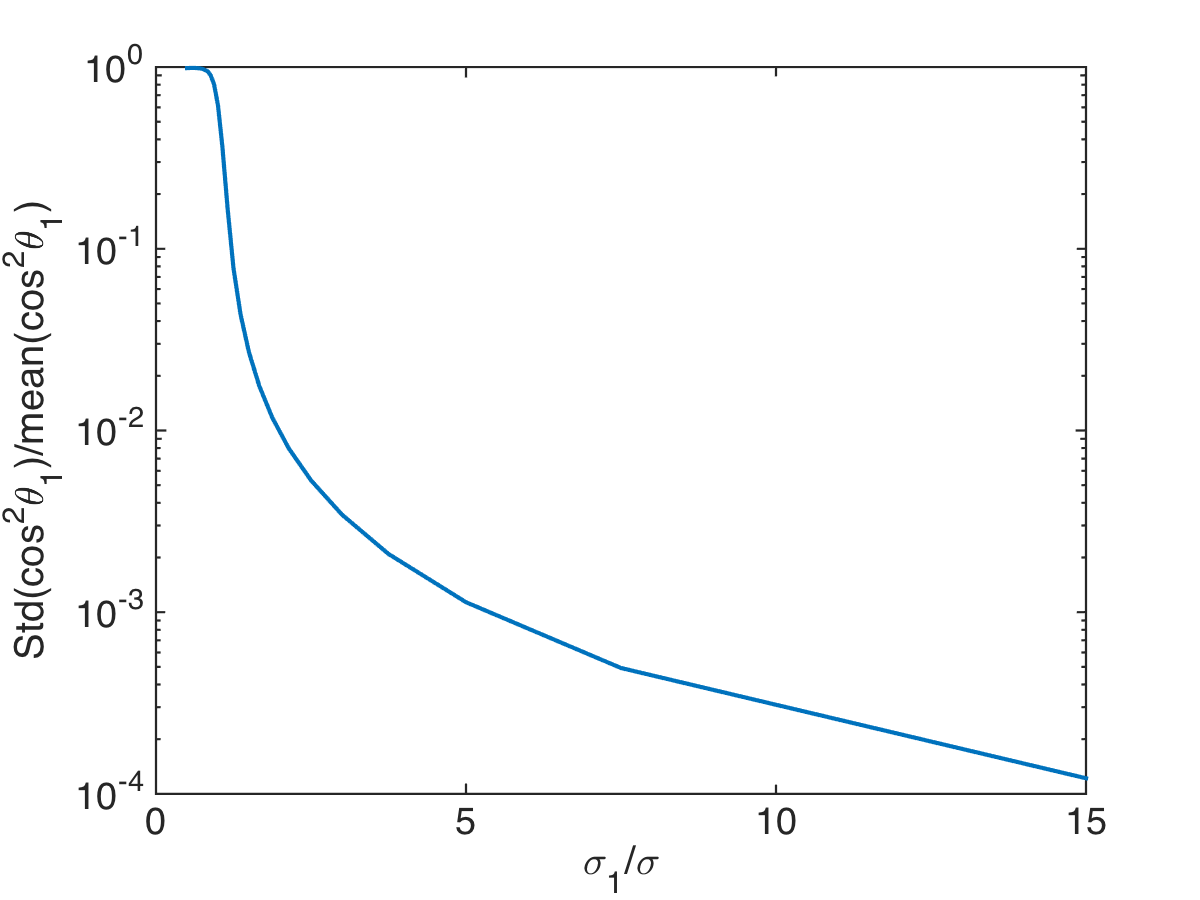}
\end{center}
\vspace{-0.1in}
\caption{Plots of $\mbox{mean}[\cos^2 \theta_1]$ (left) and
  $\mbox{std}[\cos^2 \theta_1]/\mbox{mean}[\cos^2 \theta_1]$ (right)
  vs. $\sigma_1/\sigma$, in the large aperture regime. The empirical estimates
  are shown with solid blue lines and the asymptotic estimates
  \eqref{eq:SVDTh1} are plotted with dotted red lines. }
\label{fig:SVD3}
\end{figure}

In Figure \ref{fig:SVD3} we display the empirical mean of $\cos^2
\theta_1 = |\tbu_1 \bu_1|^2$ and compare it with its asymptotic
prediction \eqref{eq:RM4}.  We also plot the standard deviation, and
note that it is much smaller than the mean when $\sigma_1 > \sigma$,
as expected from Theorem \ref{thm}. The plots look the same for the
large and small aperture regimes, and for the other angles $\theta_2$ and
$\theta_3$, when $\sigma_2$ and $\sigma_3$ are greater than $\sigma$.

\subsubsection{Incomplete measurements} 
We do not include plots for the incomplete measurements, where $\bB = \DD_{_{\hspace{-0.02in}S}}$,  because they
do not add much information. Instead, we display below the effect of
the sensing matrix $\bS$ on the singular values of
$\DD_{_{\hspace{-0.02in}S}}$.  As explained in section \ref{sect:inv.alg}, we  expect that the inversion algorithm 
will be most effective when the effective rank
$\widetilde{\mathfrak{R}}$ equals the rank $\mathfrak{R} = 3$ of $\DD_{_{\hspace{-0.02in}S}}$. 
Since $\widetilde{\mathfrak{R}}$ is  the number of singular 
values of $\DD_{_{\hspace{-0.02in}S}}$ that are larger than $2 \sigma$, the more small 
singular values $\DD_{_{\hspace{-0.02in}S}}$ has, the worse the inversion results.

 When $\bS = (\ve_1)$, the singular values of
$\tDD_{_{\hspace{-0.02in} S}}$ are
\[
\sigma_1 = 1.099, ~ ~ \sigma_2 = 0.087, ~ ~ \sigma_3 = 0.007
\]
in the large aperture regime and 
\[
\sigma_1 = 0.015, ~ ~ \sigma_2 = 1.7 \cdot 10^{-5}, ~ ~ \sigma_3 = 1.6
\cdot 10^{-8}
\]
in the small aperture regime. The results are similar for $\bS = (\ve_2)$,
and show that as predicted by the analysis in section
\ref{sect:inv.2}, only one singular value is large. The case with $\bS
= (\ve_3)$ is much worse, because all the singular values are small,
and are likely to be dominated by noise. Explicitly, we get
\[
\sigma_1 = 0.145, ~ ~ \sigma_2 = 0.065, ~ ~ \sigma_3 = 0.023
\] 
in large aperture regime and 
\[
\sigma_1 = 2.3 \cdot 10^{-5}, ~ ~ \sigma_2 = 1.3 \cdot 10^{-5}, ~ ~ \sigma_3 = 
1.9 \cdot 10^{-7}
\]
in the small aperture regime. Finally, when both cross-range components of
the electric field are measured using $\bS = (\ve_1, \ve_2)$, the
singular values of $\DD_{_{\hspace{-0.02in} S}}$ are similar to those
of $\DD$. Thus, we expect  inversion results which are comparable to those for complete measurements.

\subsubsection{Estimation of the noise level}
\begin{figure}[t]
\begin{center}
\includegraphics[width=0.35\textwidth]{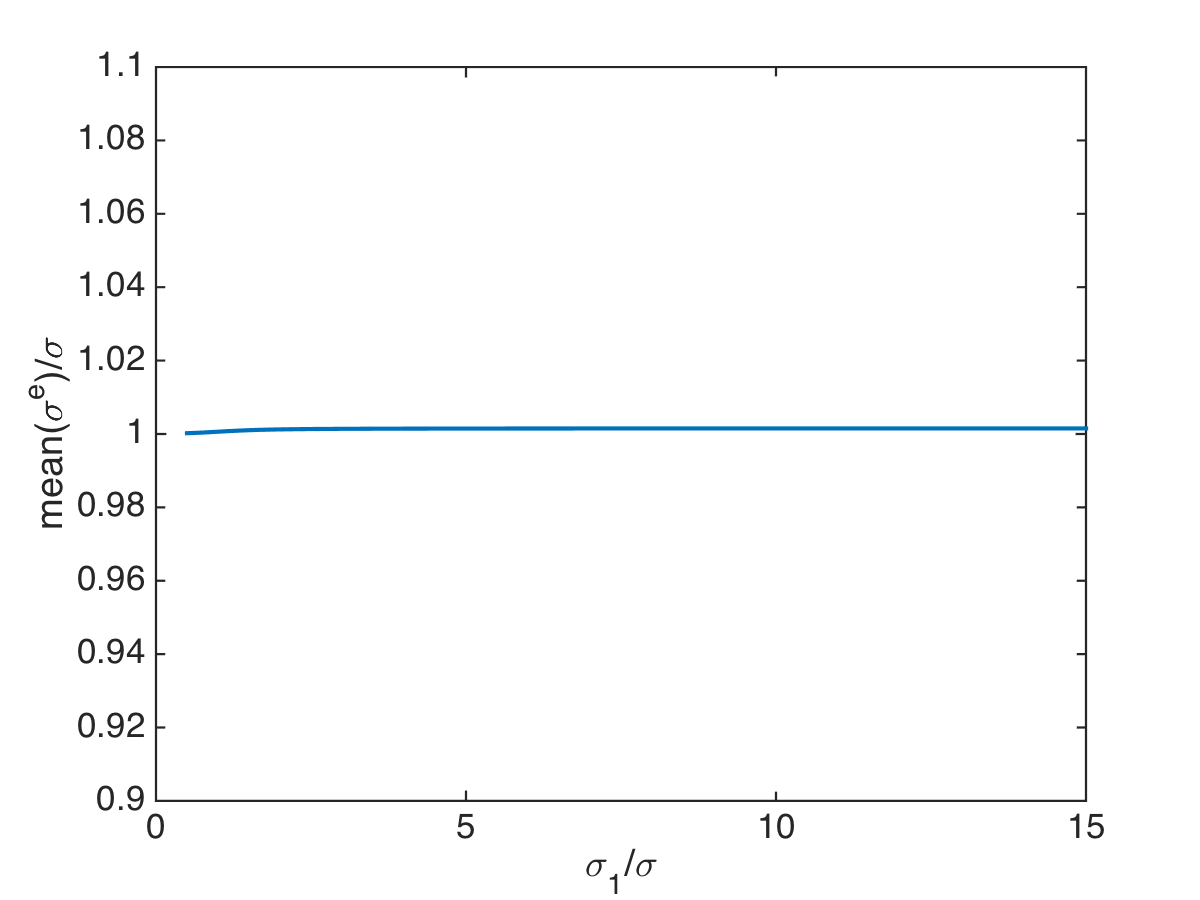}
\includegraphics[width=0.35\textwidth]{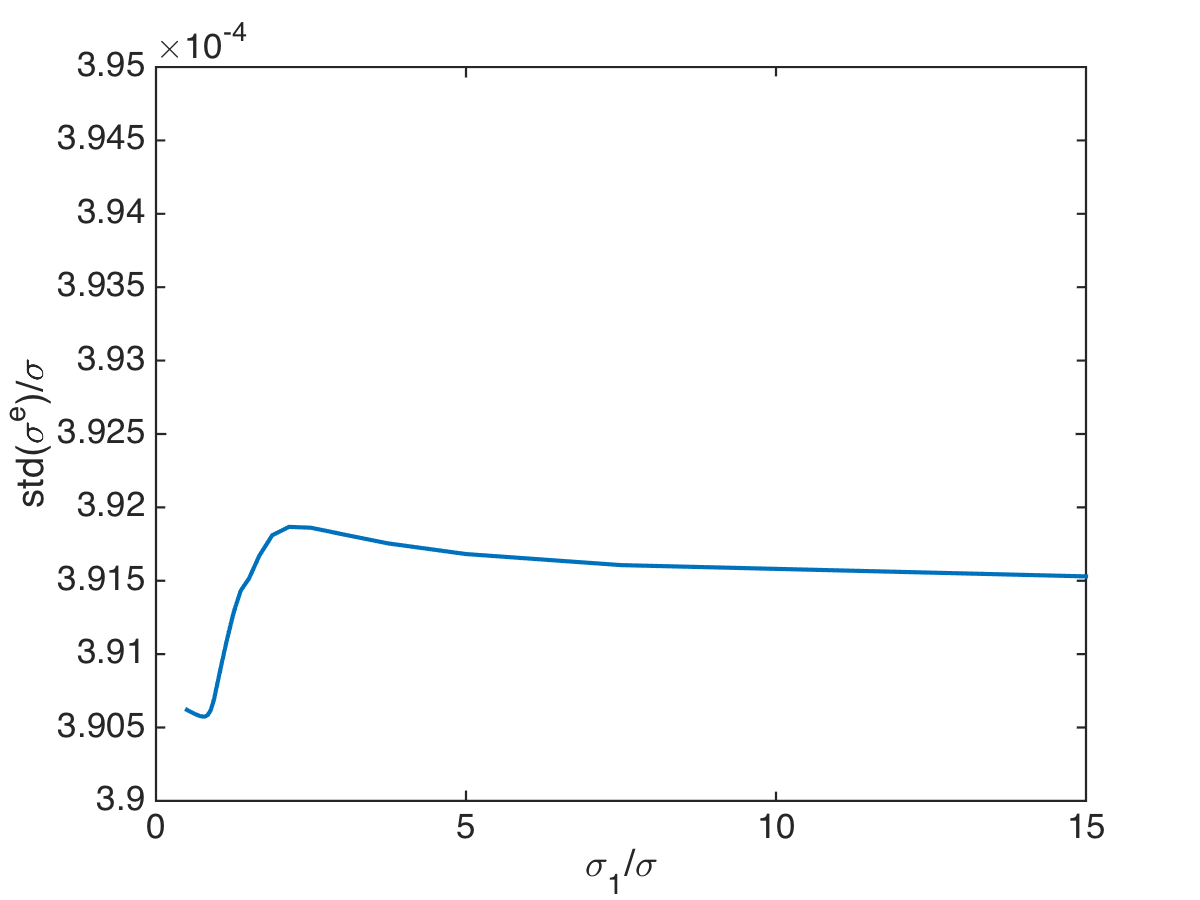}
\end{center}
\vspace{-0.1in}
\caption{Empirical mean of the noise level estimate \eqref{eq:IN2}
  (left plot) and standard deviation (right plot). The results are
  normalized by $\sigma$. The abscissa is
  $\sigma_1/\sigma$.}
\label{fig:NLEV}
\end{figure}

We display in Figure \ref{fig:NLEV} the empirical mean and standard
deviation of the unbiased estimate \eqref{eq:IN2} of the noise level,
using the complete measurement matrix $\tDD$.  We note that $ \sigma^e
\approx \mbox{mean}[\sigma^e] \approx \sigma.  $ The estimate
\eqref{eq:IN3} is very similar, because we have many sensors. For the
same reason we obtain similar estimates of $\sigma$ from incomplete
measurements.
If the number of sensors were smaller, then the refined estimation methods proposed in \cite{garnier2014applications} could be used.

\subsection{Inversion results for one inclusion}
\label{sect:num.2}
We present results for the inclusion at $\vy_1 =
(\lambda,-\lambda,L)$, with reflectivity tensor \eqref{eq:rhoT1}, in both the
large and small aperture regimes. We compare the MUSIC method
\eqref{eq:IN6} with the imaging function \eqref{eq:IN13}, and estimate
the reflectivity tensor as in section \ref{sect:getref}.

\begin{figure}[h!]
\begin{center}
\includegraphics[width=0.4\textwidth]{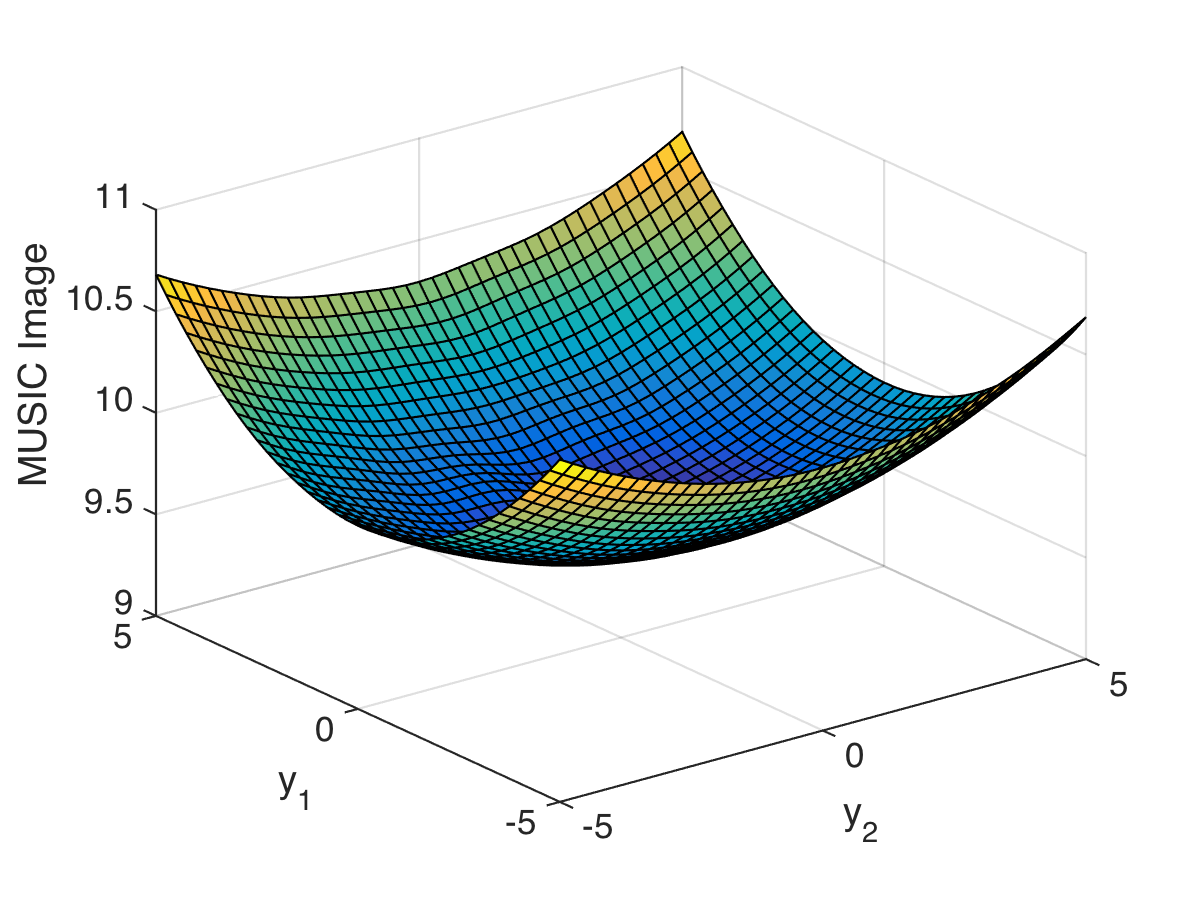}
\includegraphics[width=0.4\textwidth]{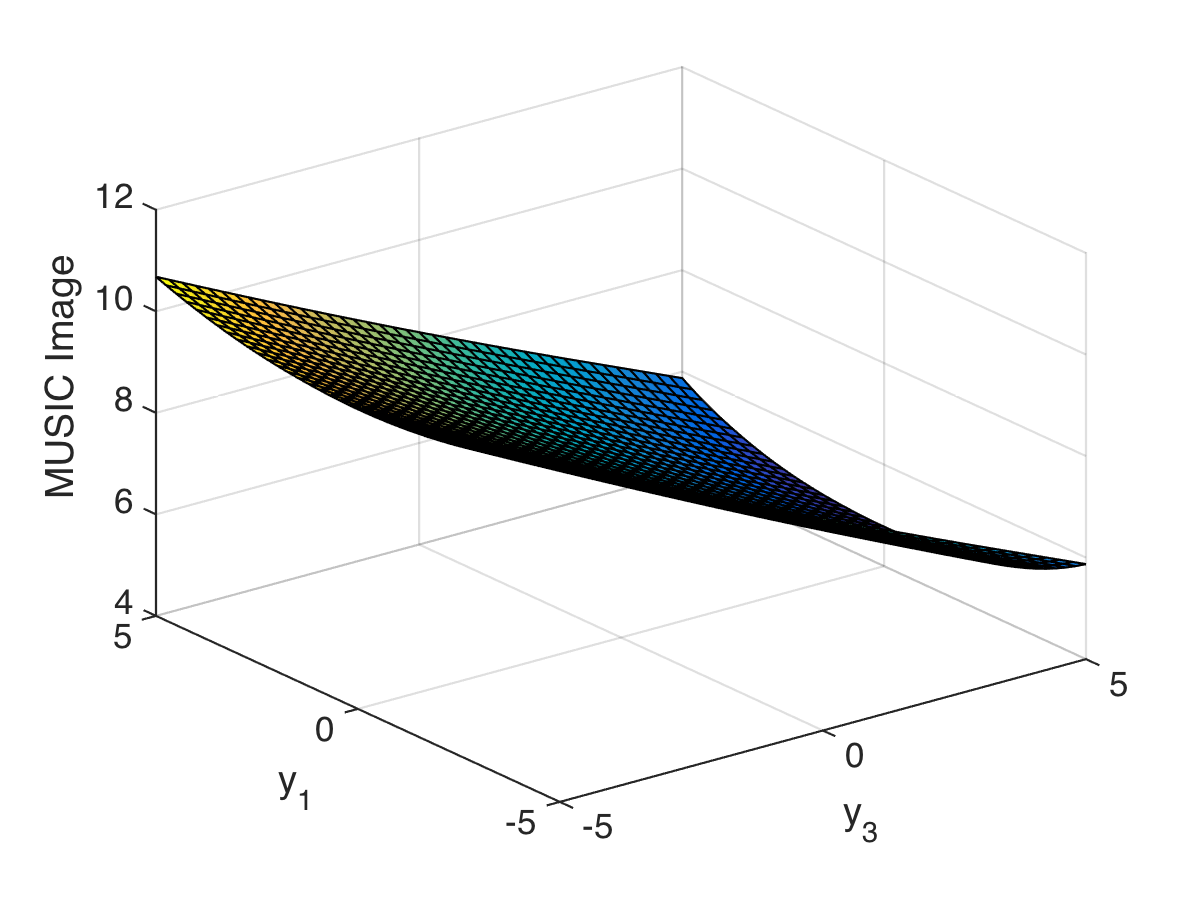}
\\ \vspace{-0.1in} \includegraphics[width=0.4\textwidth]{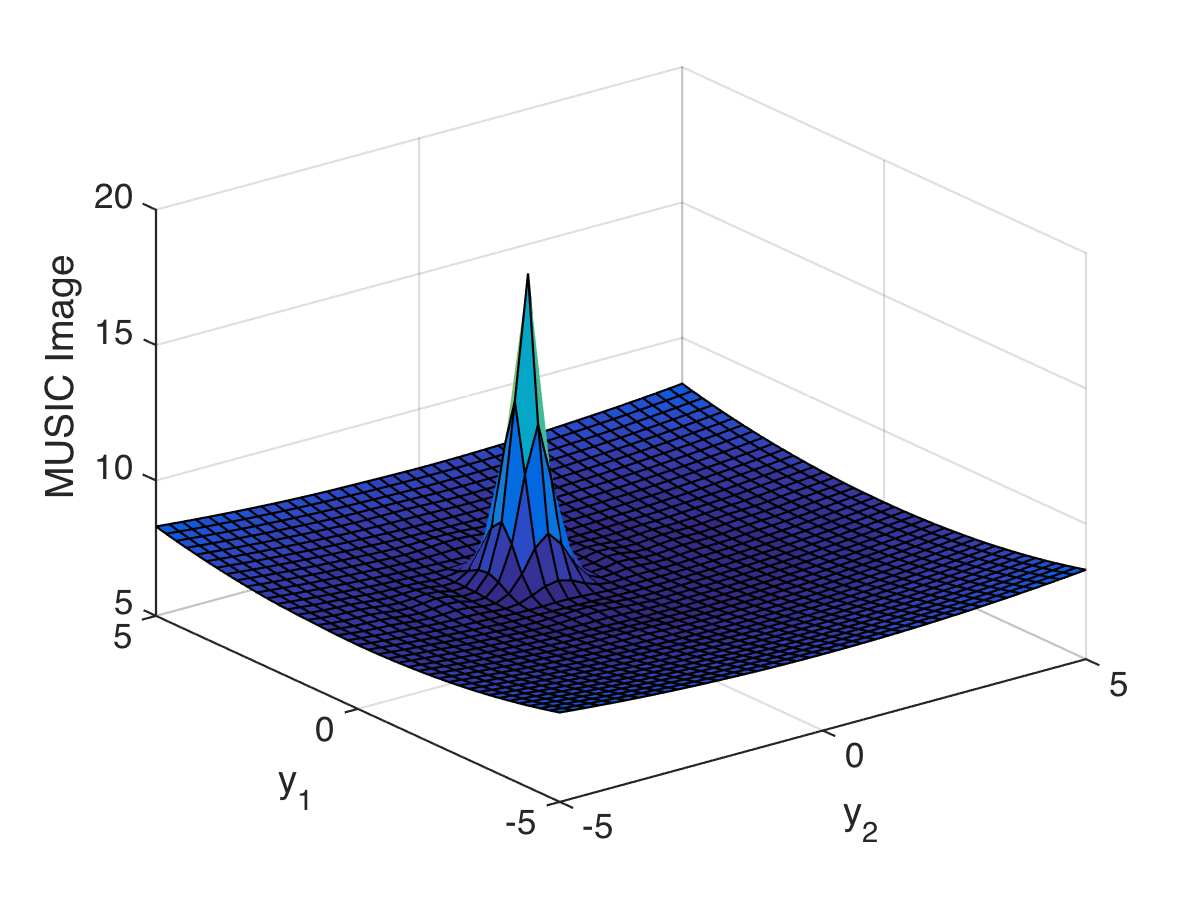}
\includegraphics[width=0.4\textwidth]{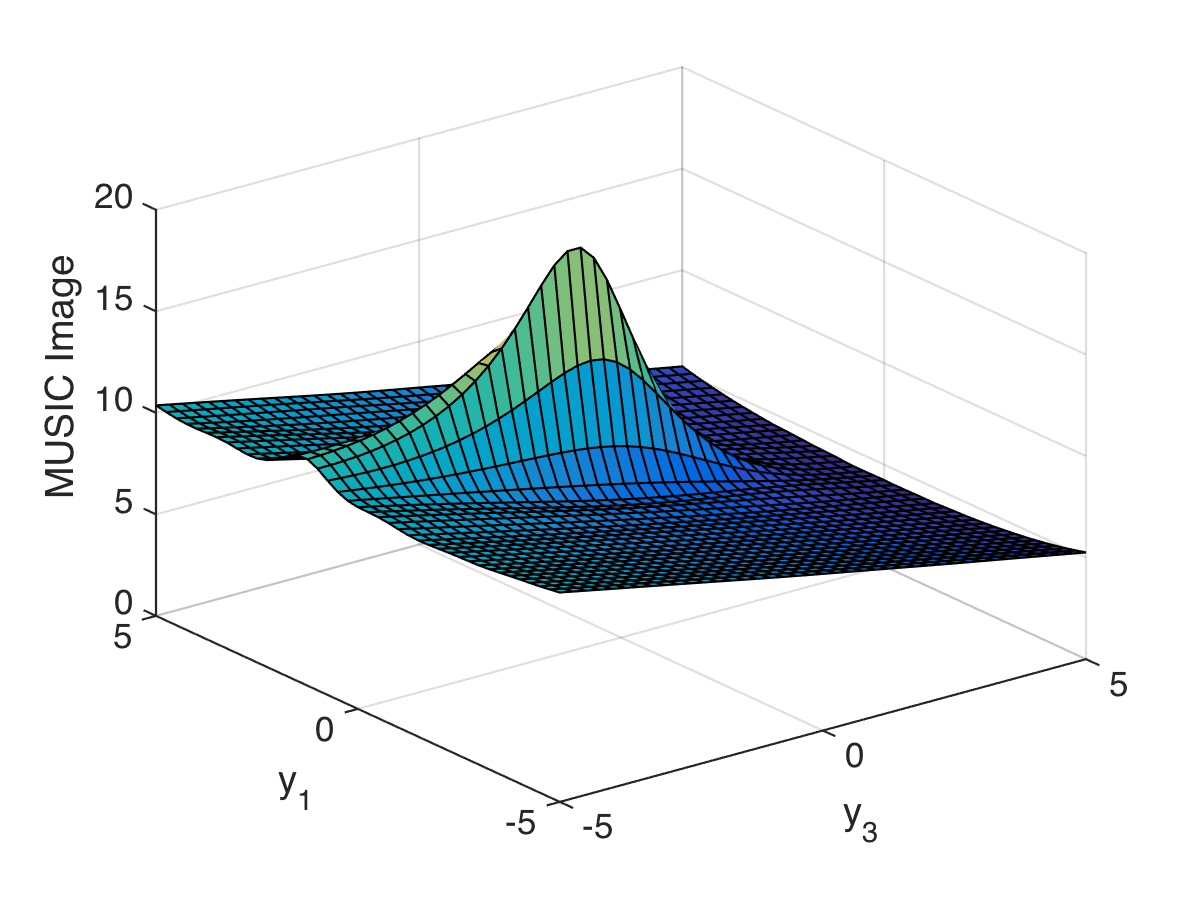}
\\
\end{center}
\vspace{-0.1in}
\caption{MUSIC imaging of one inclusion in the large aperture regime and
  complete measurements. We display $\mathcal{I}_{_{MUSIC}}(\vy)$ in
  the plane $y_3 = L$ (left) and $y_2 = -\lambda$ (right). The axes
  are in units of the wavelength. The top row is for $75\%$ noise and
  the bottom  row for $25\%$.}
\label{fig:MUSIC}
\end{figure}
\subsubsection{Complete measurements} 
\begin{figure}[t]
\begin{center}
\includegraphics[width=0.4\textwidth]{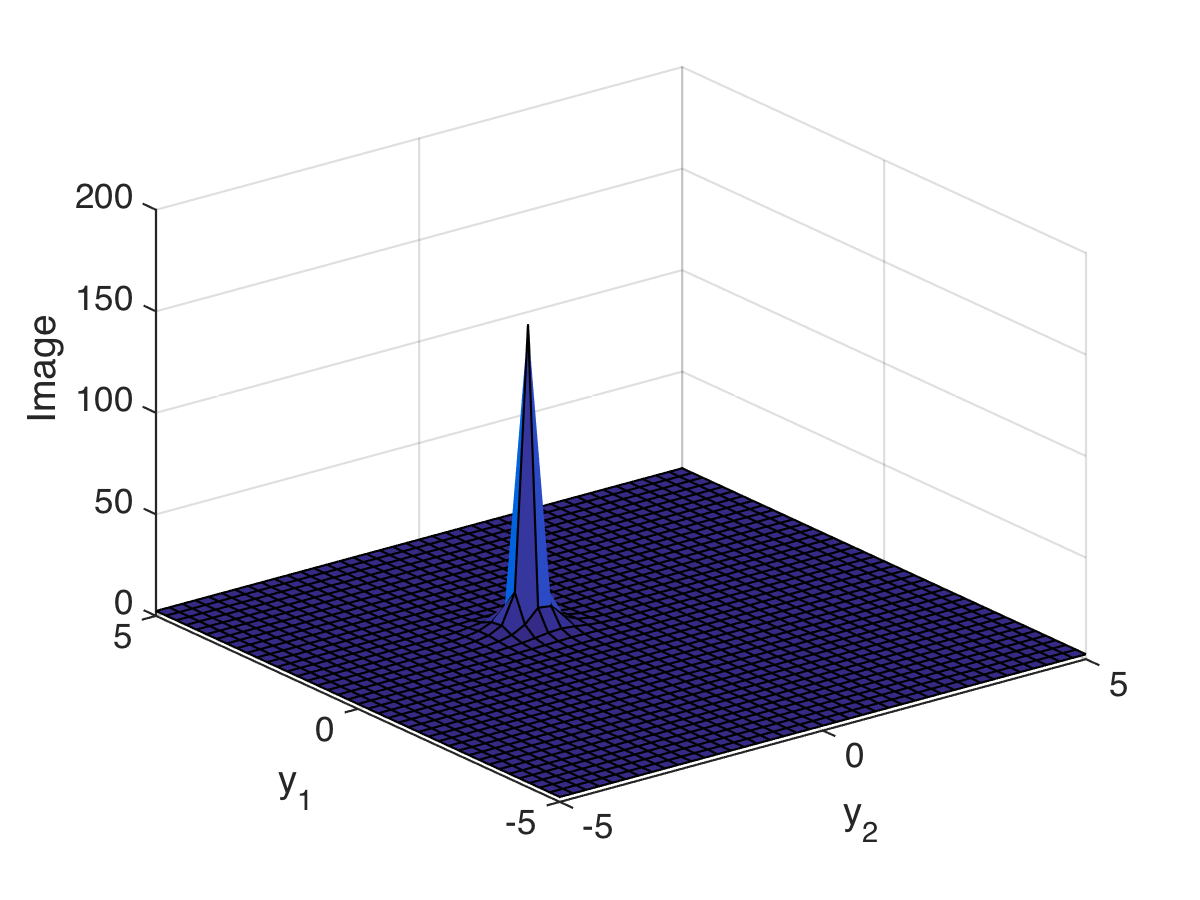}
\includegraphics[width=0.4\textwidth]{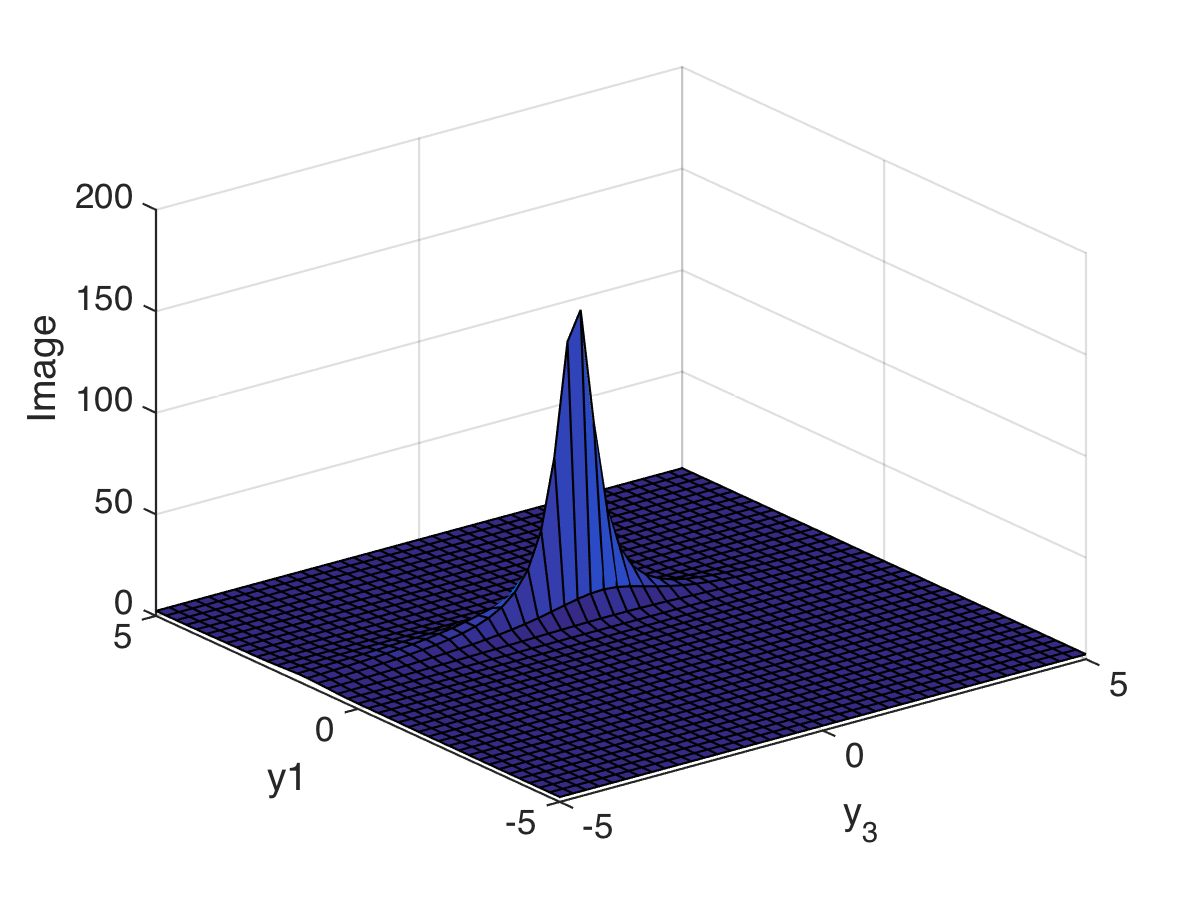} 
\end{center}
\vspace{-0.1in}
\caption{Imaging function \eqref{eq:IN13} in the large aperture regime, for 
  complete measurements and $75\%$ noise. We display $\mathcal{I}(\vy)$ in
  the plane $y_3 = L$ (left) and $y_2 = -\lambda$ (right). The axes
  are in units of the wavelength. }
\label{fig:imagesNF}
\end{figure}

Let us begin with Figure \ref{fig:MUSIC},  which shows the surface plots of  $\mathcal{I}_{_{MUSIC}}(\vy)$ in the planes $y_3 = L$
and $y_2 = -\lambda$, respectively. These are obtained in the large aperture
regime, for two realizations of the noise matrix: one for $25\%$ noise
(bottom row) and the other for $75\%$ noise (top row). At the higher
noise level the MUSIC function fails to localize the inclusion, as it
has no peak. For comparison, we display in Figure \ref{fig:imagesNF}
the results obtained with our imaging function \eqref{eq:IN13}, for
the same realization at $75\%$ noise. We do not show the images for
$25\%$ noise because they are similar. The images are a significant
improvement over those in Figure \ref{fig:MUSIC}, and they clearly
localize the inclusion in range and cross-range.

\begin{figure}[h]
\begin{center}
\includegraphics[width=0.3\textwidth]{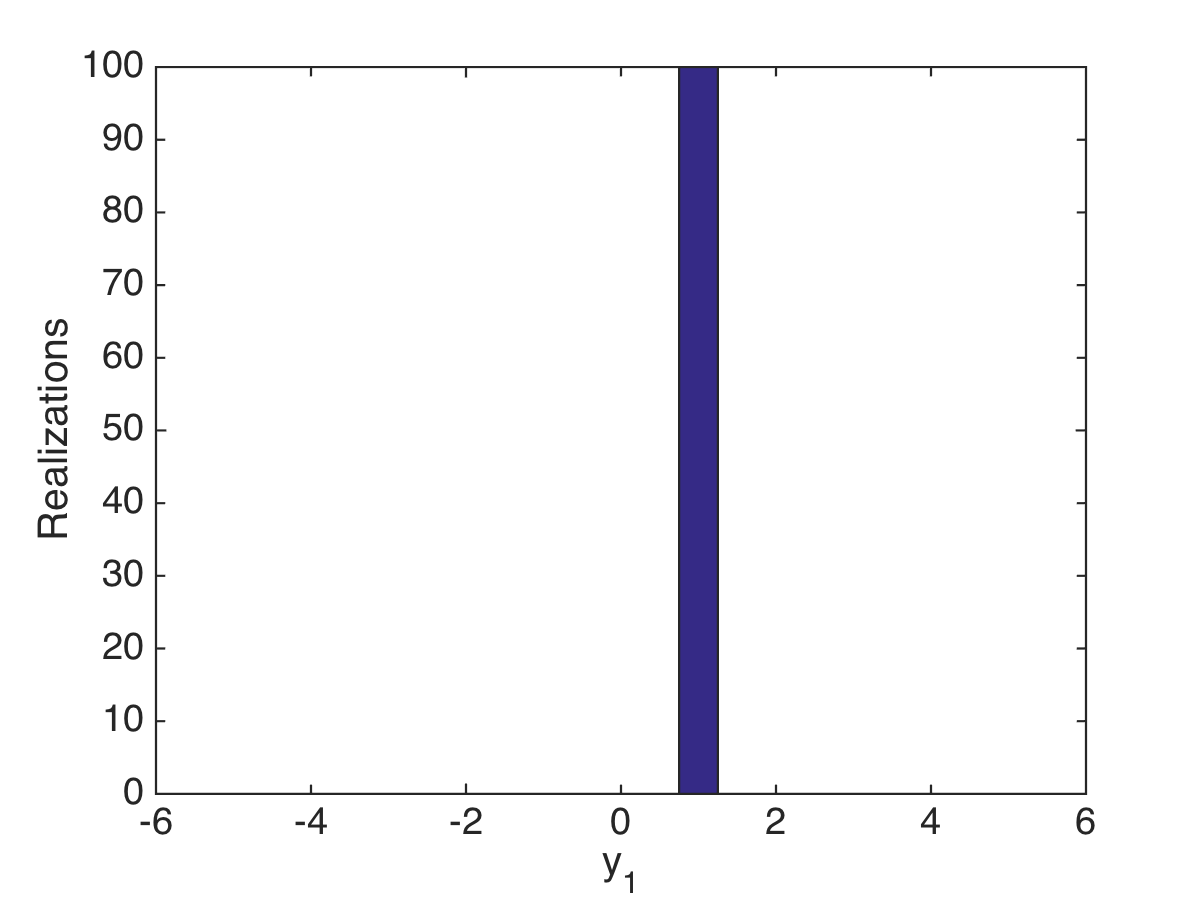}
\includegraphics[width=0.3\textwidth]{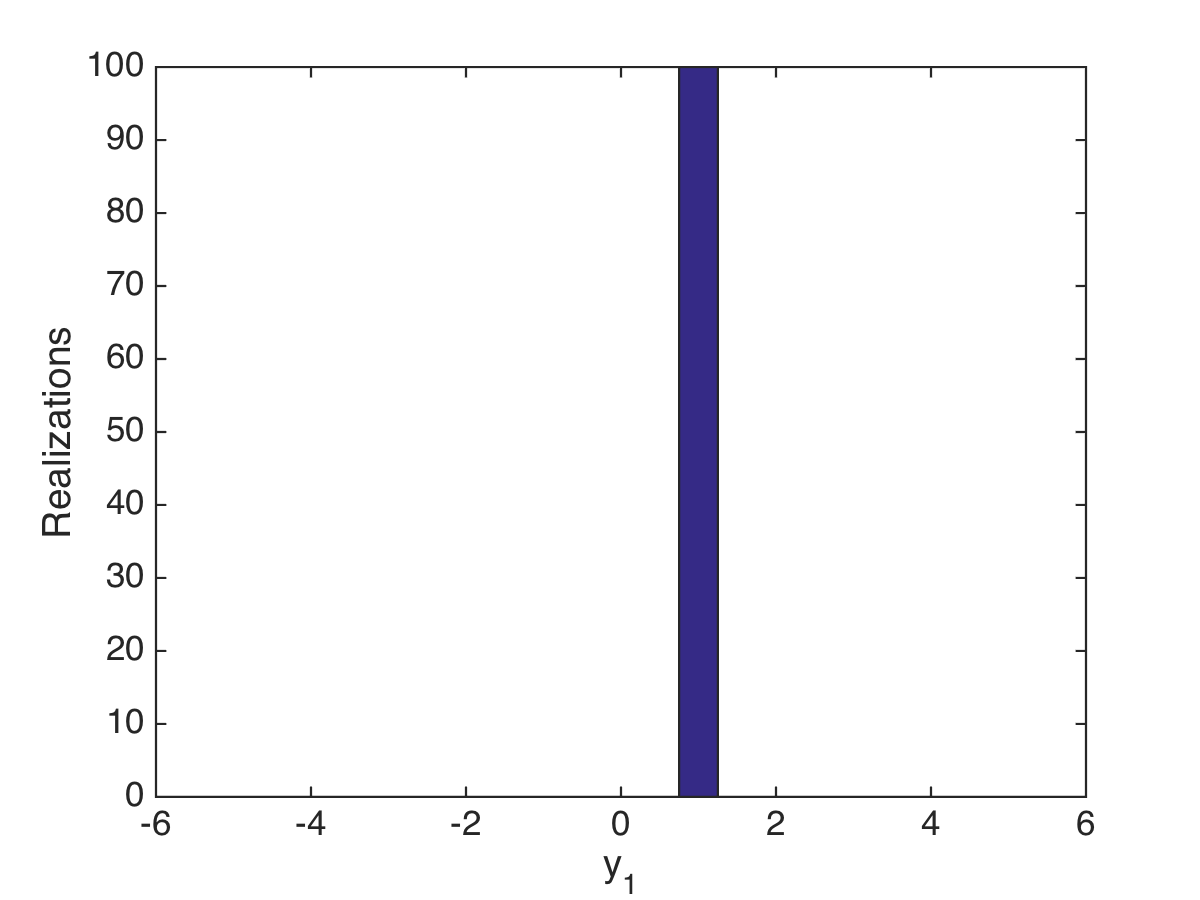}
\includegraphics[width=0.3\textwidth]{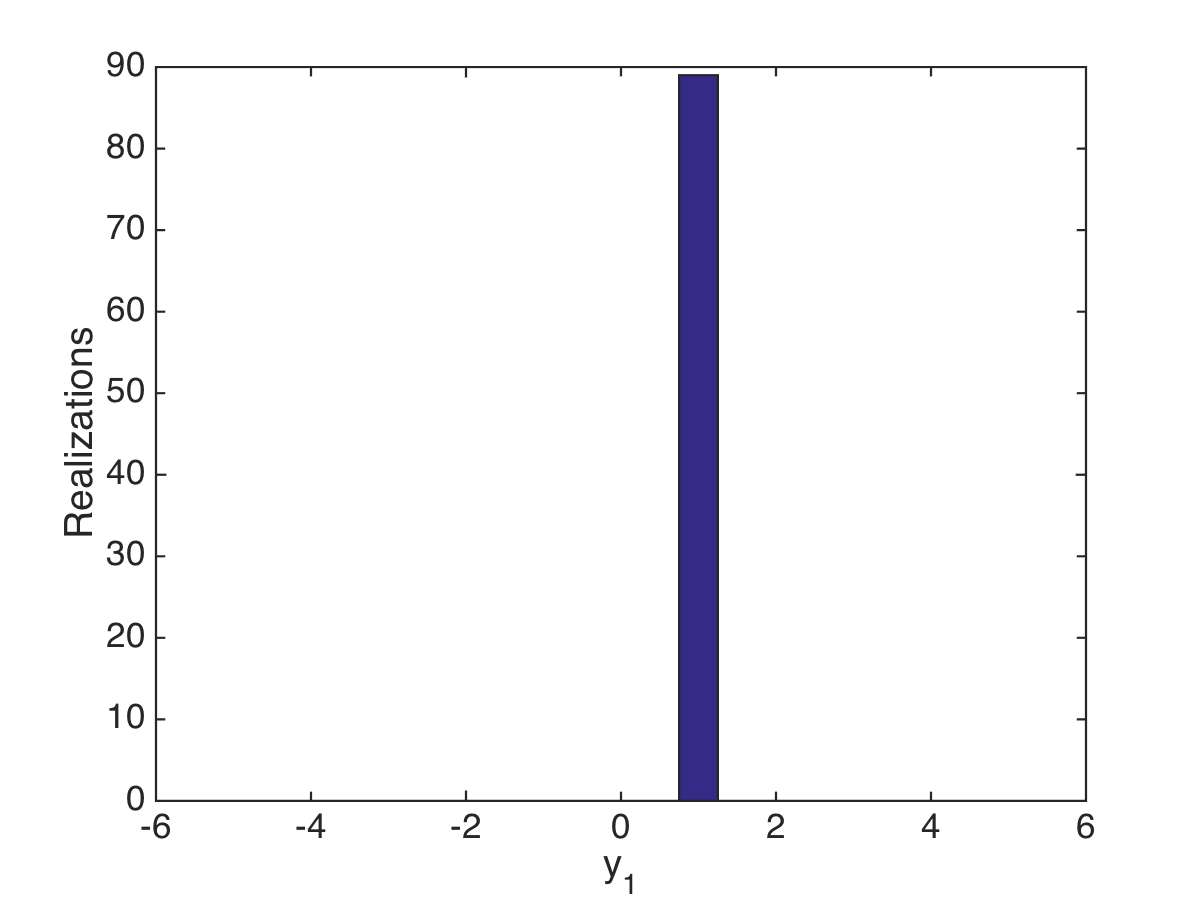}
\\ \vspace{-0.02in}\includegraphics[width=0.3\textwidth]{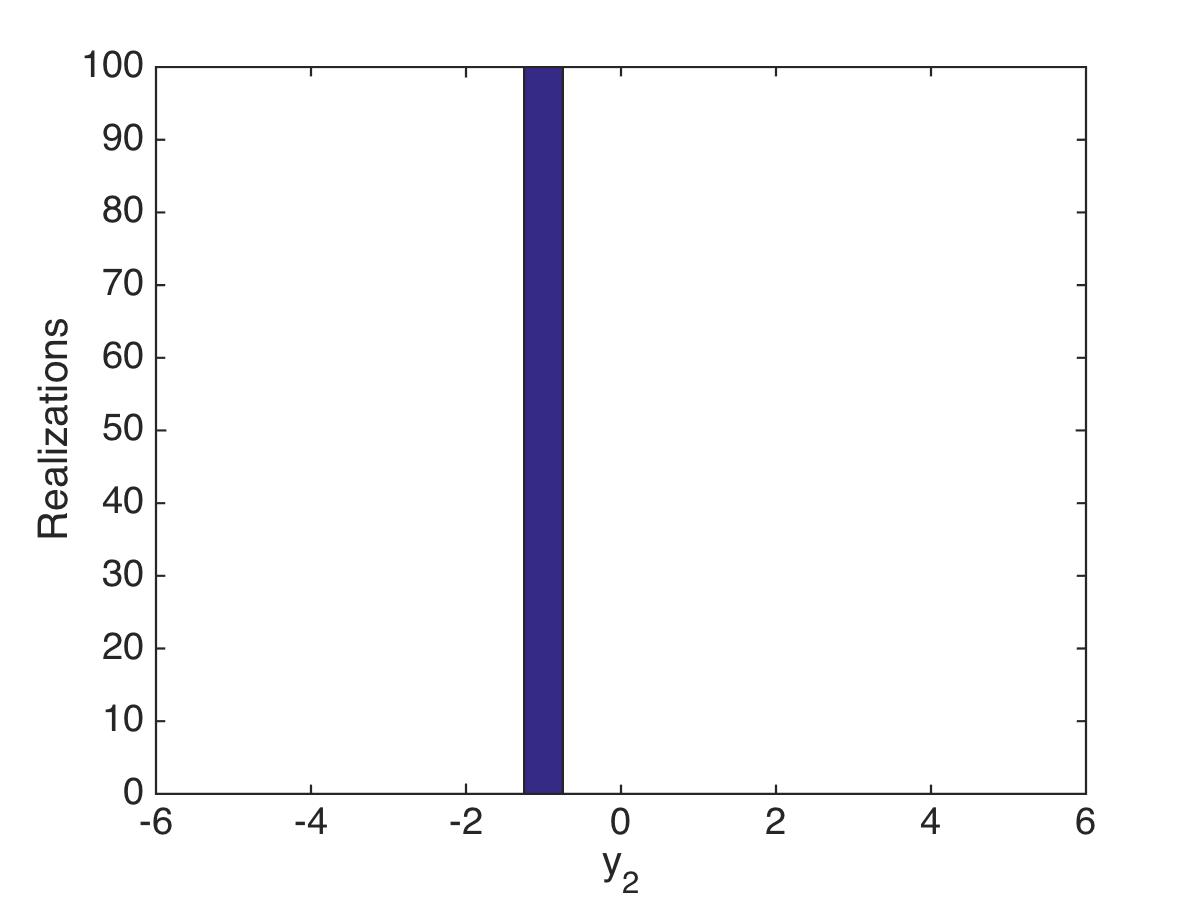}
\includegraphics[width=0.3\textwidth]{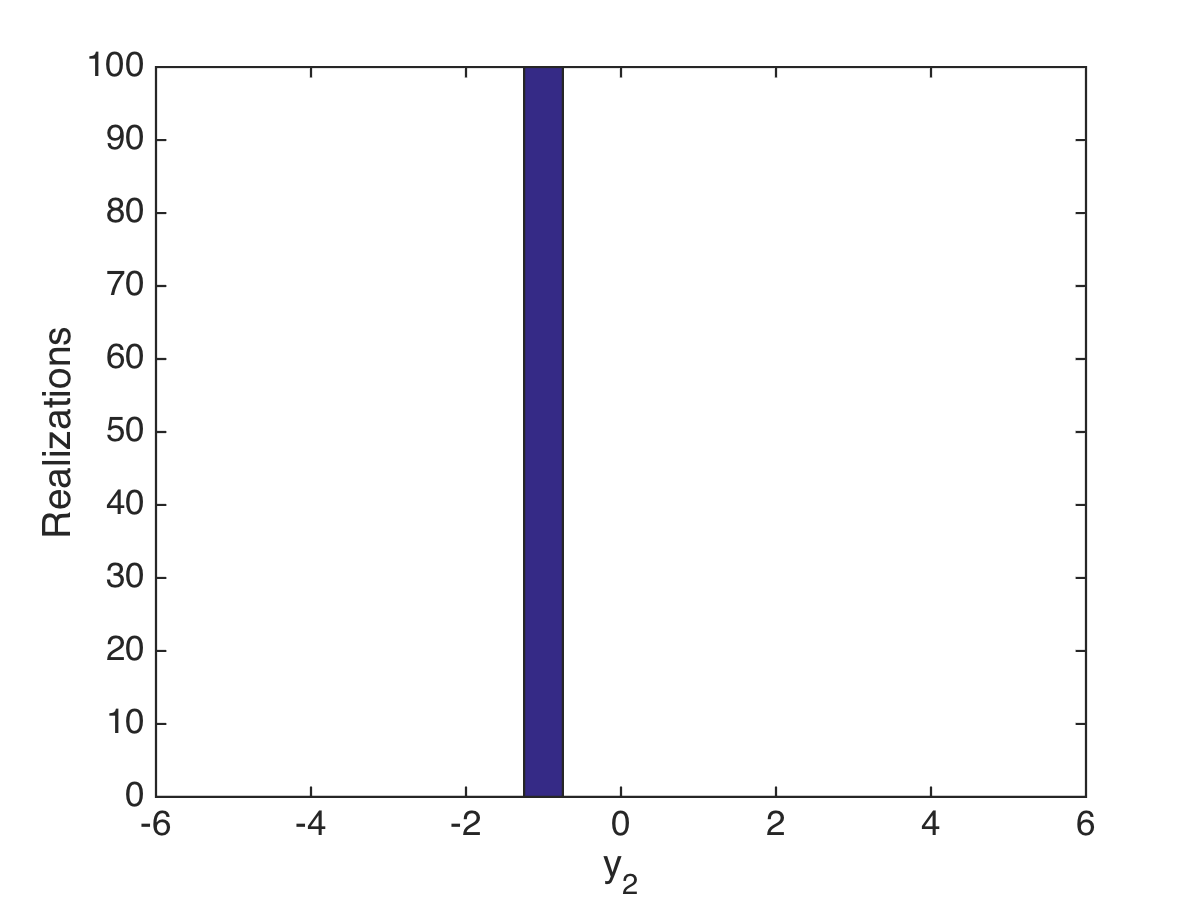}
\includegraphics[width=0.3\textwidth]{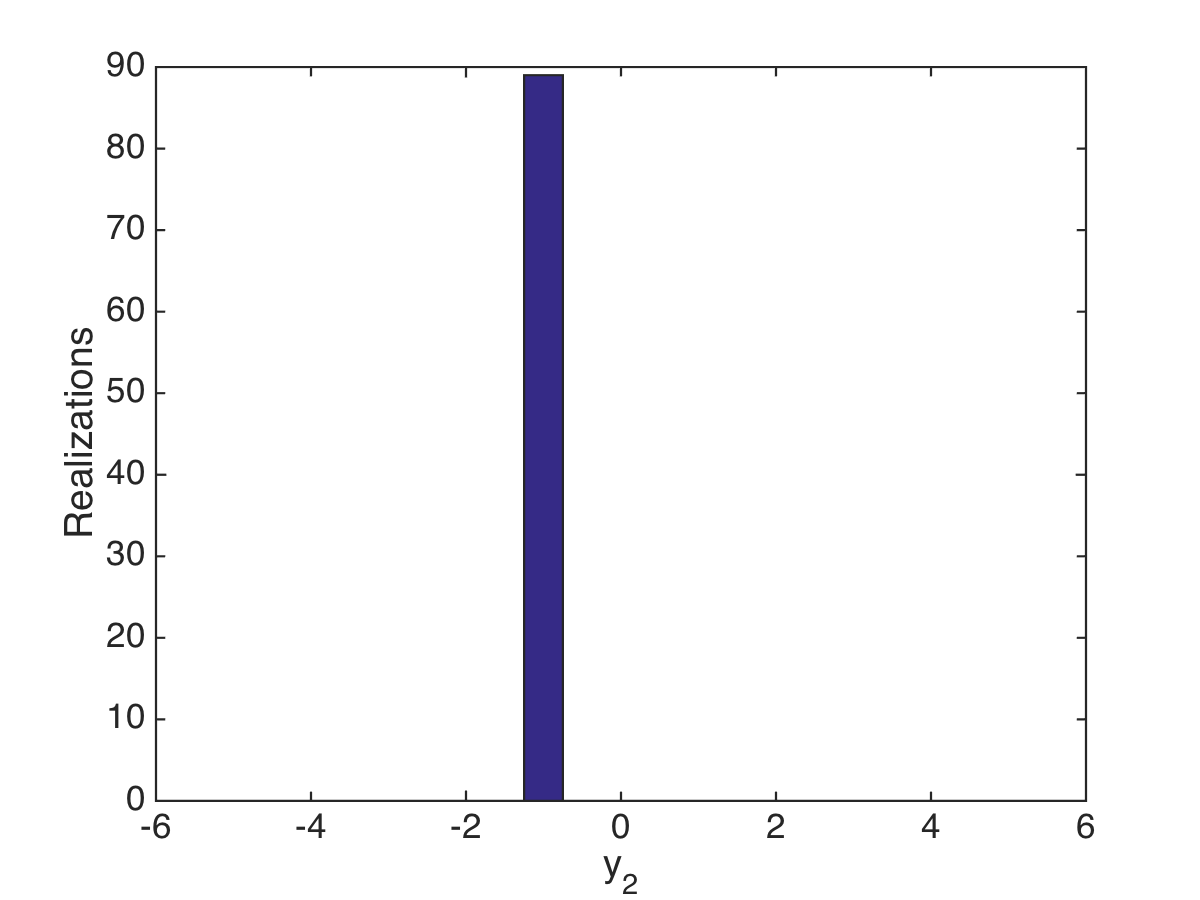}
\\ \vspace{-0.02in}\includegraphics[width=0.3\textwidth]{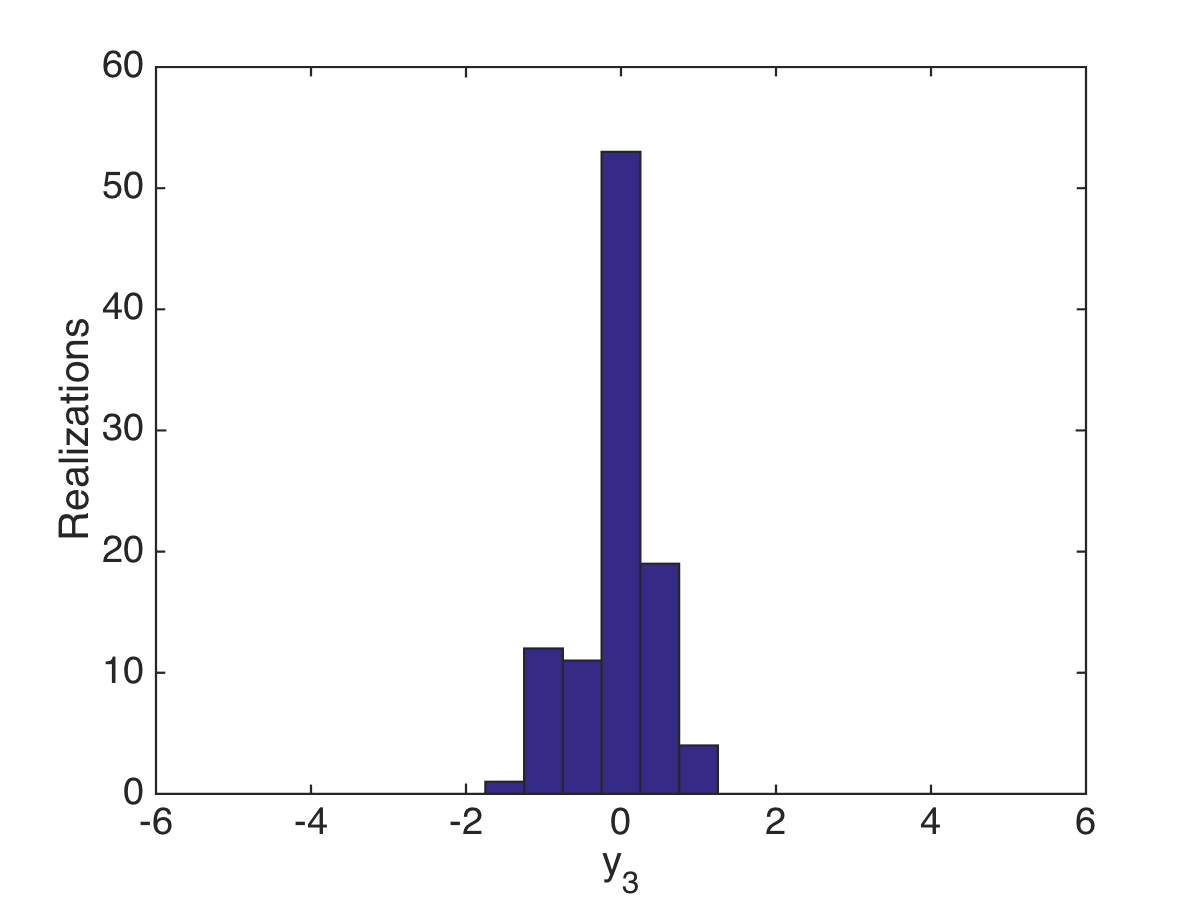}
\includegraphics[width=0.3\textwidth]{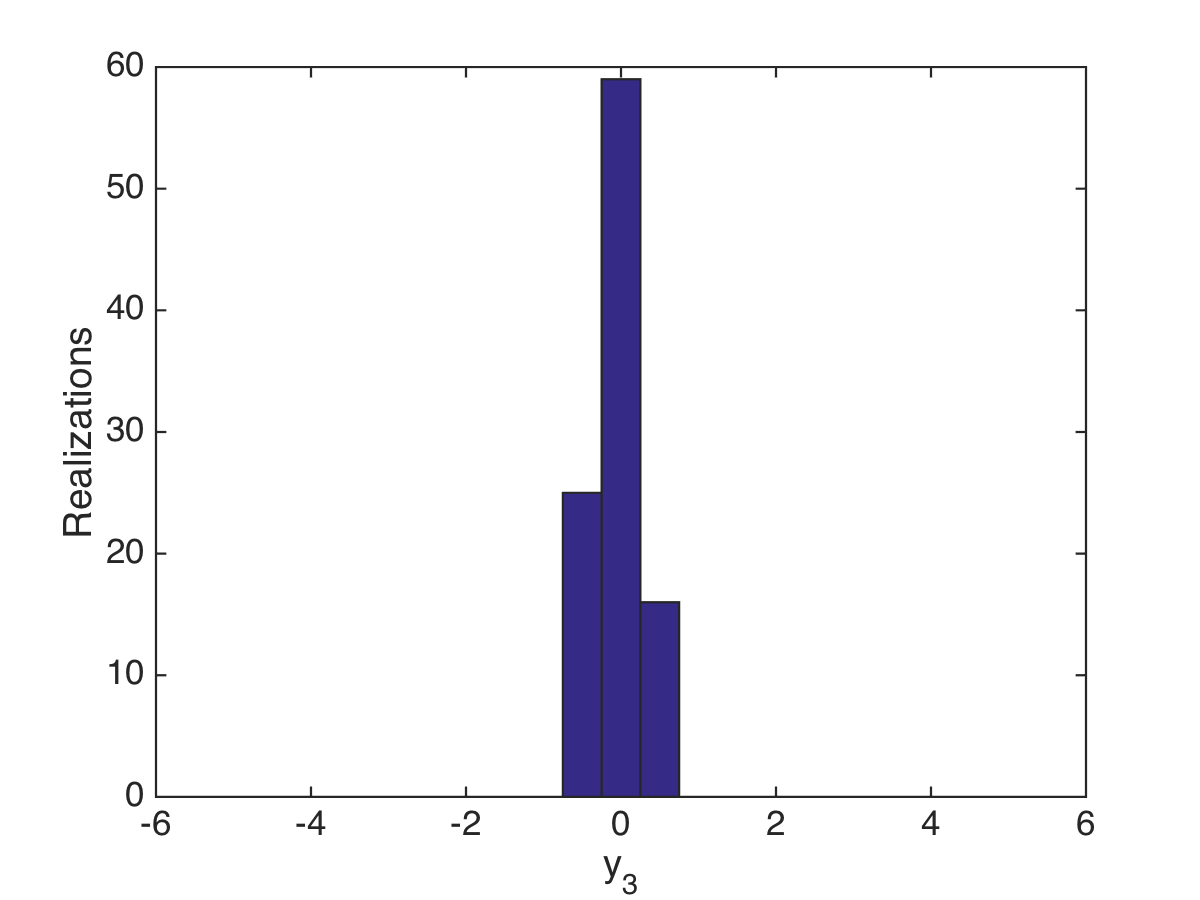}
\includegraphics[width=0.3\textwidth]{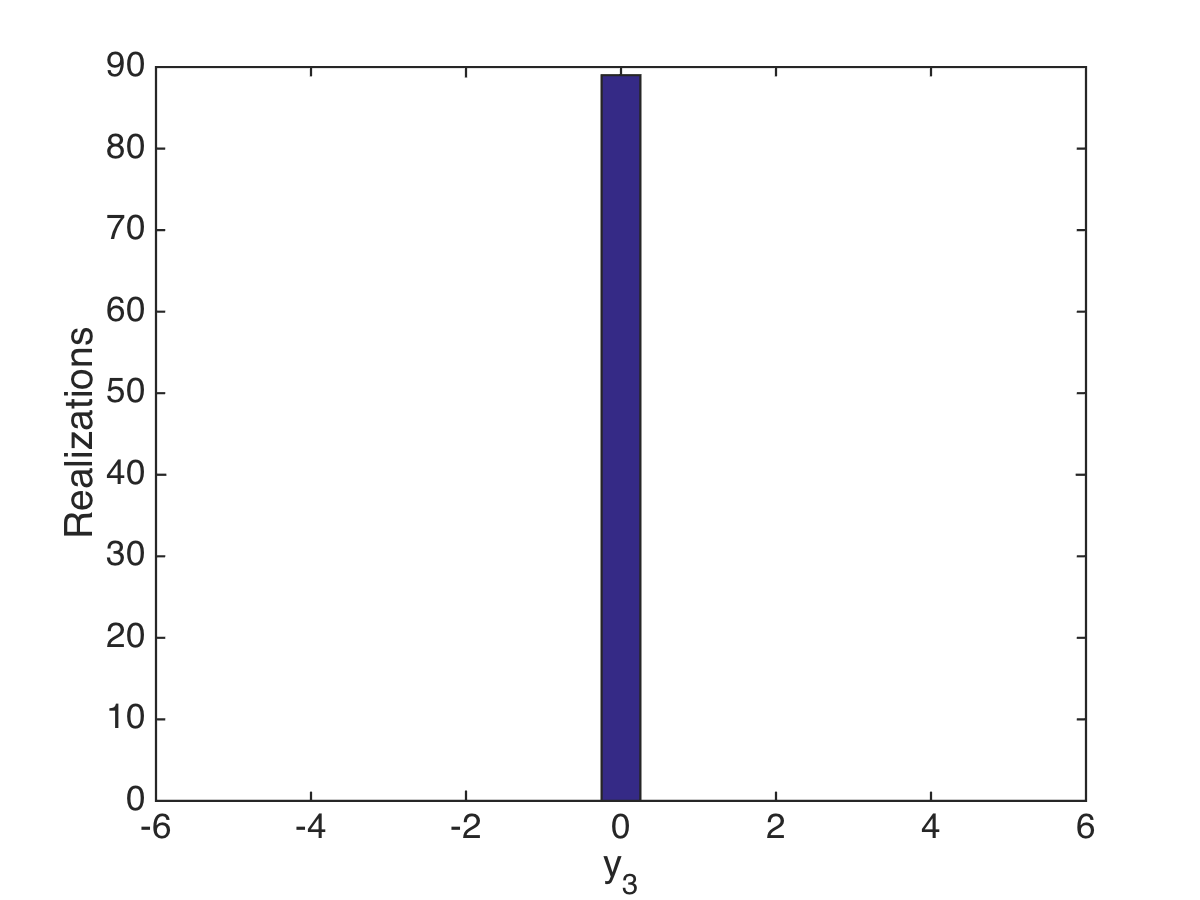}
\end{center}
\vspace{-0.1in}
\caption{The histograms of the peak location in $y_1$ (top), $y_2$
  (middle) and $y_3$ (bottom) of the imaging function \eqref{eq:IN13}
  in the large aperture regime.  The left column is for $\sigma/\sigma_1 =
  75\%$, the middle column for $50\%$ and the right column for
  $25\%$. The abscissa is in units of the wavelength.}
\label{fig:HIST_NF}
\end{figure}

To illustrate the statistical stability of our imaging function
\eqref{eq:IN13}, we display  in Figure \ref{fig:HIST_NF} the histograms of
the peak locations along the three coordinate axes, for one hundred
realizations of the noise at $25\%$, $50\%$ and $75\%$ noise level. The effective rank of 
the noisy data matrix $\tDD$ is 
$\widetilde{\mathfrak{R}} = 1$  at $75\%$ noise level and 
$\widetilde{\mathfrak{R}} = 2$ for the weaker noise. We
note that the focusing in cross-range is perfect for all
realizations. The range focusing suffers slightly at the higher noise
levels, but the error is at most $1.5\la$.

\begin{figure}[h]
\begin{center}
\includegraphics[width=0.3\textwidth]{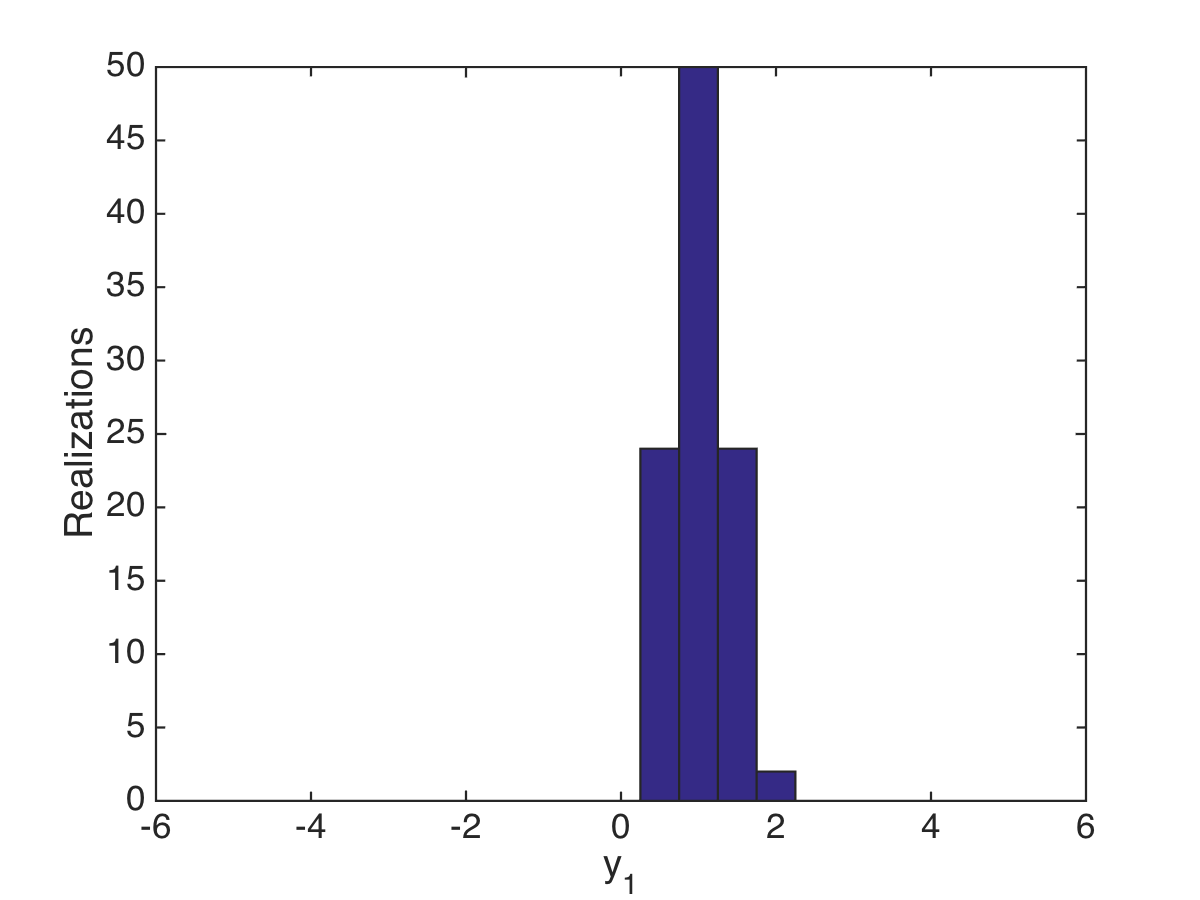}
\includegraphics[width=0.3\textwidth]{Figures/FF_COMPLETE/SIGMA05/HIST_MAXIMA1}
\includegraphics[width=0.3\textwidth]{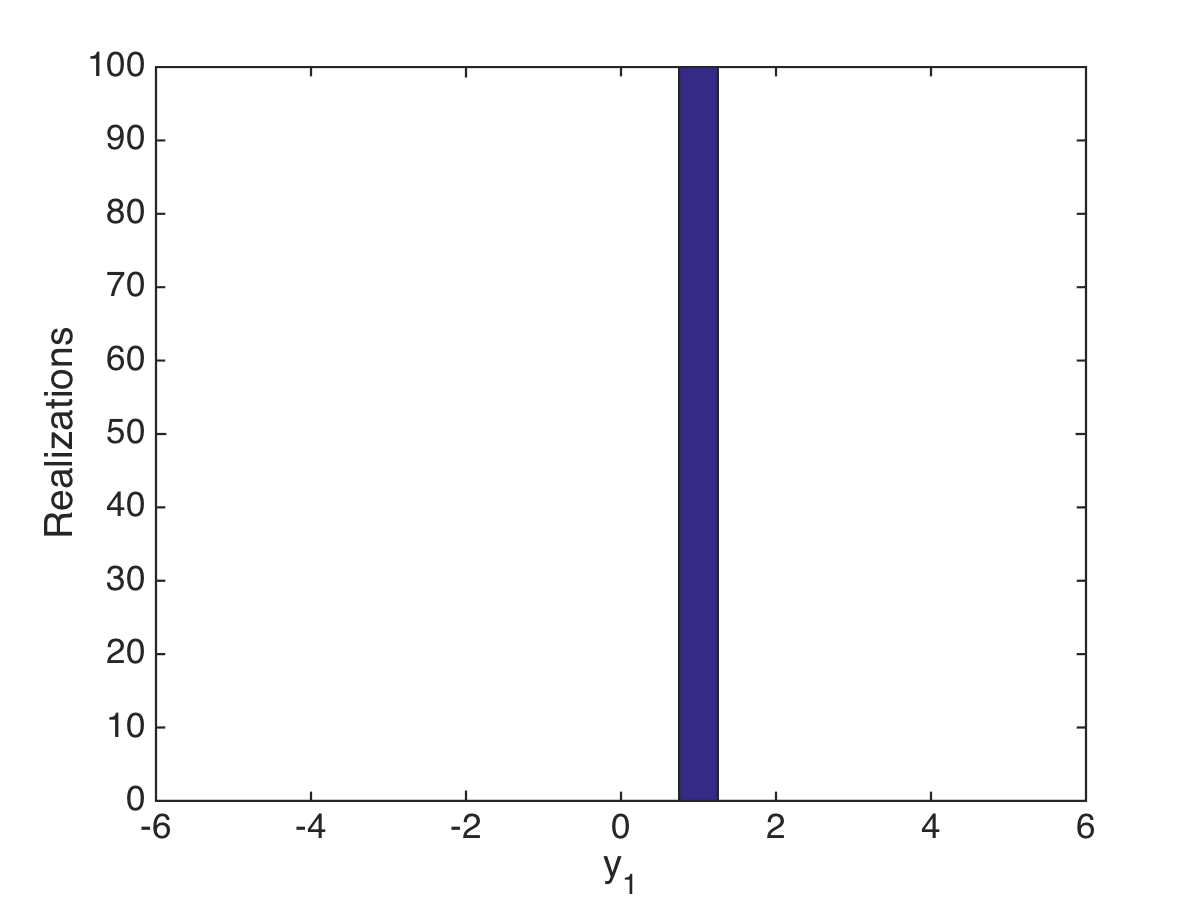} \\
\vspace{-0.02in}\includegraphics[width=0.3\textwidth]{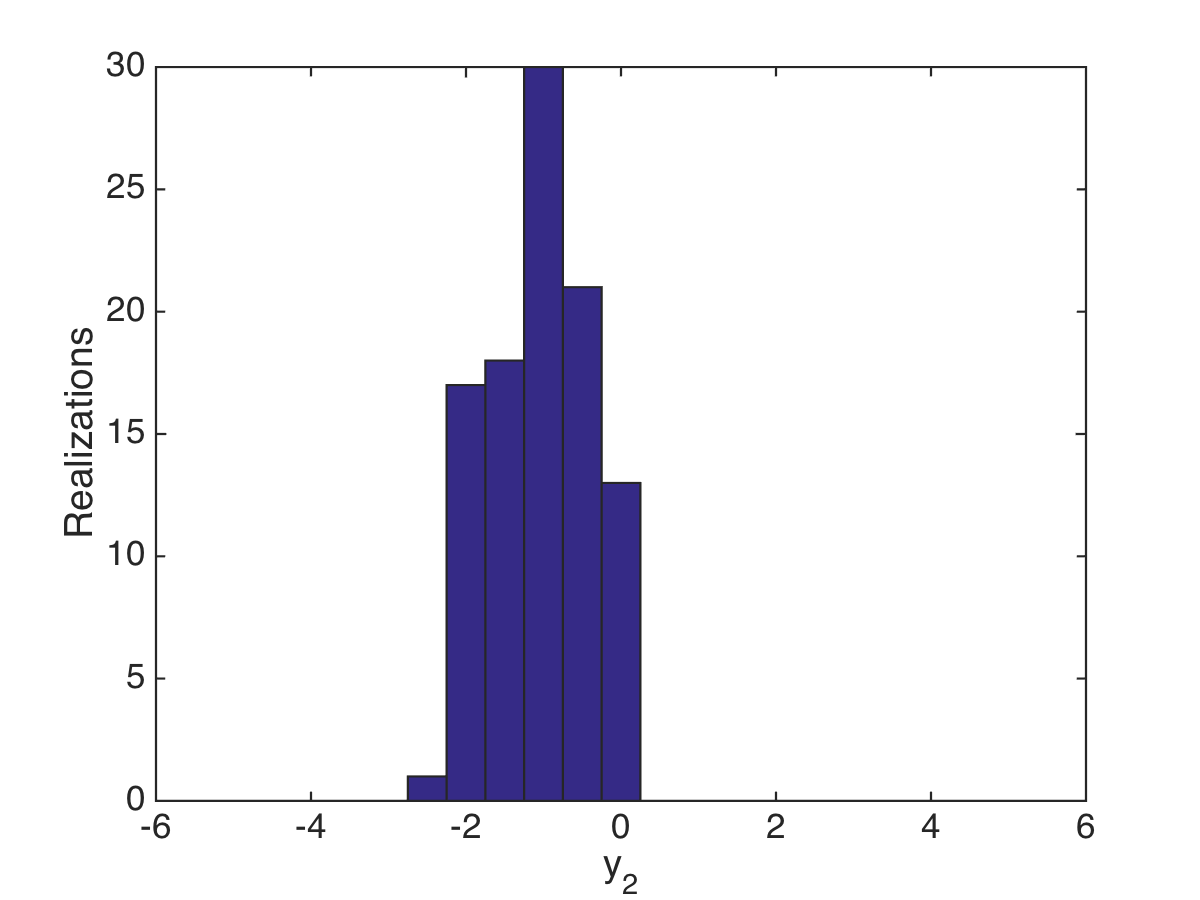}
\includegraphics[width=0.3\textwidth]{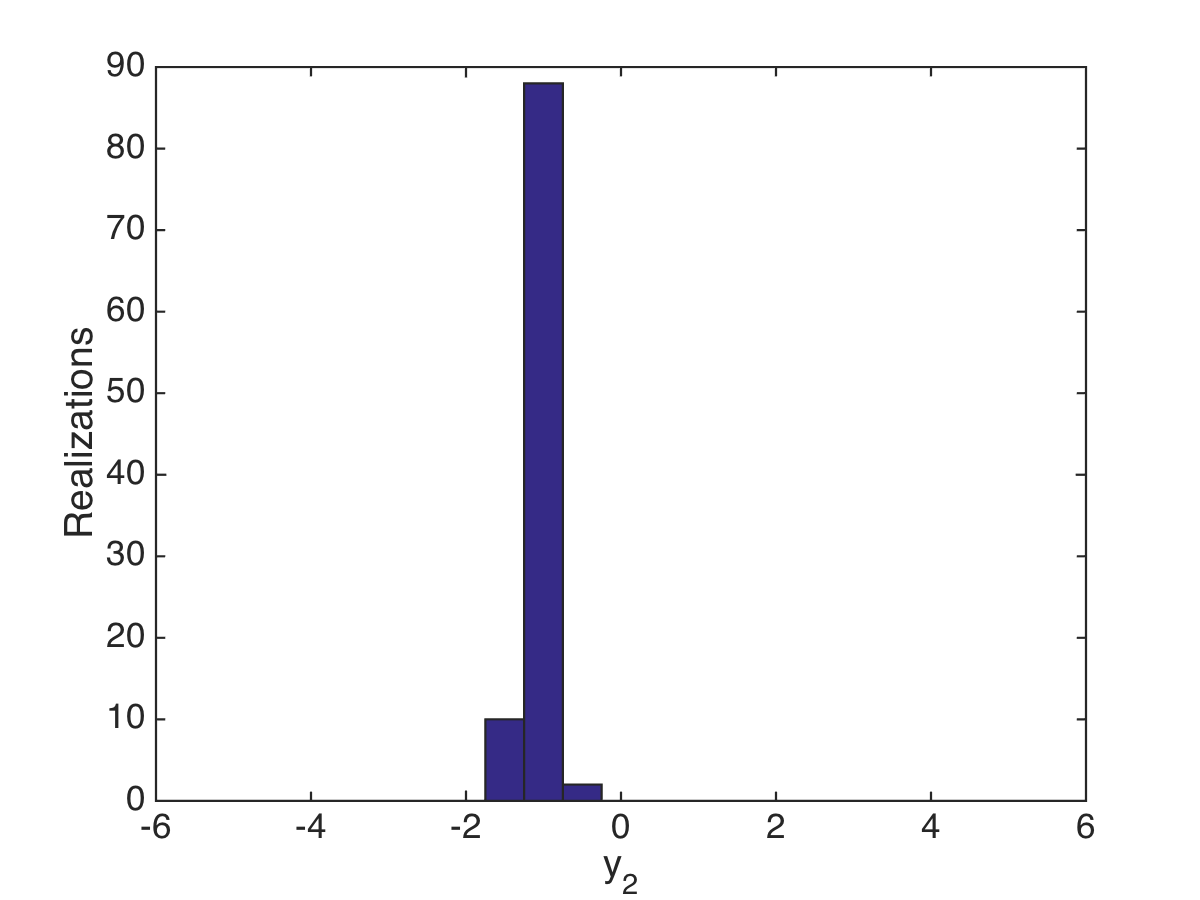}
\includegraphics[width=0.3\textwidth]{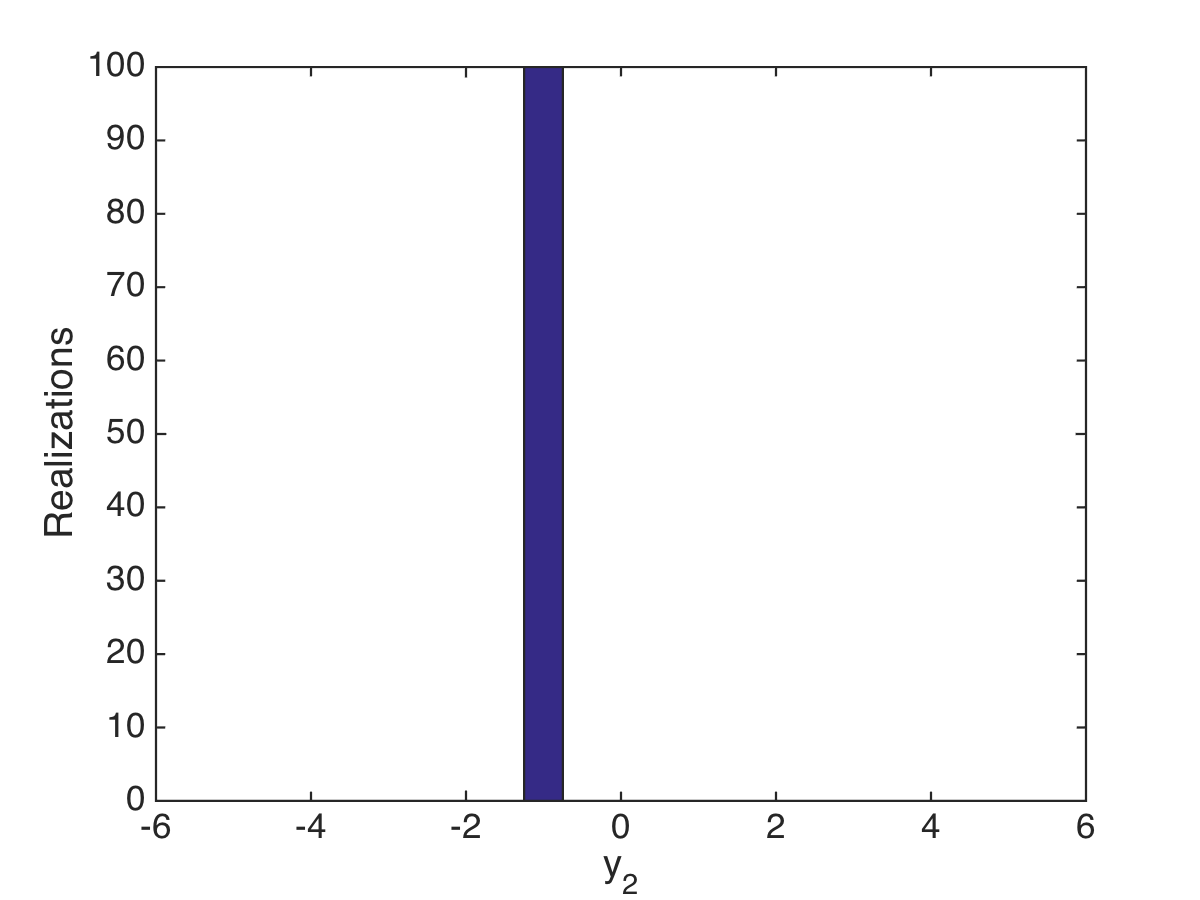} \\
\vspace{-0.02in}\includegraphics[width=0.3\textwidth]{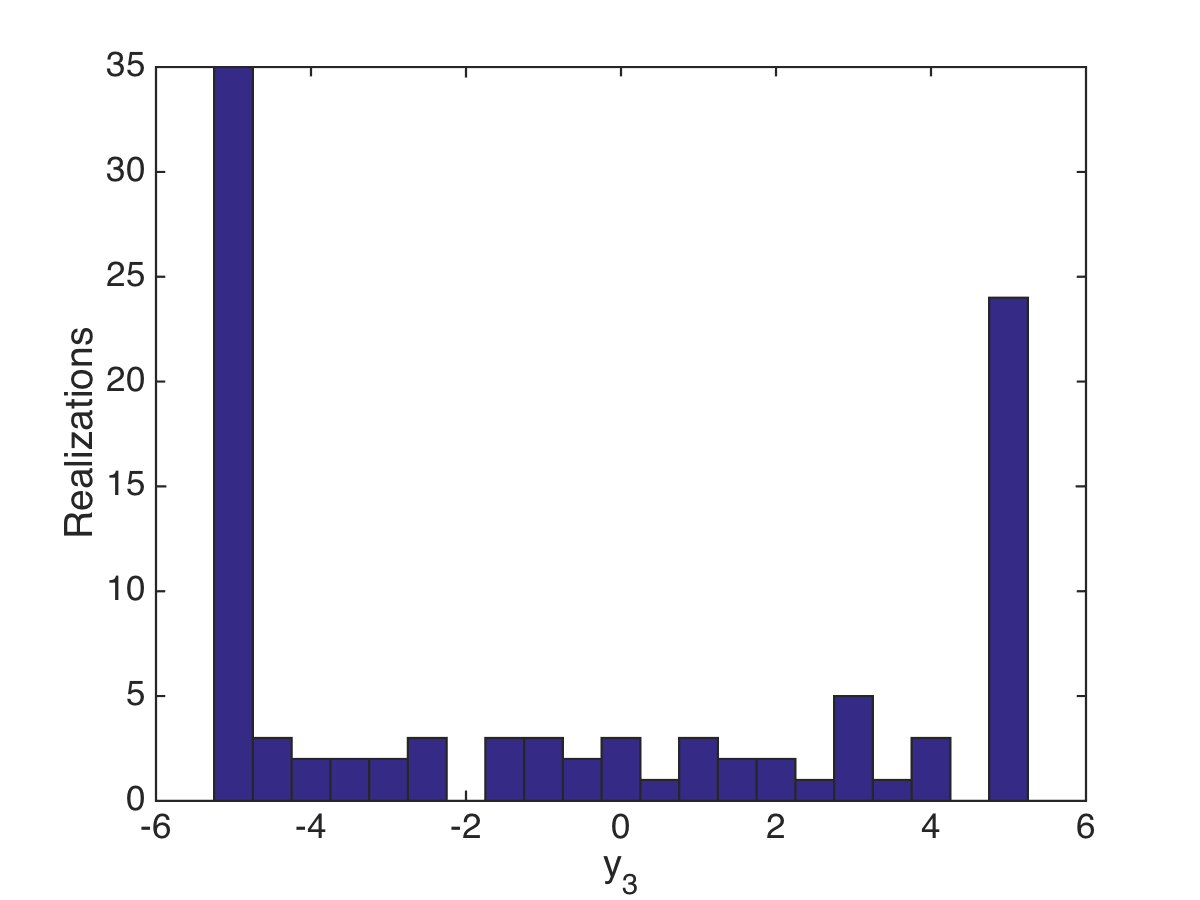}
\includegraphics[width=0.3\textwidth]{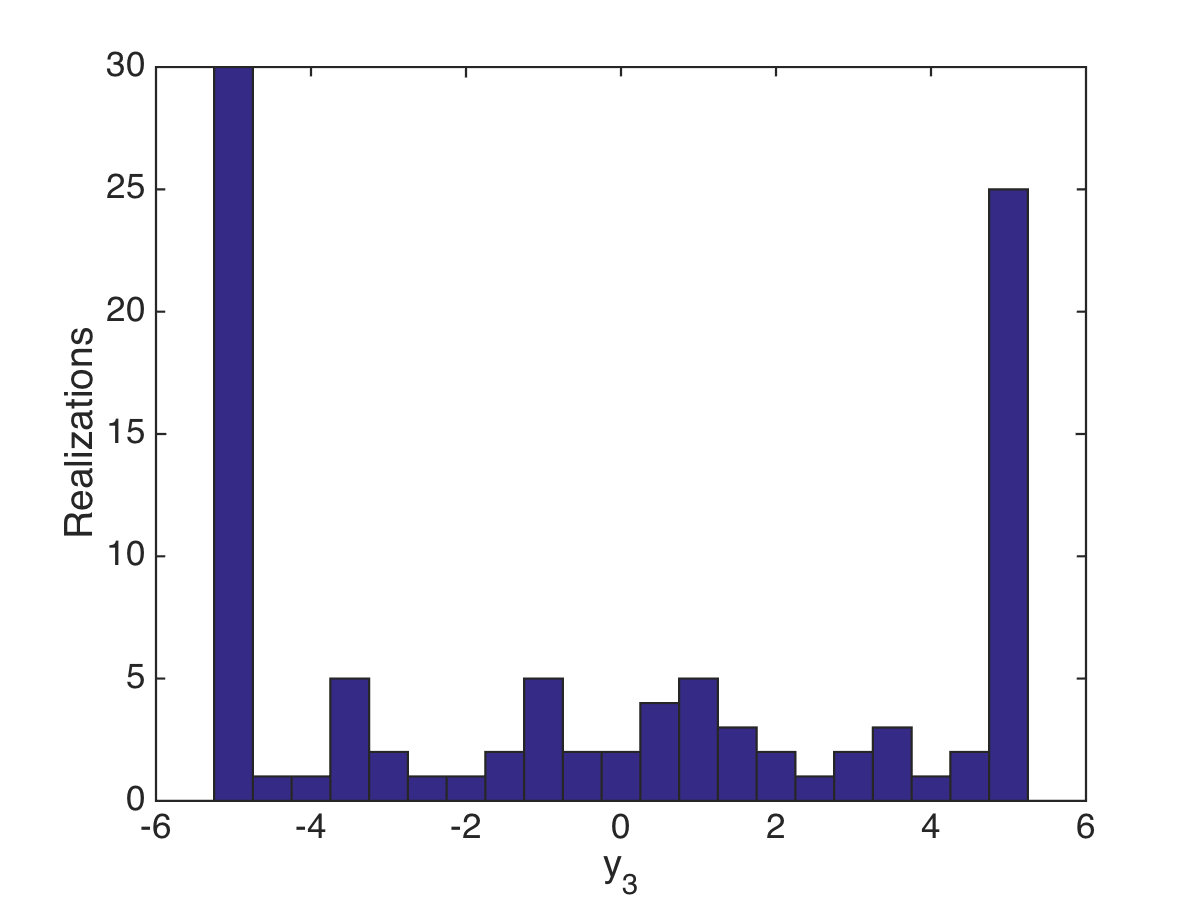}
\includegraphics[width=0.3\textwidth]{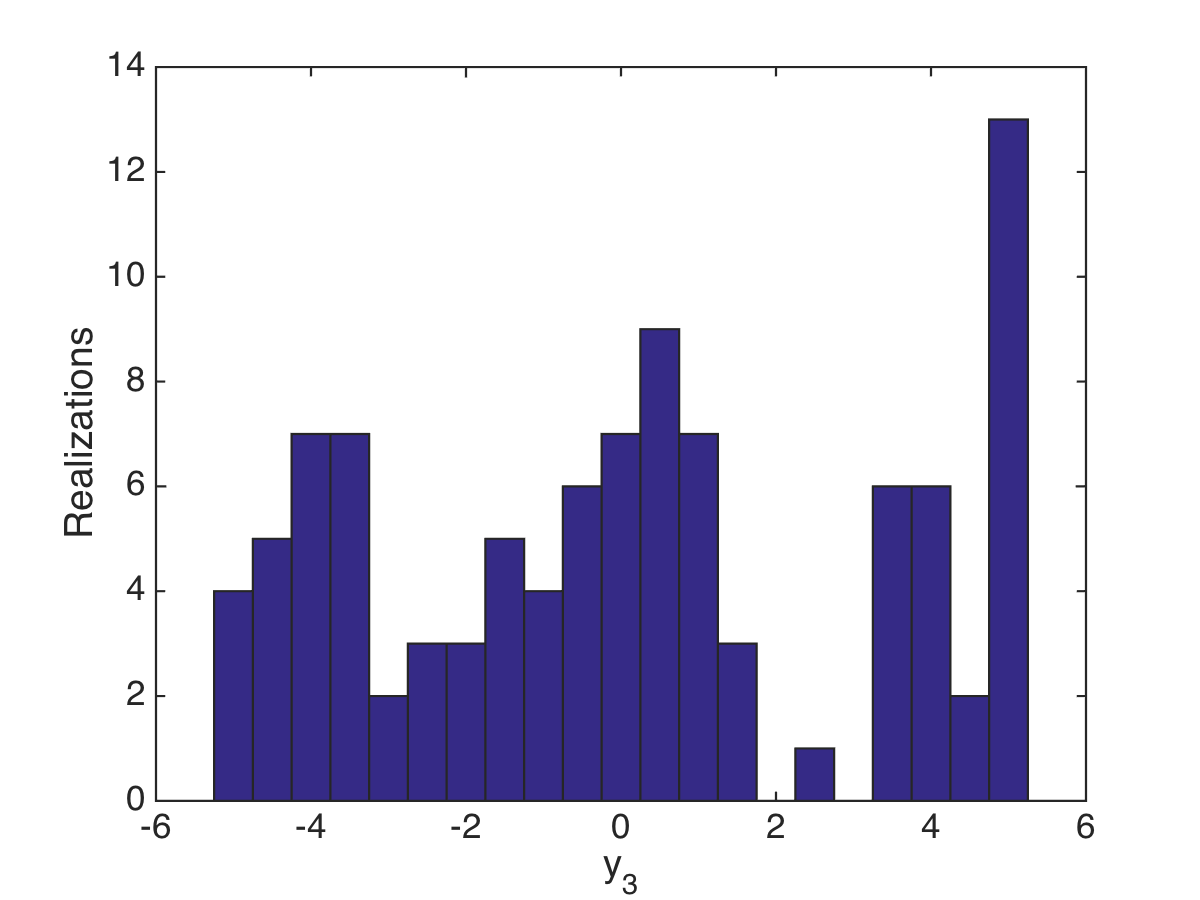}
\end{center}
\vspace{-0.1in}
\caption{The histograms of the peak location in $y_1$ (top), $y_2$
  (middle) and $y_3$ (bottom) of the imaging function \eqref{eq:IN13}
  in the small aperture regime.  The left column is for $\sigma/\sigma_1 =
  50\%$, the middle column for $25\%$ and the right column for
  $10\%$. The abscissa is in units of the wavelength.  }
\label{fig:HIST_FF}
\end{figure}
Imaging is more difficult in the small aperture, as seen in Figure
\ref{fig:HIST_FF}. Although the cross-range localization remains
robust up to $50\%$ noise, and the error is bounded by $\lambda$,
which is an improvement over the resolution limit $\lambda L/a = 10
\lambda$ of migration imaging, the range localization fails even at
$10\%$ noise.  This is because the effective rank of $\tDD$ is $\widetilde{\mathfrak R} = 2$ 
for all the simulations in Figure \ref{fig:HIST_FF}, and as explained in section \ref{sect:disc}, 
the projection matrices on the subspace spanned by
the first two  singular vectors  are approximately
independent of the range of the inclusion. To determine the range we
need an array of larger aperture, or two arrays that view the
inclusion from different directions. Alternatively, we may probe the
medium with broad-band pulses and estimate the range from travel times
of the scattered returns.

\begin{figure}[h]
\begin{center}
\includegraphics[width=0.3\textwidth]{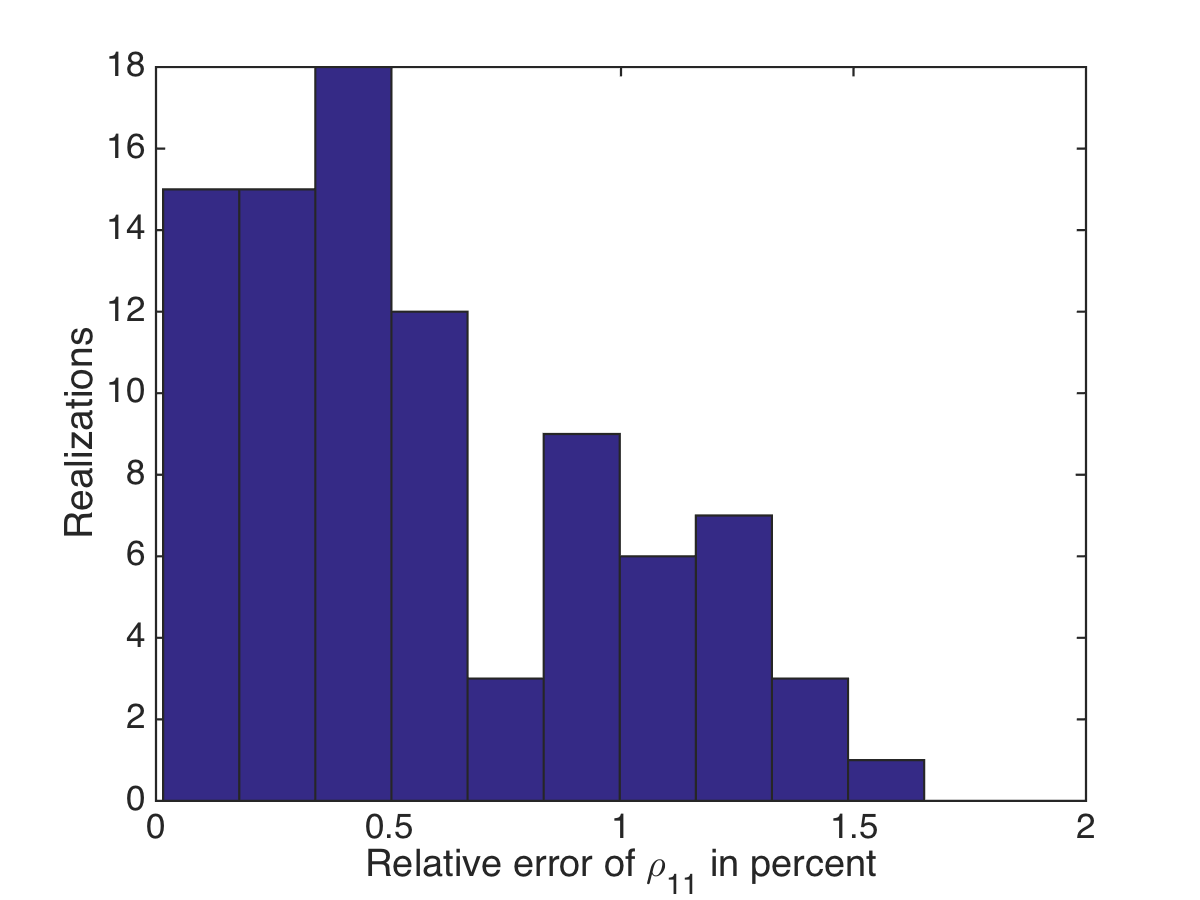}
\includegraphics[width=0.3\textwidth]{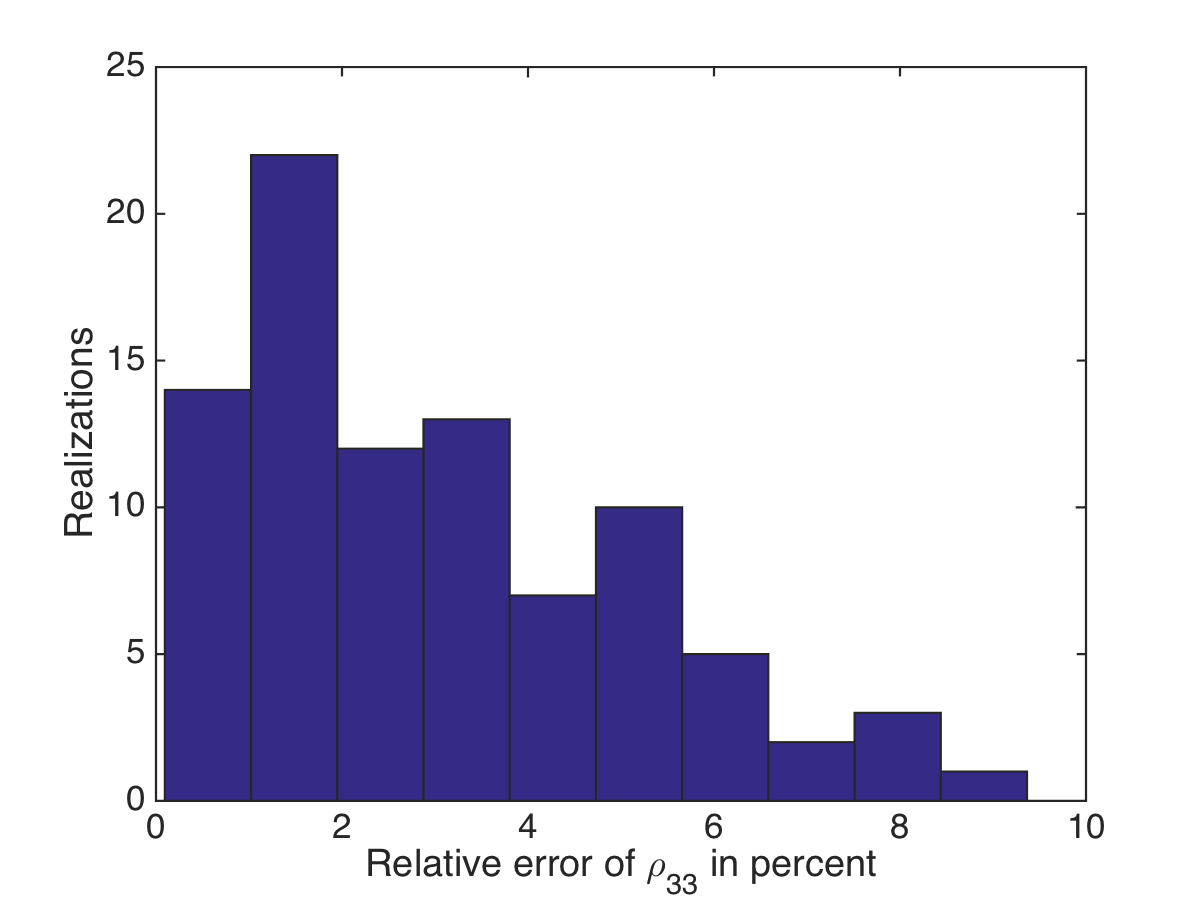}
\includegraphics[width=0.3\textwidth]{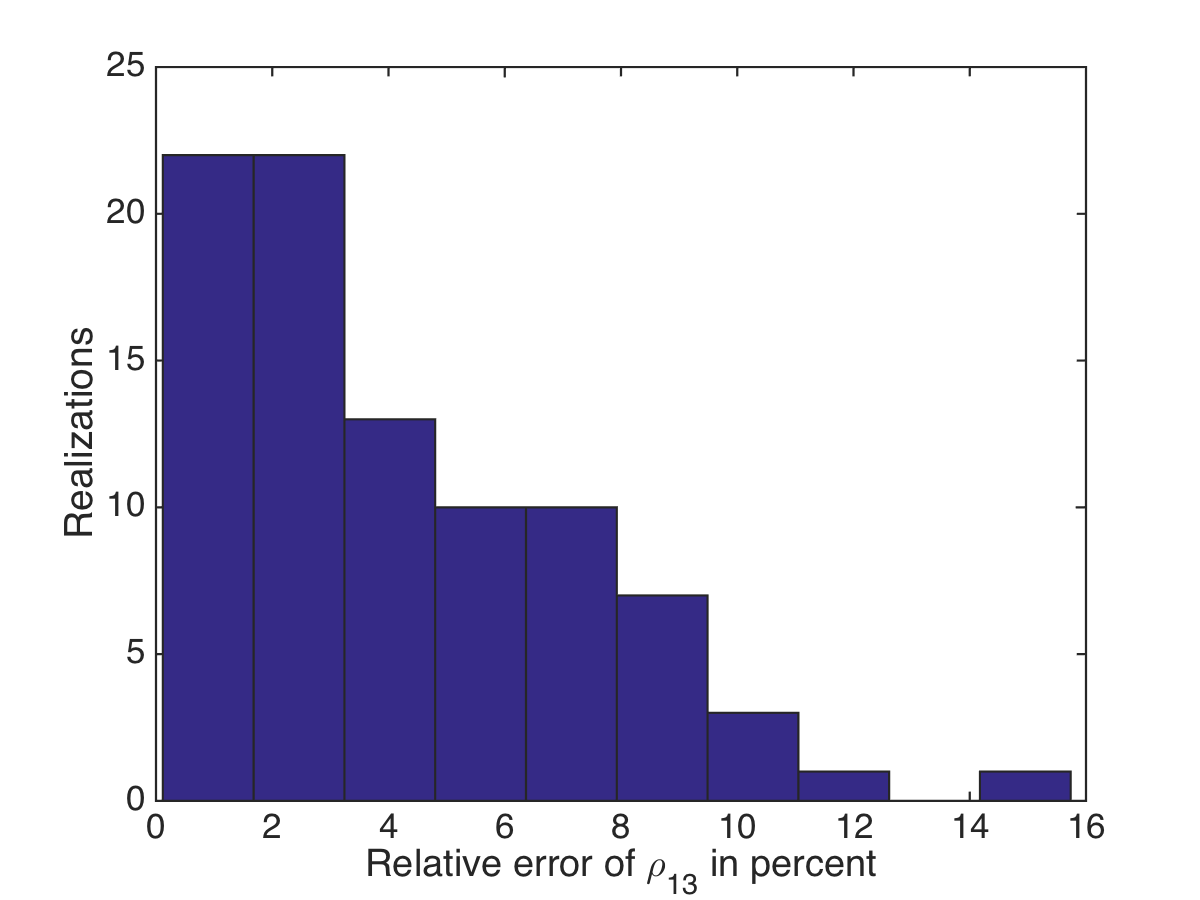} \\
\vspace{-0.01in}\includegraphics[width=0.3\textwidth]{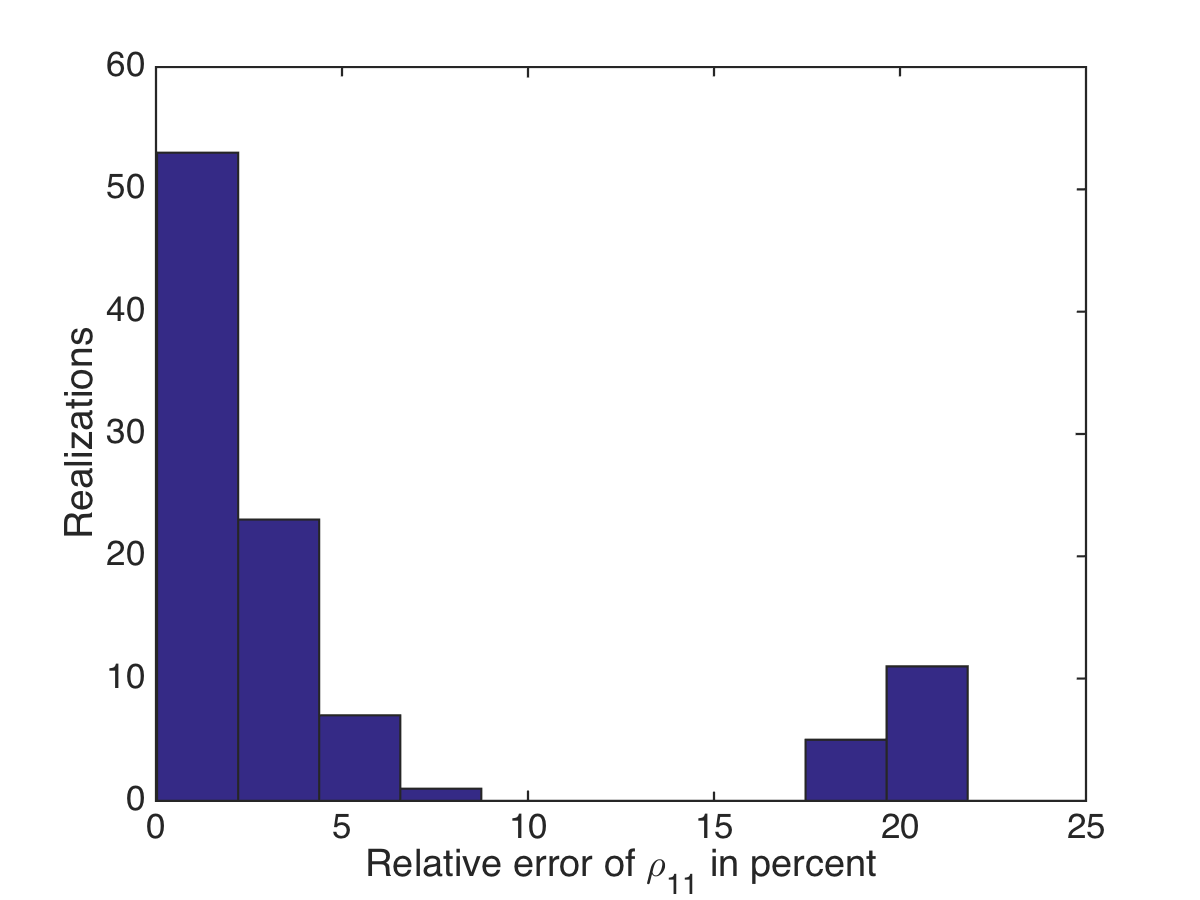}
\includegraphics[width=0.3\textwidth]{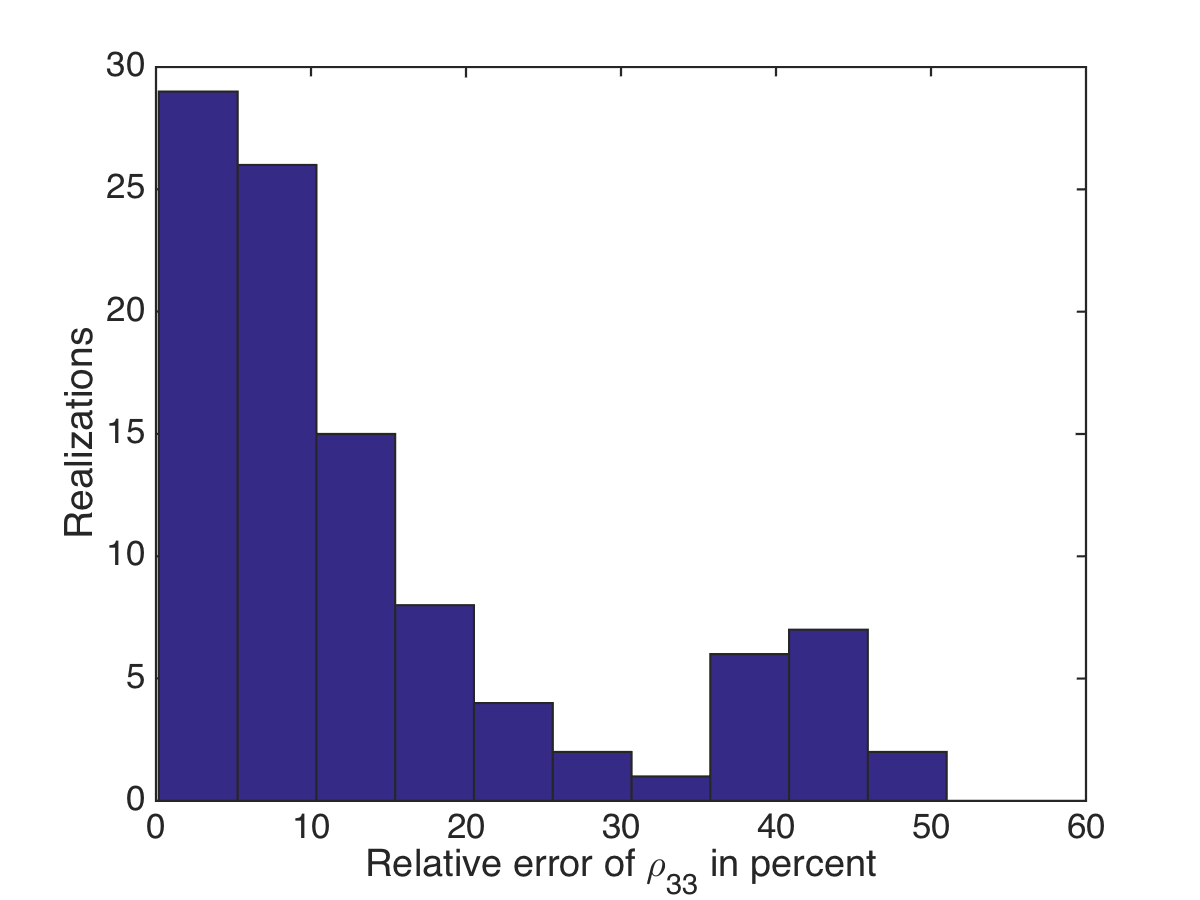}
\includegraphics[width=0.3\textwidth]{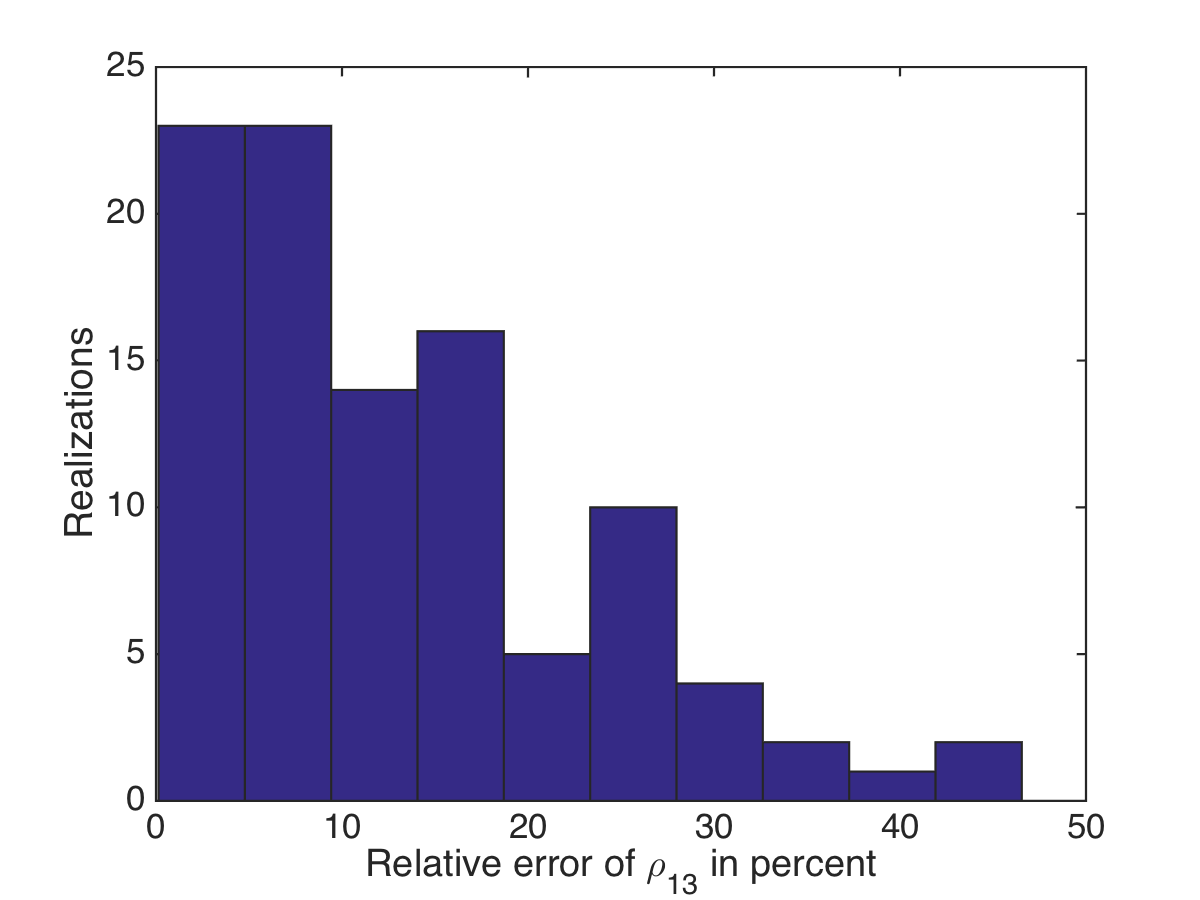} 
\end{center}
\vspace{-0.1in}
\caption{Histograms of the relative errors of the components
  $\rho_{11}$ (left), $\rho_{33}$ (middle) and $\rho_{13}$ (right) of
  the reflectivity tensor $\brho$, for the large aperture regime and complete
  measurements at $25\%$ noise (top row) and $50\%$ noise (bottom
  row). The abscissa is in percent.}
\label{fig:HIST_NF_RHO}
\end{figure}

\begin{figure}[h]
\begin{center}
\includegraphics[width=0.3\textwidth]{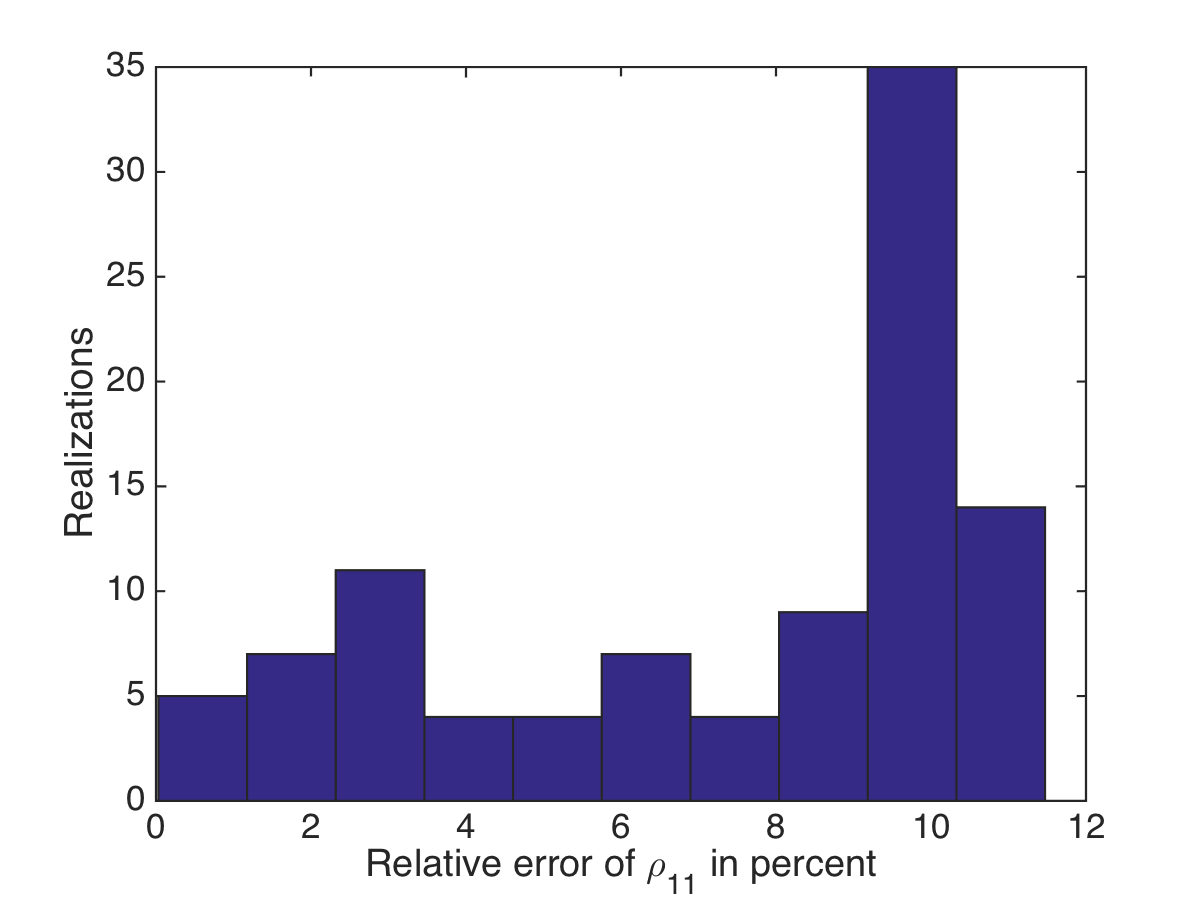}
\includegraphics[width=0.3\textwidth]{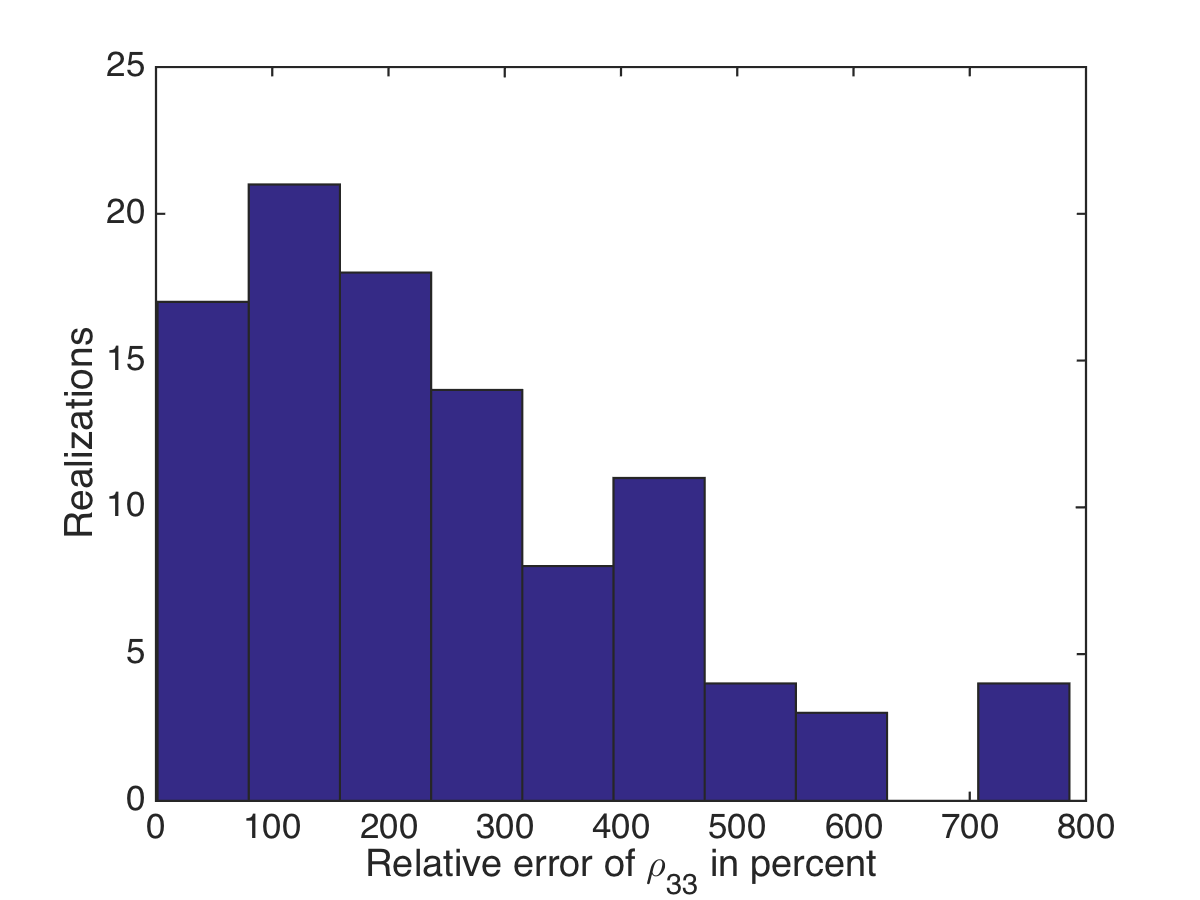} 
\includegraphics[width=0.3\textwidth]{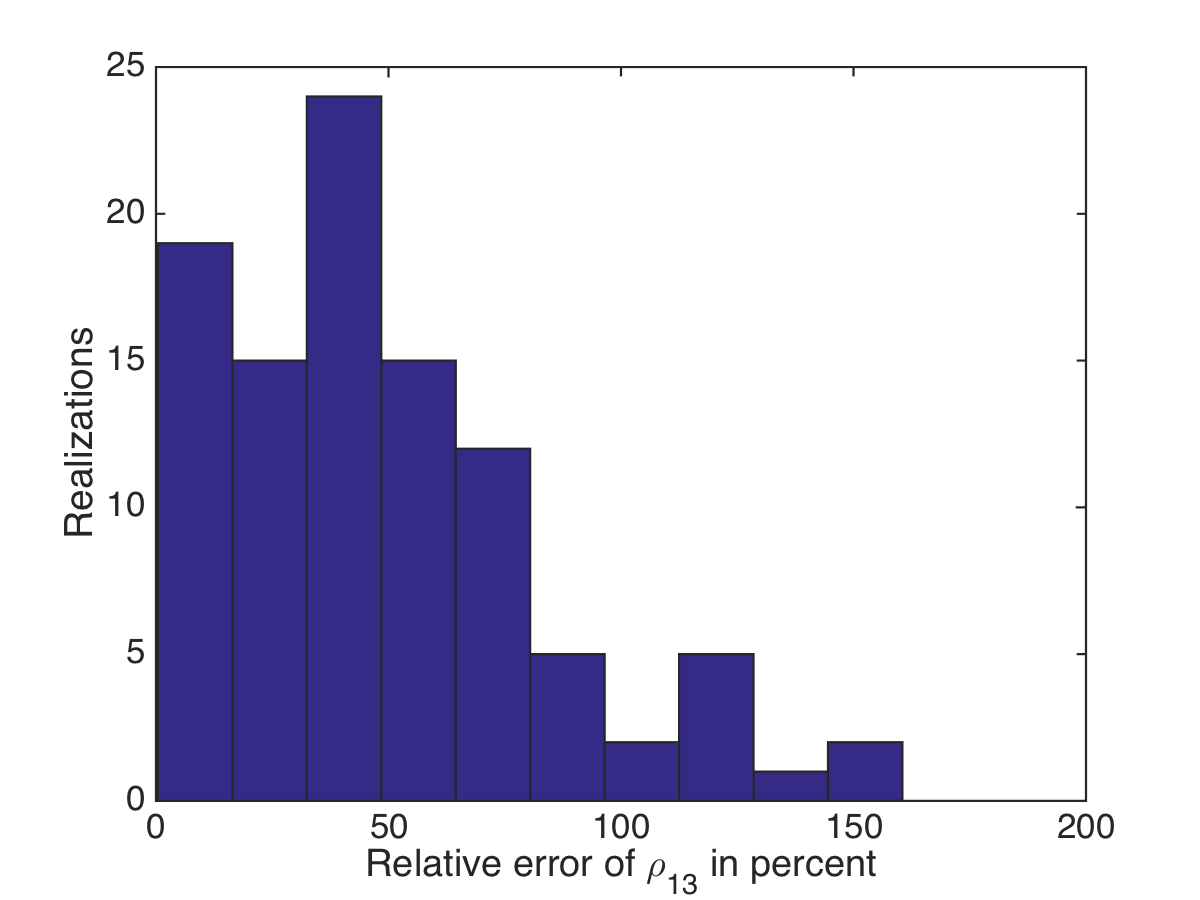}
\end{center}
\vspace{-0.1in}
\caption{Histograms of the relative errors of the components
  $\rho_{11}$ (left), $\rho_{33}$ (middle) and $\rho_{13}$ (right) of
  the reflectivity tensor $\brho$, for the small aperture regime and
  complete measurements at $25\%$ noise. The abscissa is in percent.}
\label{fig:HIST_FF_RHO}
\end{figure}

The histograms of the relative errors in the estimation of the reflectivity tensor
are displayed in Figures \ref{fig:HIST_NF_RHO}, for the large aperture regime
and at $25\%$ and $50\%$ noise levels. The estimates are obtained
using equations \eqref{eq:IN22} and \eqref{eq:IN23}, and the peak
location of the imaging function \eqref{eq:IN13}. The plots in the
left column are for $\rho_{11}$, in the middle column for $\rho_{33}$,
and in the right column for $\rho_{13}$. The results for $\rho_{22}$
and $\rho_{12}$ are similar to those for $\rho_{11}$, and the results
for $\rho_{23}$ are similar to those for $\rho_{13}$.  We note that
the errors are below $10\%$ at $25\%$ noise and naturally, they increase with
the noise level. As expected, the errors are larger for  $\rho_{33}$ and
$\rho_{13}$, but not by a big factor. This is because the condition
number of the search matrix  $\GG(\vy)$, defined as the ratio of its
largest and smallest singular values, equals $2.6$ in the large aperture
regime. The estimates are worse in the small aperture regime, where the
condition number of $\GG(\vy)$ is $23.3$, as shown in Figure
\ref{fig:HIST_FF_RHO}. Here we also have errors due to the poor
estimates of the range component of the inclusion location.

\subsubsection{Incomplete measurements}
If the sensors measure a single component of the scattered electric
field, say along $\ve_q$, the best choice is for $q = 1$ or $2$. The
inversion results are very poor when $q = 3$, even at low levels of
noise, as expected from the discussion in section \ref{sect:inv.2}. We
show here results for the sensing matrix $\bS = (\ve_1)$, in the
large aperture regime, and compare them with those in Figures
\ref{fig:HIST_NF} and \ref{fig:HIST_NF_RHO}. We do not show results in
the small aperture regime because they do not add more information, and
they are, as expected, worse than in the large aperture regime.

\begin{figure}[h]
\begin{center}
\includegraphics[width=0.3\textwidth]{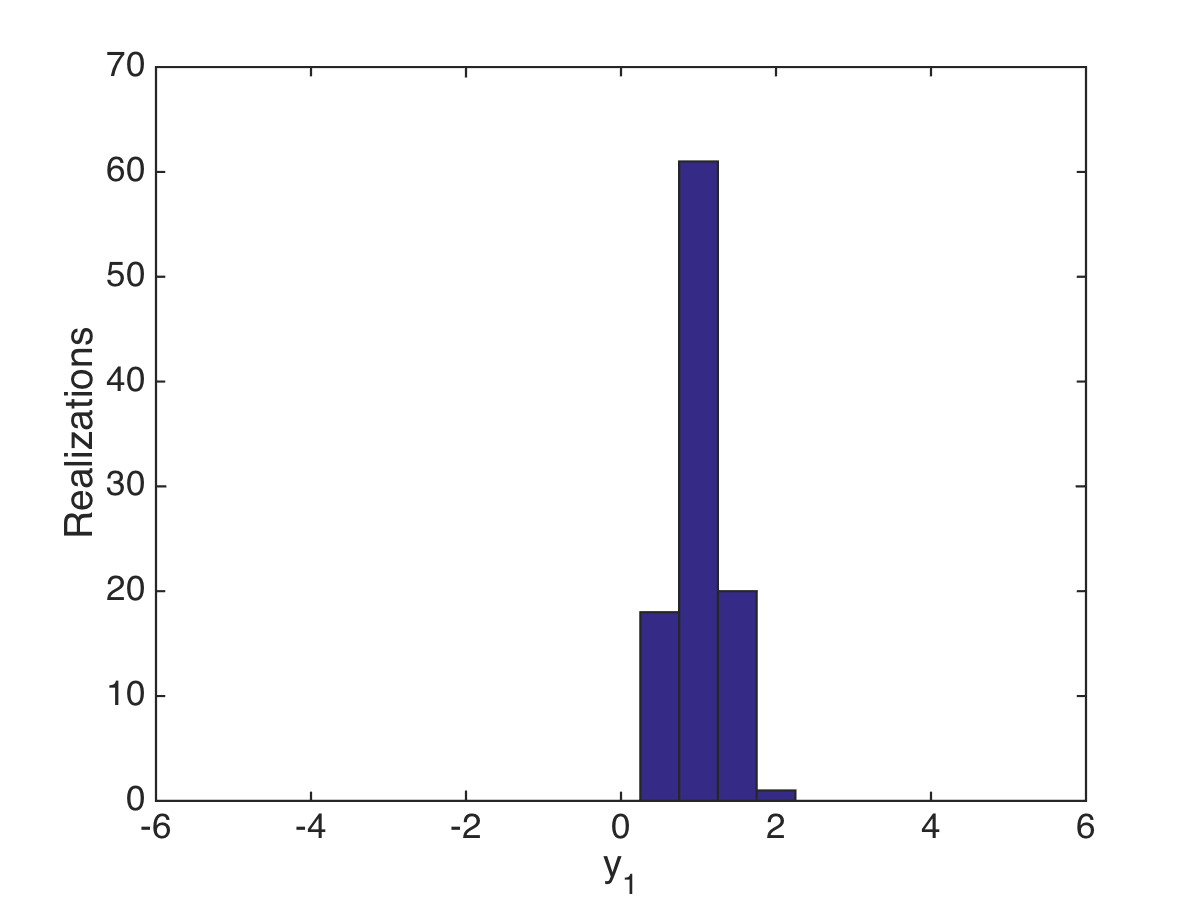}
\includegraphics[width=0.3\textwidth]{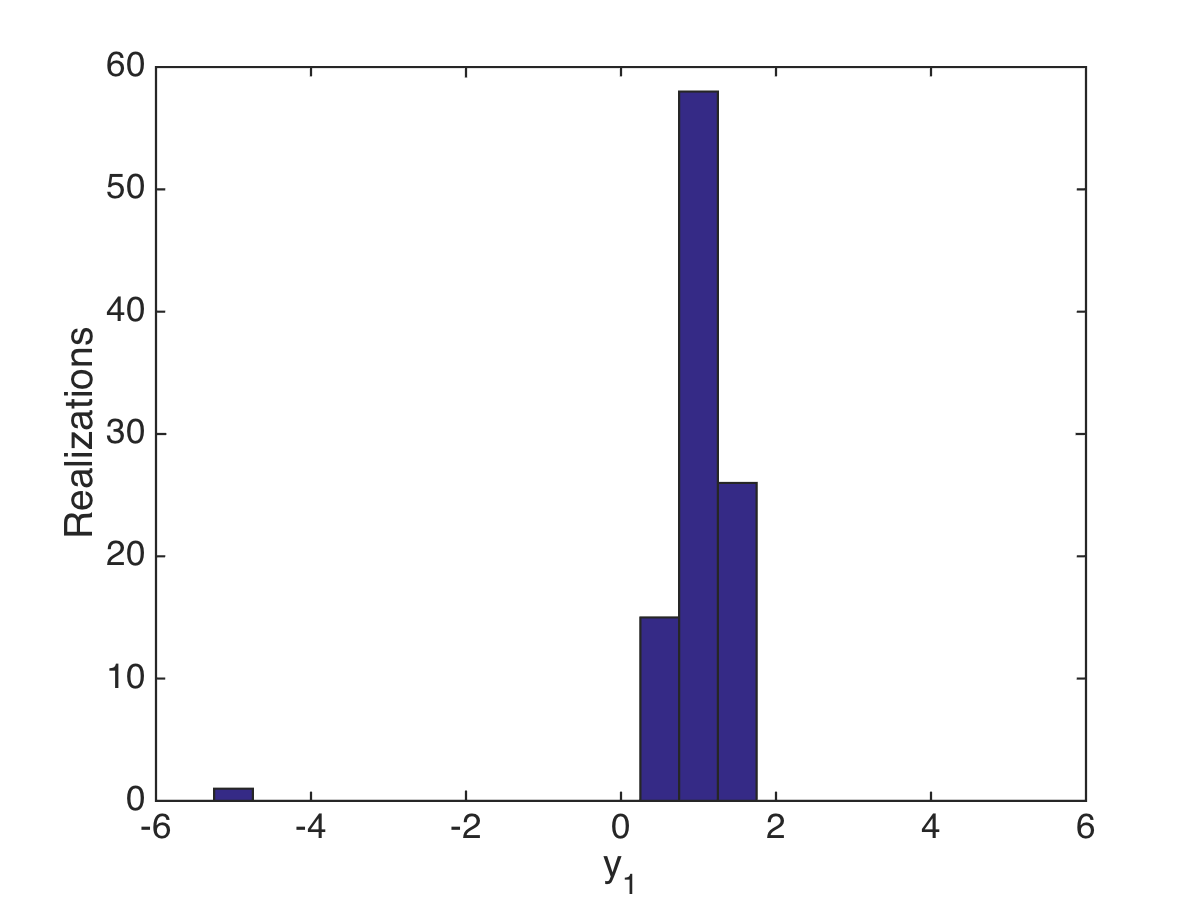}
\includegraphics[width=0.3\textwidth]{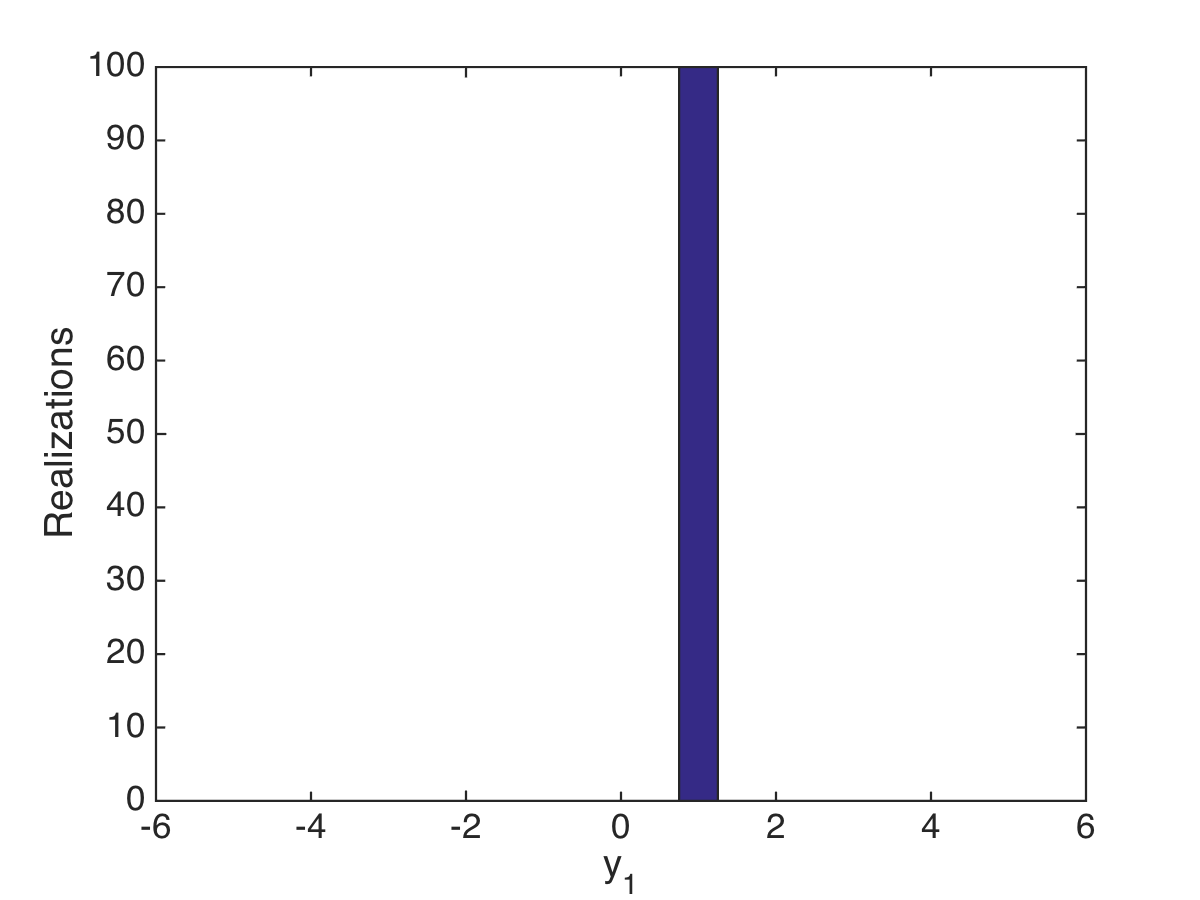} \\
\vspace{-0.02in}\includegraphics[width=0.3\textwidth]{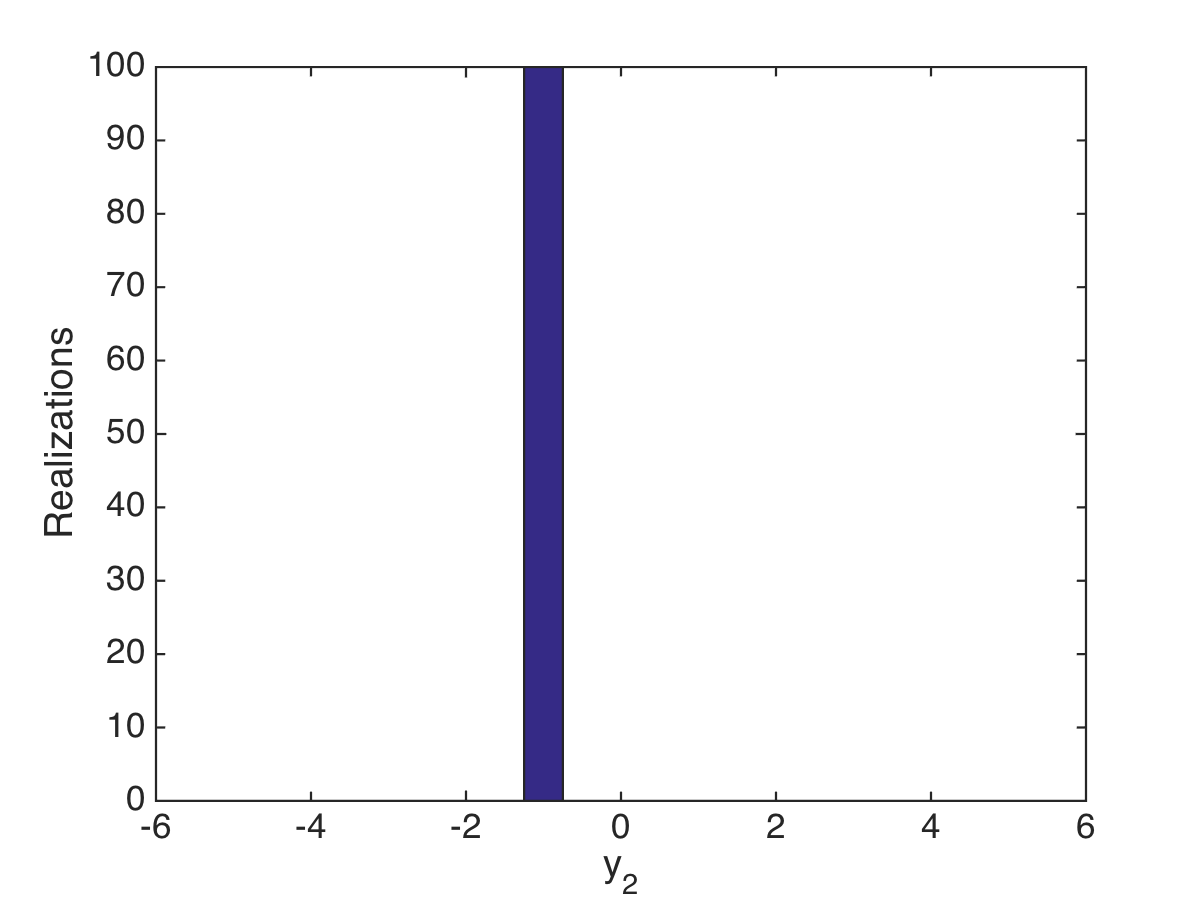}
\includegraphics[width=0.3\textwidth]{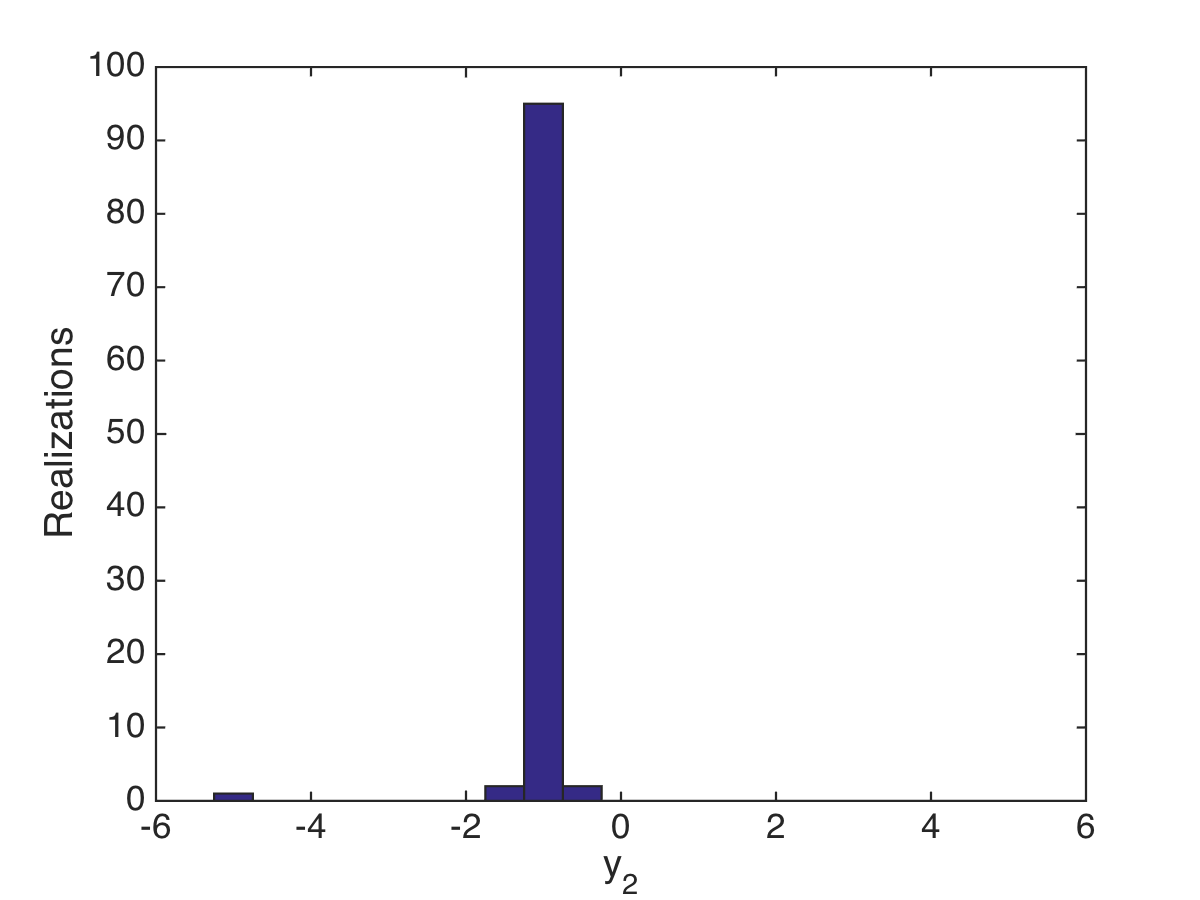}
\includegraphics[width=0.3\textwidth]{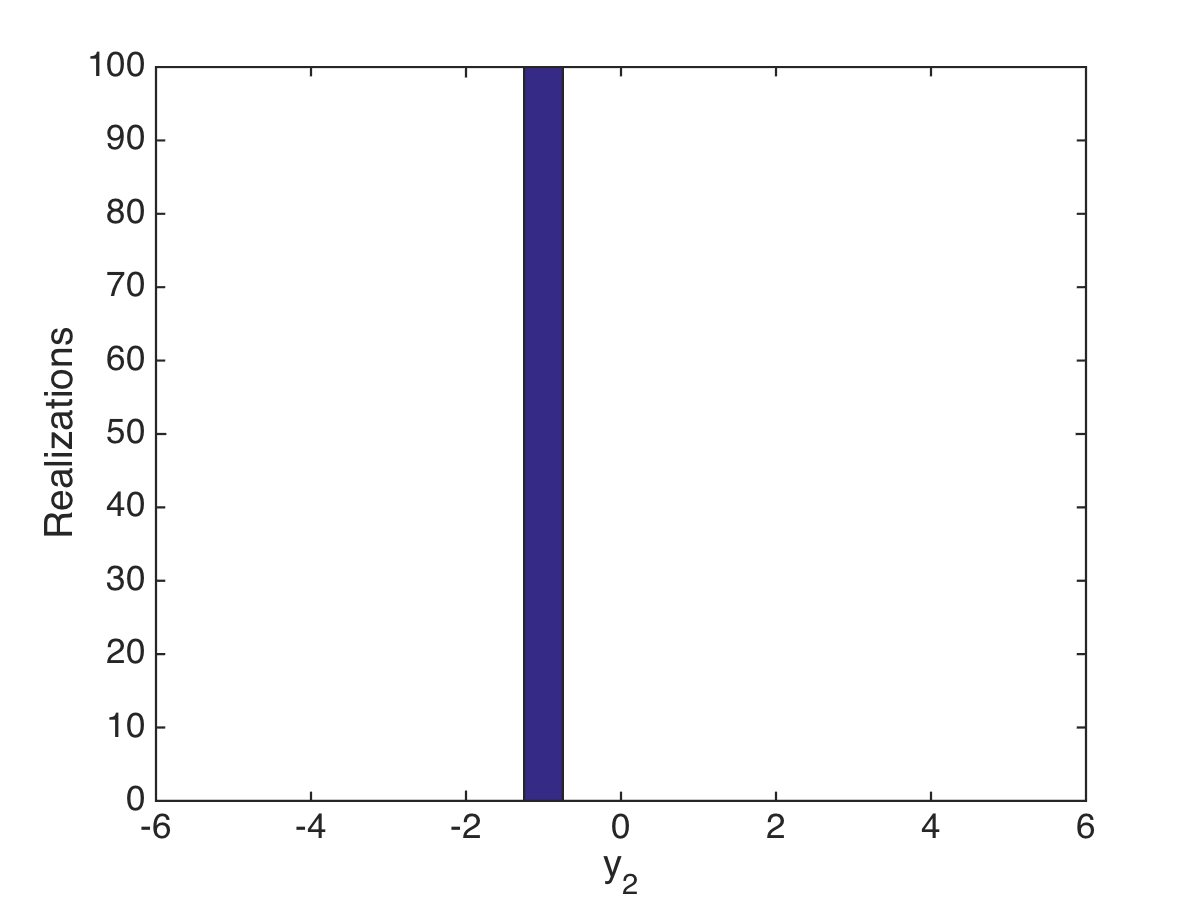} \\
\vspace{-0.02in}\includegraphics[width=0.3\textwidth]{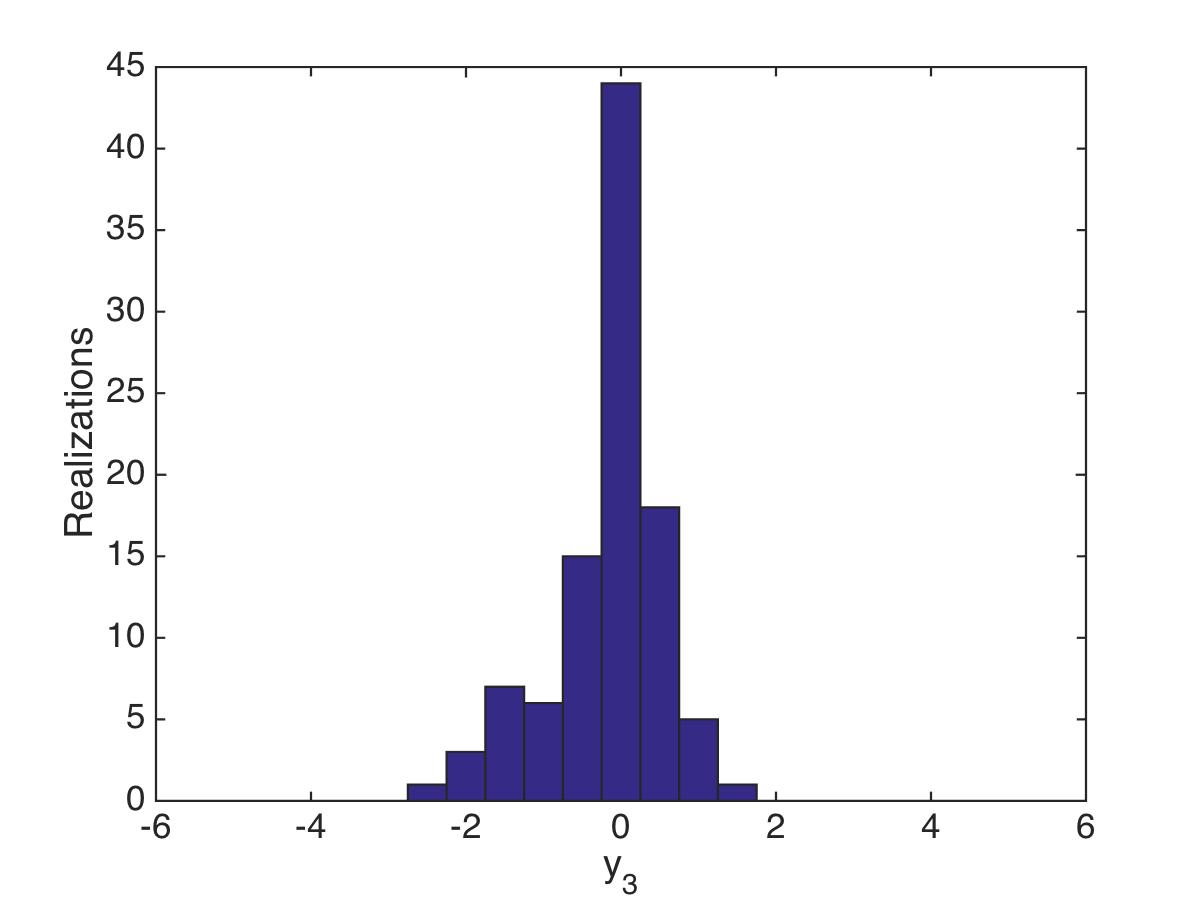}
\includegraphics[width=0.3\textwidth]{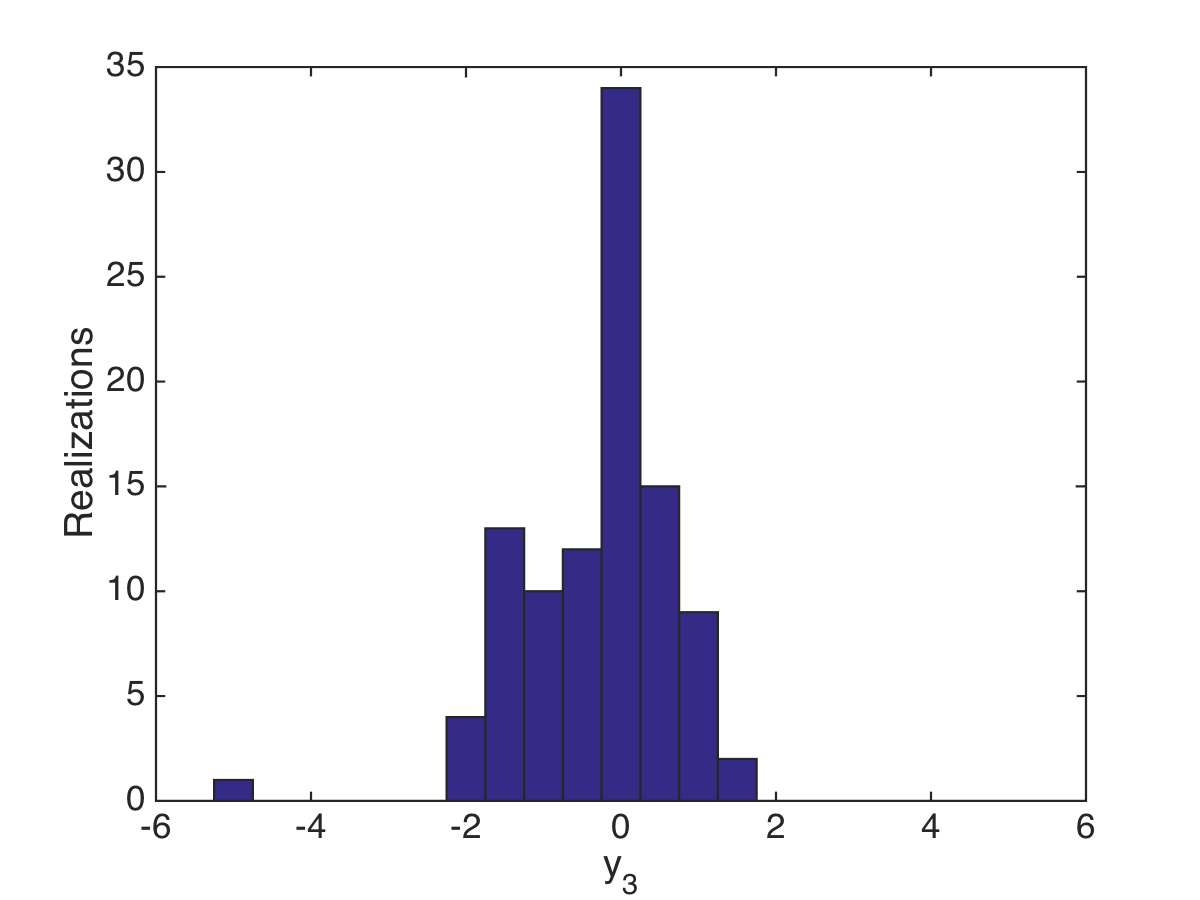}
\includegraphics[width=0.3\textwidth]{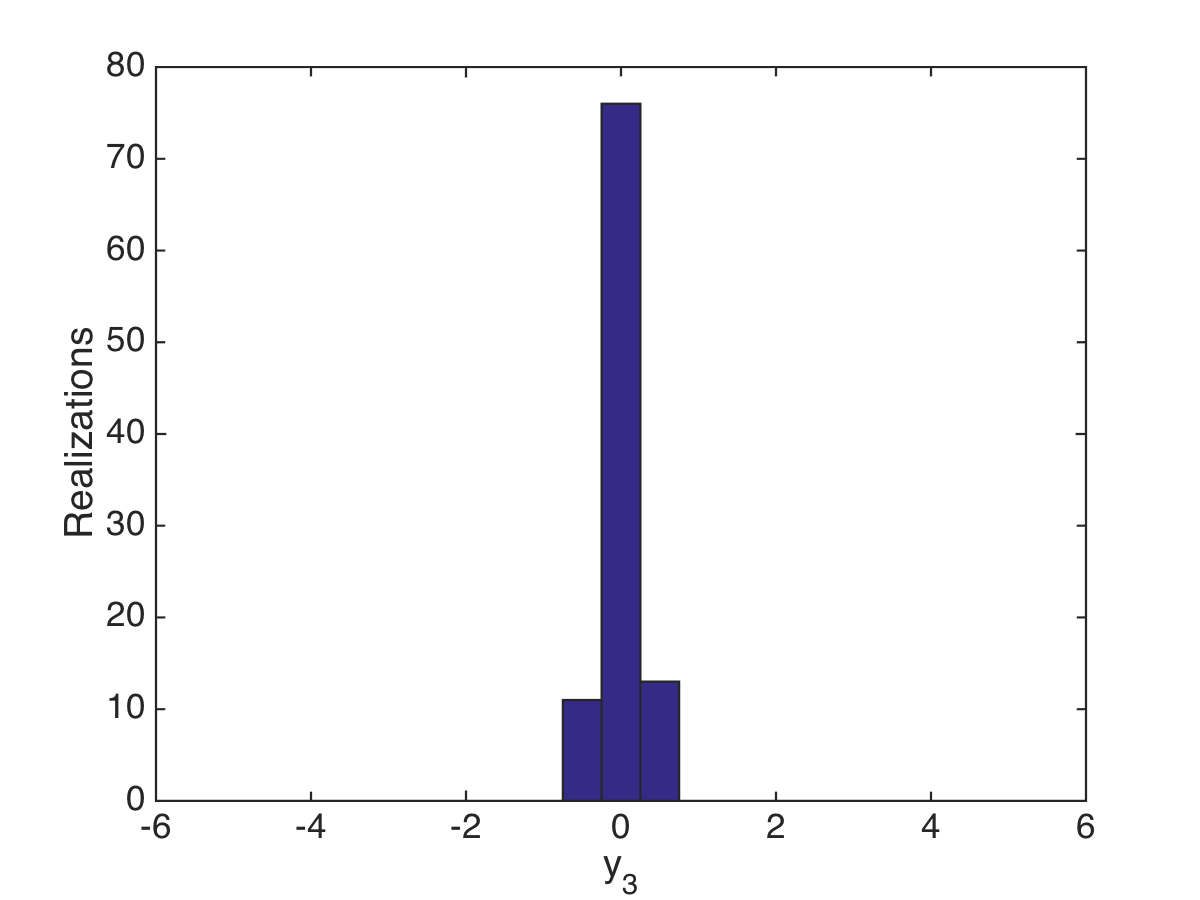}
\end{center}
\vspace{-0.1in}
\caption{The histograms of the peak location in $y_1$ (top), $y_2$
  (middle) and $y_3$ (bottom) of the imaging function \eqref{eq:IN13},
  in the large aperture regime, for incomplete measuremenst with $\bS =
  (\ve_1)$.  The left column is for $\sigma/\sigma_1 = 75\%$, the
  middle column for $50\%$ and the right column for $25\%$. The
  abscissa is in units of the wavelength.}
\label{fig:HIST_NF_PARTU1}
\end{figure}

\begin{figure}[t]
\begin{center}
\includegraphics[width=0.3\textwidth]{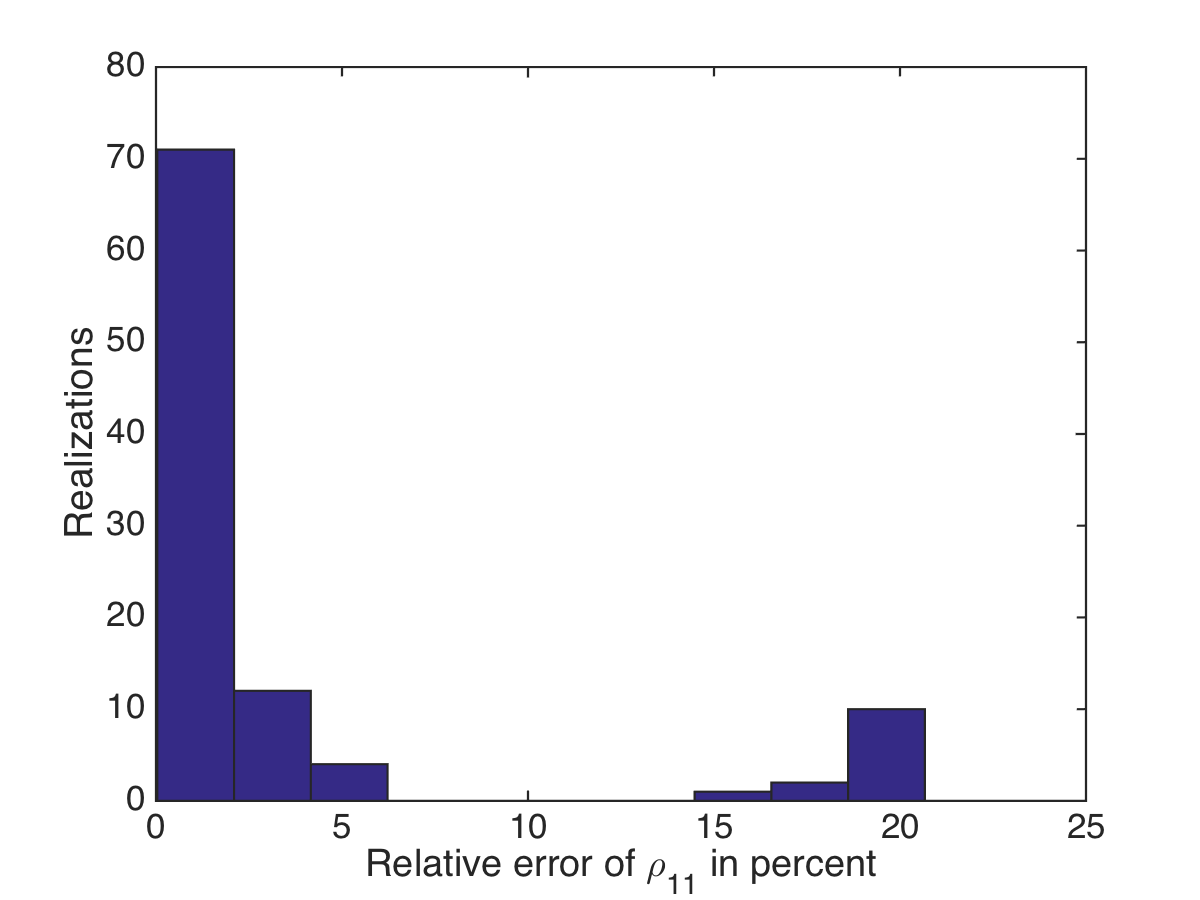}
\includegraphics[width=0.3\textwidth]{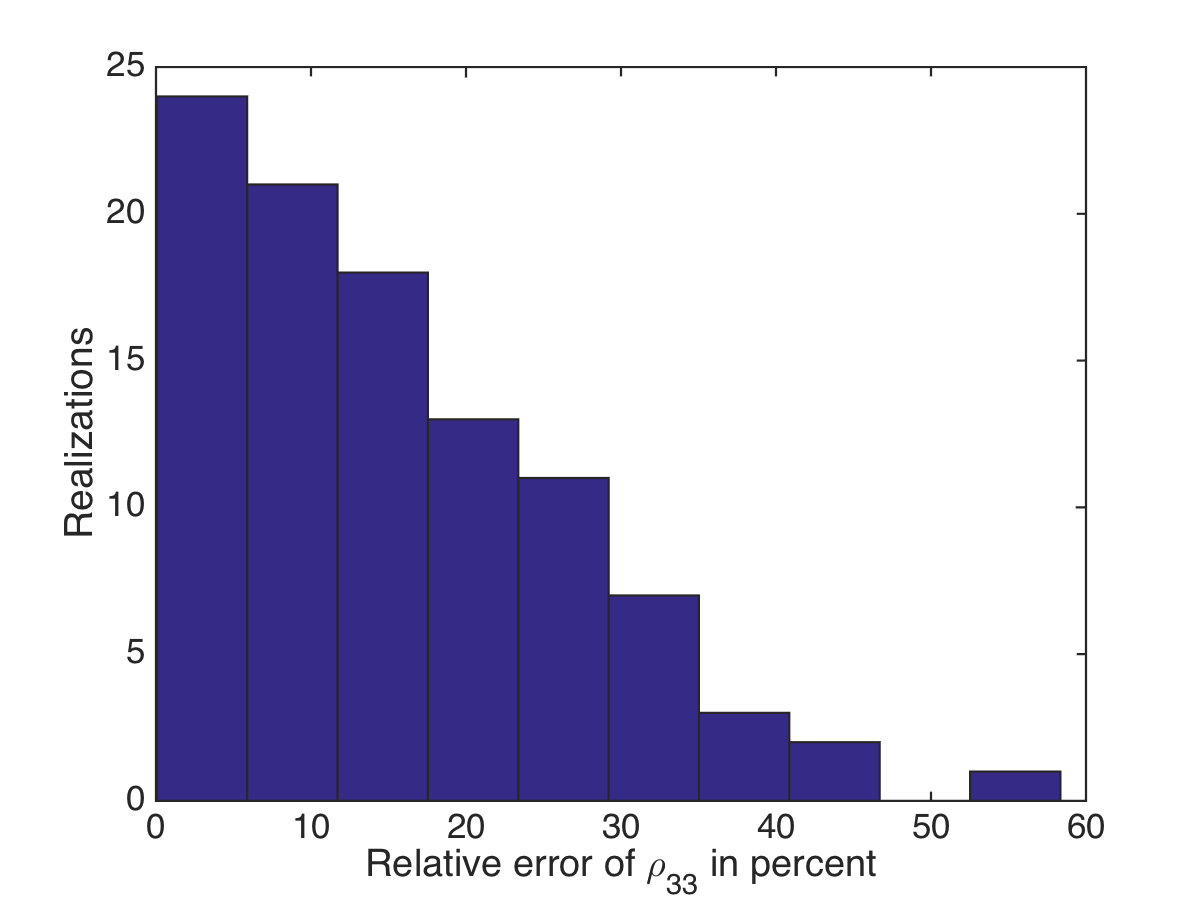} 
\includegraphics[width=0.3\textwidth]{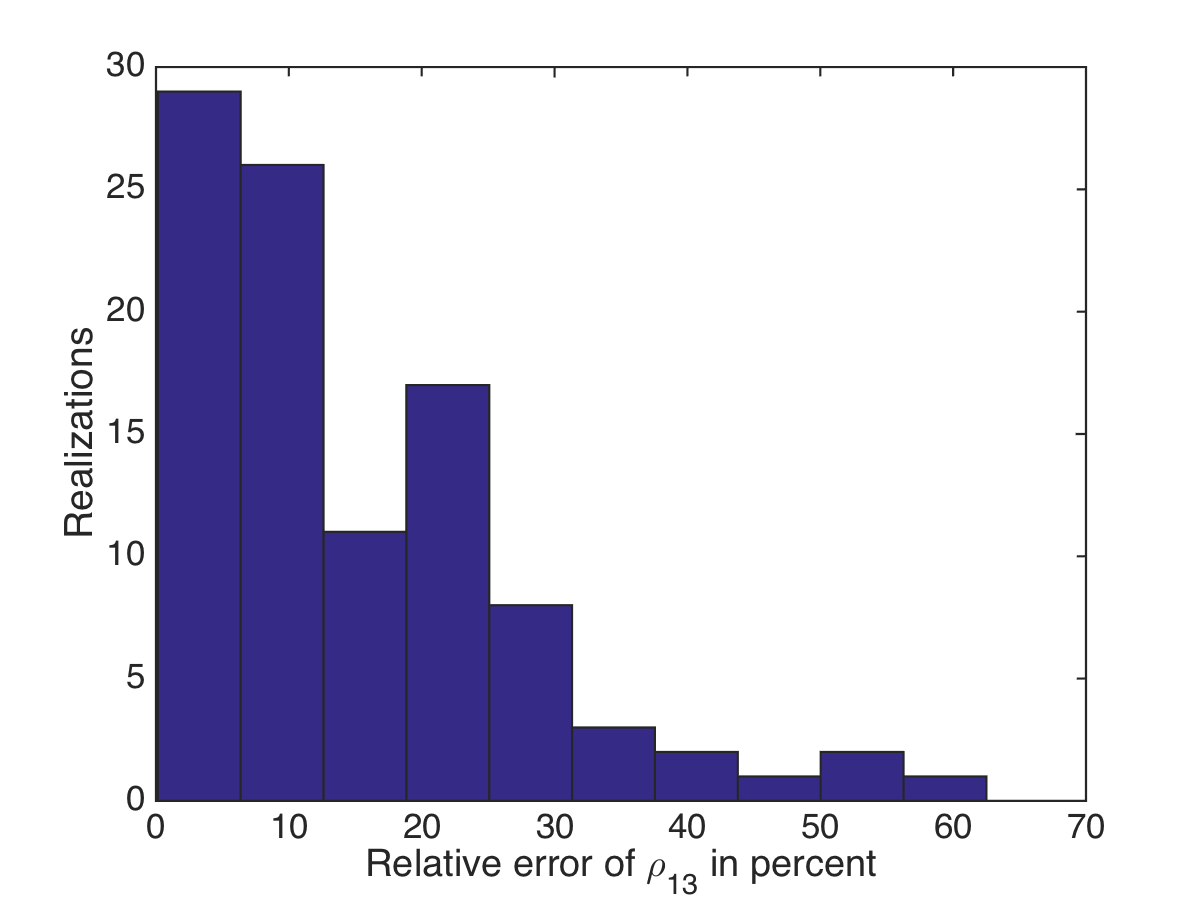}
\end{center}
\vspace{-0.1in}
\caption{Histograms of the relative errors of the components
  $\rho_{11}$ (left), $\rho_{33}$ (middle) and $\rho_{13}$ (right) of
  the reflectivity tensor $\brho$, in the large aperture regime, for
  incomplete measuremenst with $\bS = (\ve_1)$ and $25\%$ noise. The
  abscissa is in percent.}
\label{fig:HIST_NFRHO_PARTU1}
\end{figure}

The histogram plots in Figure \ref{fig:HIST_NF_PARTU1} show that the
inclusion localization is worse than in the complete measurement case
(Figure \ref{fig:HIST_NF}), but not dramatically so. However, the
errors in the estimation of the reflectivity tensor are much larger, even at the $25\%$ noise
level. Compare the results in Figure \ref{fig:HIST_NFRHO_PARTU1} to
those on the top row of Figure \ref{fig:HIST_NF_RHO}.  This demonstrates the
benefit of measuring more than one component of the scattered
electric field.

We do not display inversion results from incomplete data with the sensing
matrix $\bS = (\ve_1,\ve_2)$, because they are comparable to those in
the complete measurements case. This is expected from the discussion
in section \ref{sect:inv.2}.

\begin{figure}[h]
\begin{center}
\includegraphics[width=0.52\textwidth]{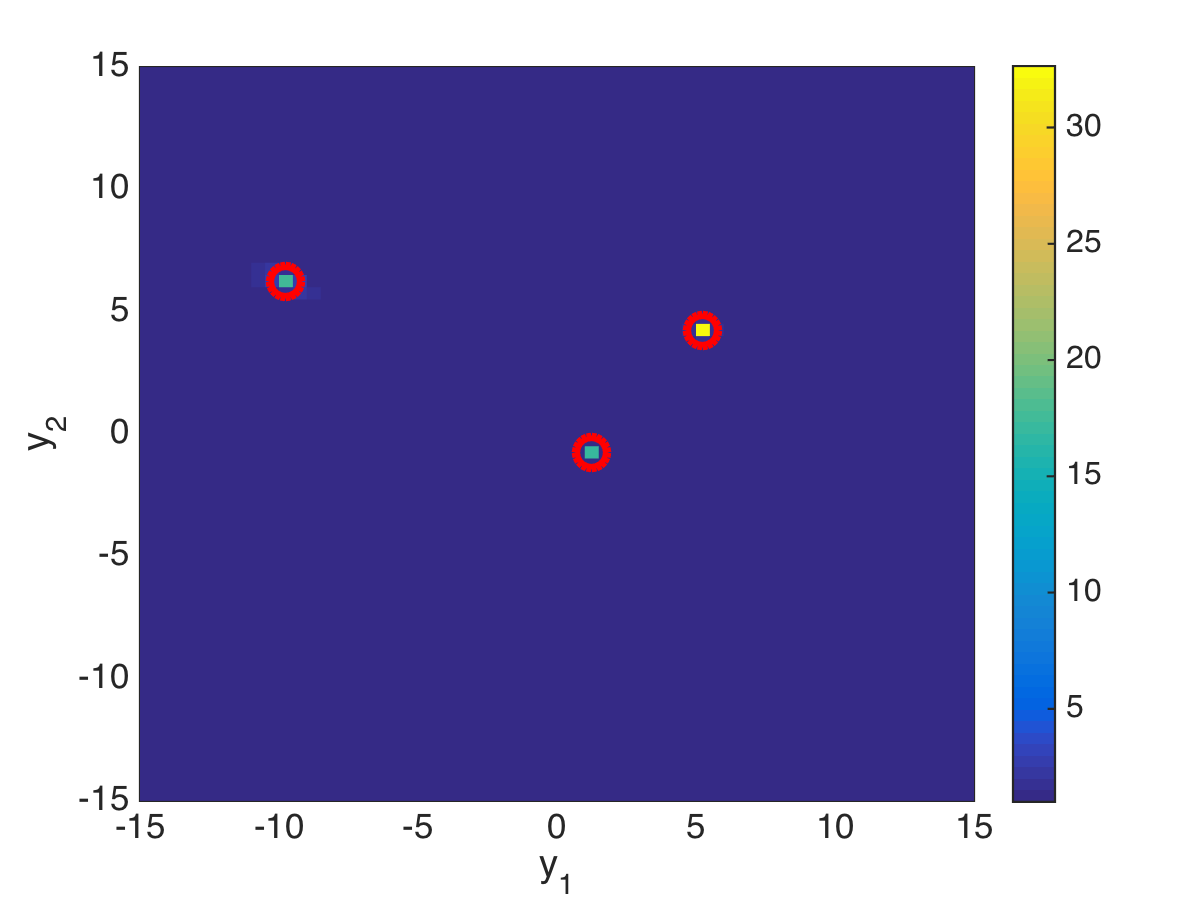}
\hspace{-0.25in}\includegraphics[width=0.52\textwidth]{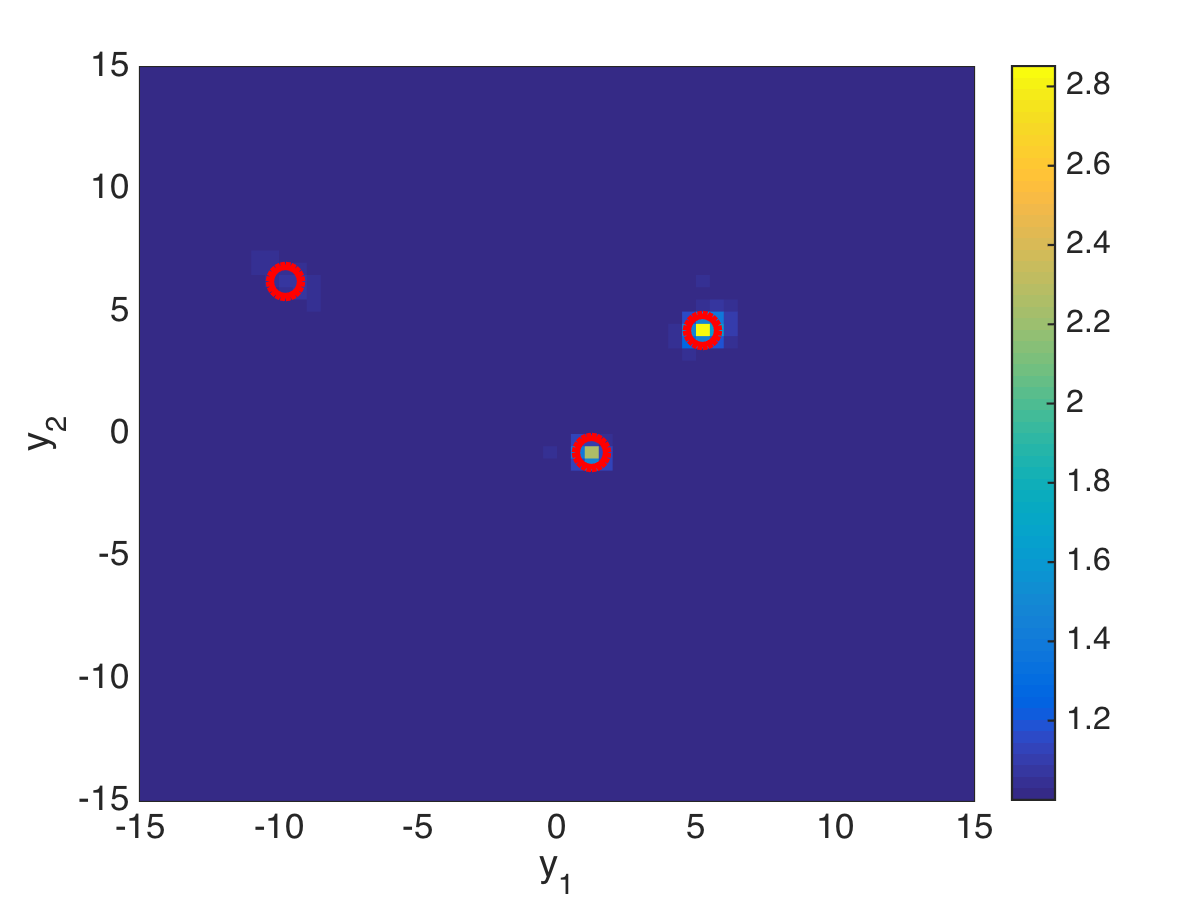}
\end{center}
\vspace{-0.1in}
\caption{Image of three inclusions in the large aperture regime, for complete measurements and $25\%$ noise (left) and 
$75\%$ noise (right). The display is in the plane $y_3 = L$. The axes are in units of the wavelength. The locations of the inclusions
are indicated with red circles.}
\label{fig:3T_NEAR_XY}
\end{figure}
\subsection{Inversion results for multiple inclusions}
\label{sect:num.3}
We present  results for three inclusions  at 
$\vy_1 = (\la,-\la,L)$, $\vy_2 = (-10 \la, 6 \la, L)$ and $\vy_3  = (5 \la, 4 \la, L)$. They are localized using the 
imaging function \eqref{eq:IN19}, and their  reflectivities  are estimated as in \eqref{eq:IN23}. The first inclusion has 
the reflectivity  \eqref{eq:rhoT1} and for the other two we have
\begin{equation}
\alpha^{-3} \brho_2 = \left(\begin{matrix} 91.12 & -29.43& -12.04 \\
-29.43& 69.99& - 8.53\\
-12.04 & -8.53 & 111.60
\end{matrix} \right).
\label{eq:rhoT2}
\end{equation}
and 
\begin{equation}
\alpha^{-3} \brho_3 = \left(\begin{matrix} 99.37 & 5.35 & -26.89 \\
5.35& 101.12 & -10.37 \\
-26.89 & -10.37 & 137.06
\end{matrix} \right).
\label{eq:rhoT3}
\end{equation}
In the  large aperture regime the   matrix  $\DD$ has rank $\mathfrak{R} = 9$, and singular values
\[
1.76 = \sigma_1 > \sigma_2 > \ldots > \sigma_5 = 1.14, ~~ 0.57 = \sigma_6 > \sigma_7 >\sigma_8 = 0.18, ~ ~ \sigma_9 = 0.07.
\]
The rank of $\DD$ is still nine in the small aperture case, but the matrix has three very small singular values that are
indistinguishable from noise at levels as small as $1\%$.
\[
0.034 = \sigma_1 > \sigma_2 > \ldots > \sigma_6 = 0.012, ~~ 6.39 \cdot 10^{-5}= \sigma_7 >  \sigma_8 > \sigma_9 = 2.12\cdot 10^{-5}.
\]

For brevity we  present only results in the complete measurement case,  and a few realizations of the noise. The histograms of the estimated peak locations and 
of the errors in the estimation of the reflectivities are similar to those in the previous section. 

\begin{figure}[t]
\begin{center}
\includegraphics[width=0.34\textwidth]{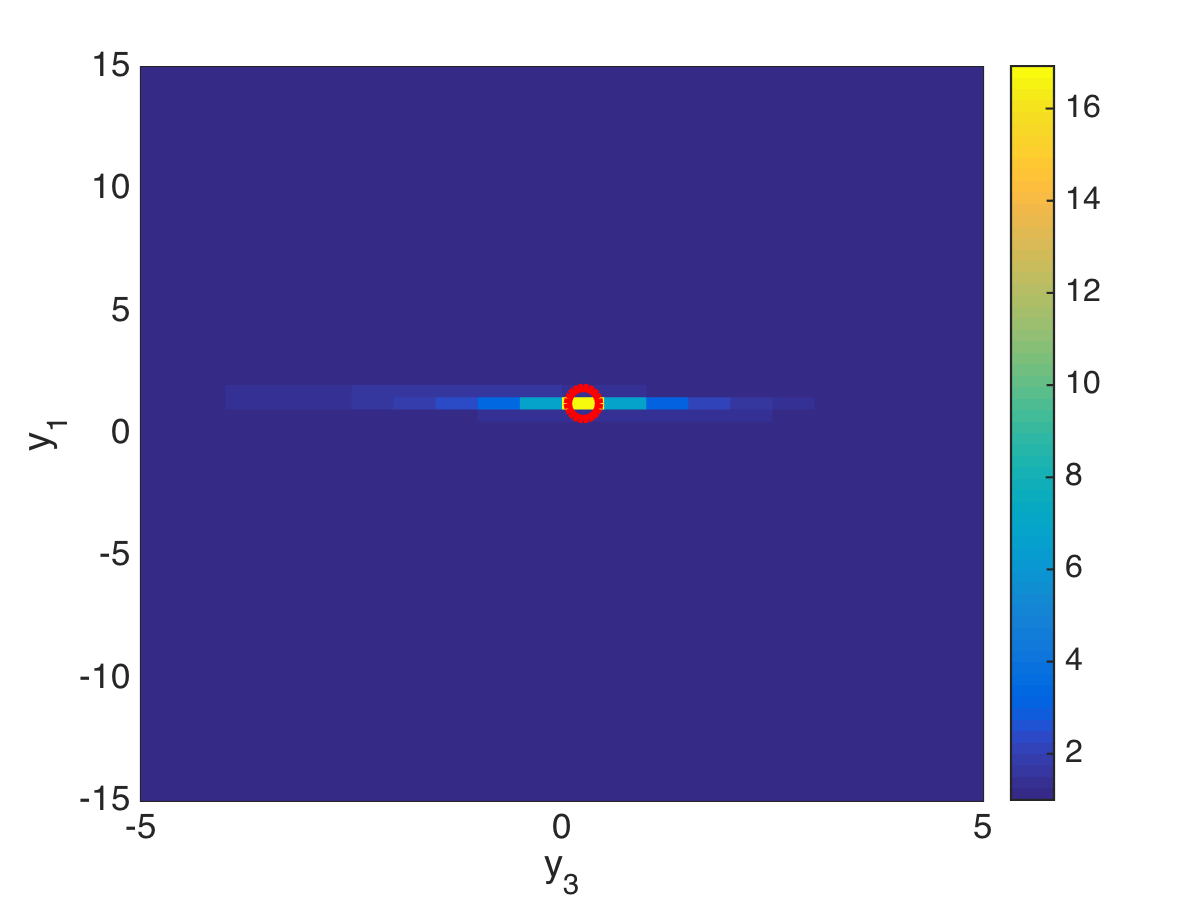}
\hspace{-0.2in}
\includegraphics[width=0.34\textwidth]{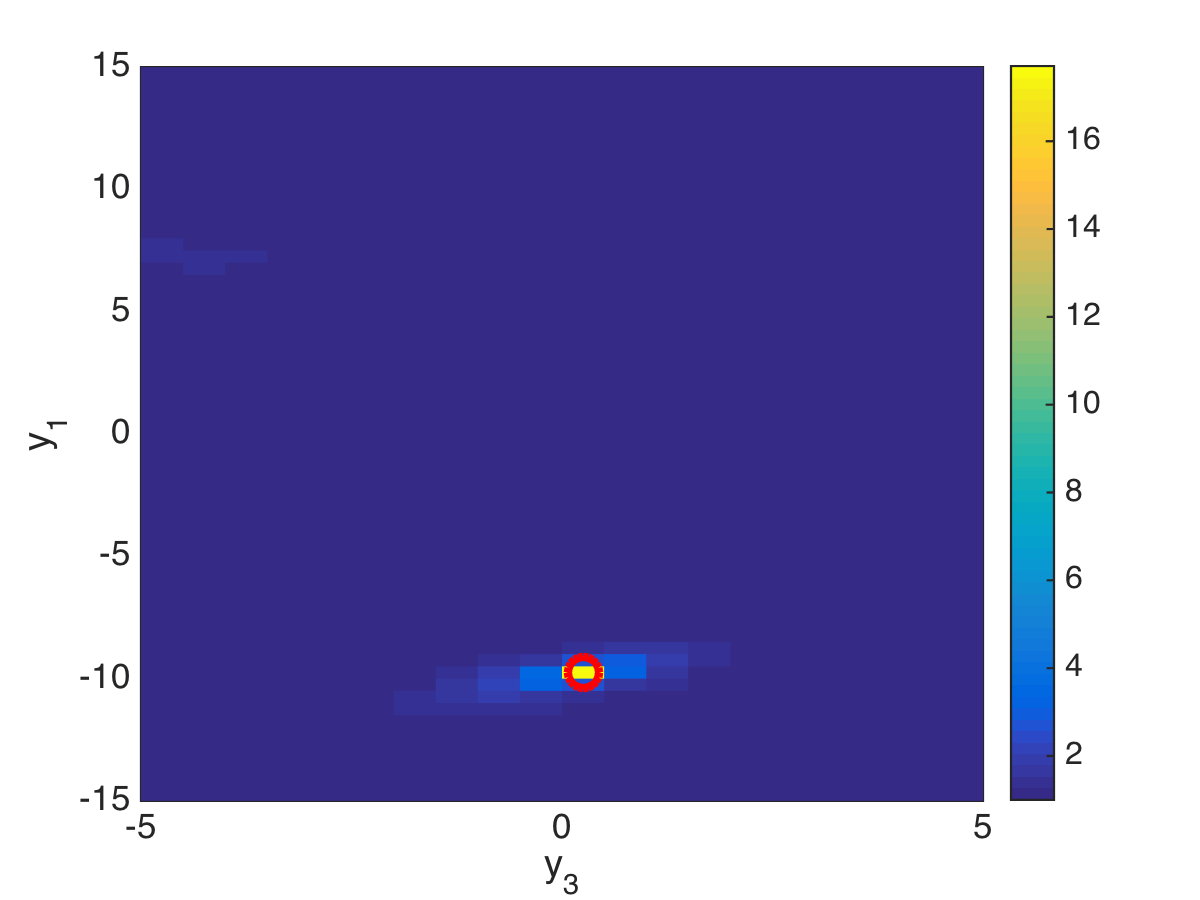}
\hspace{-0.2in}
\includegraphics[width=0.34\textwidth]{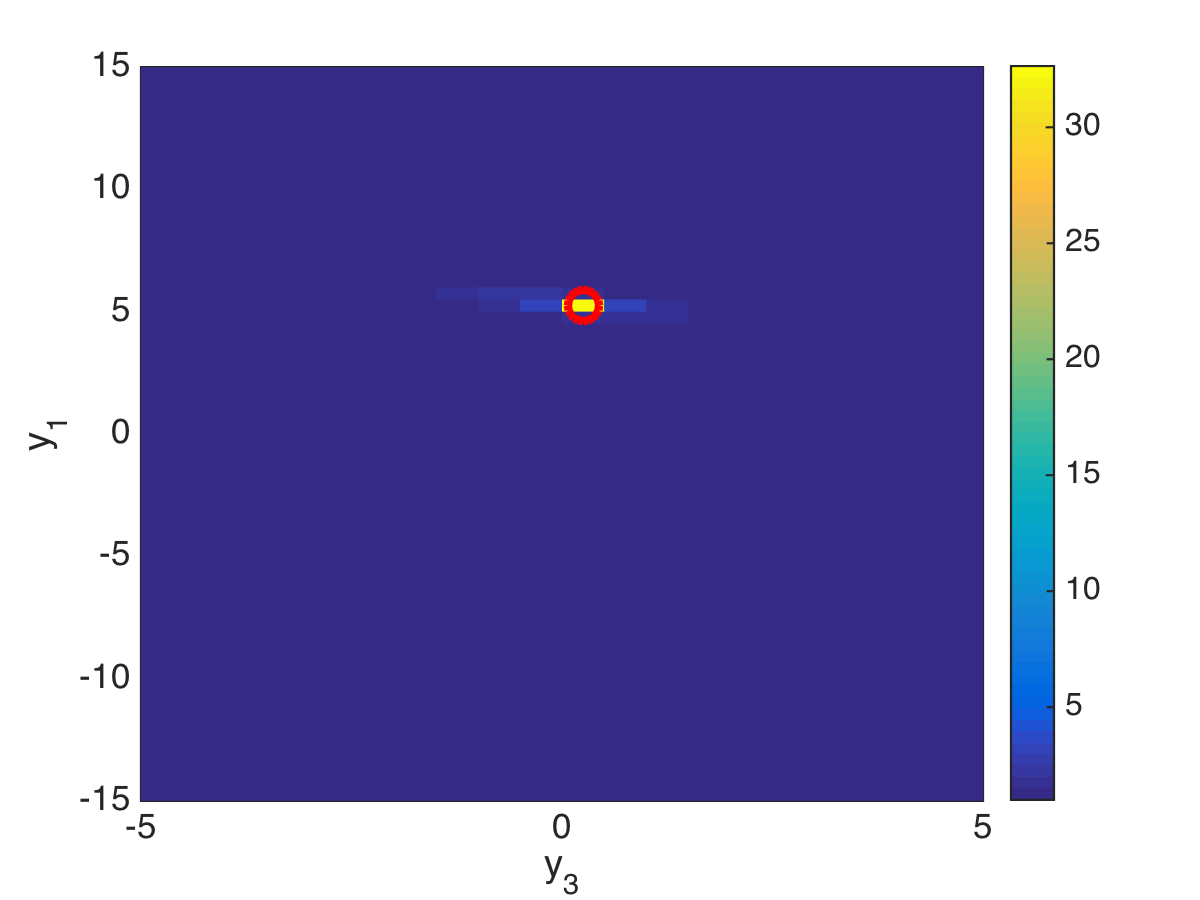} \\
\includegraphics[width=0.34\textwidth]{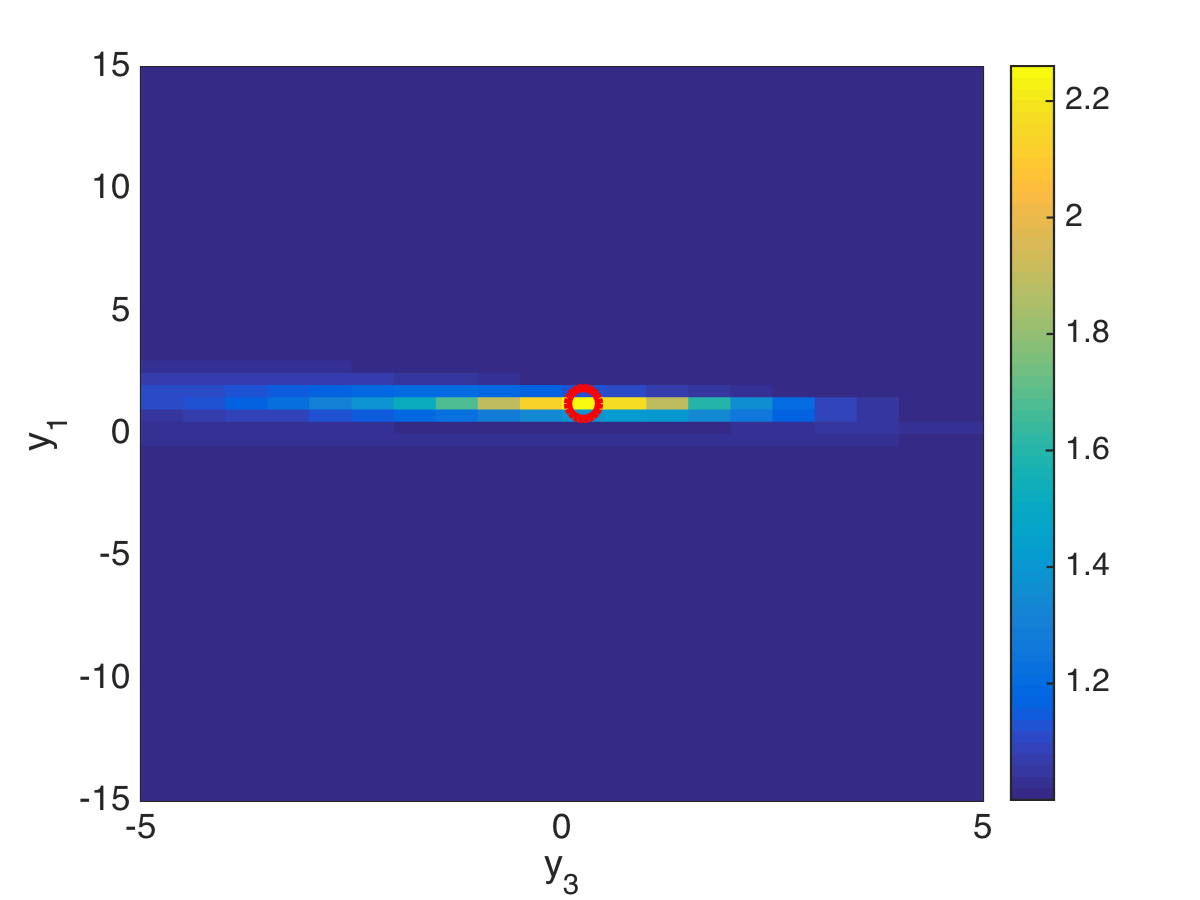}
\hspace{-0.2in}
\includegraphics[width=0.34\textwidth]{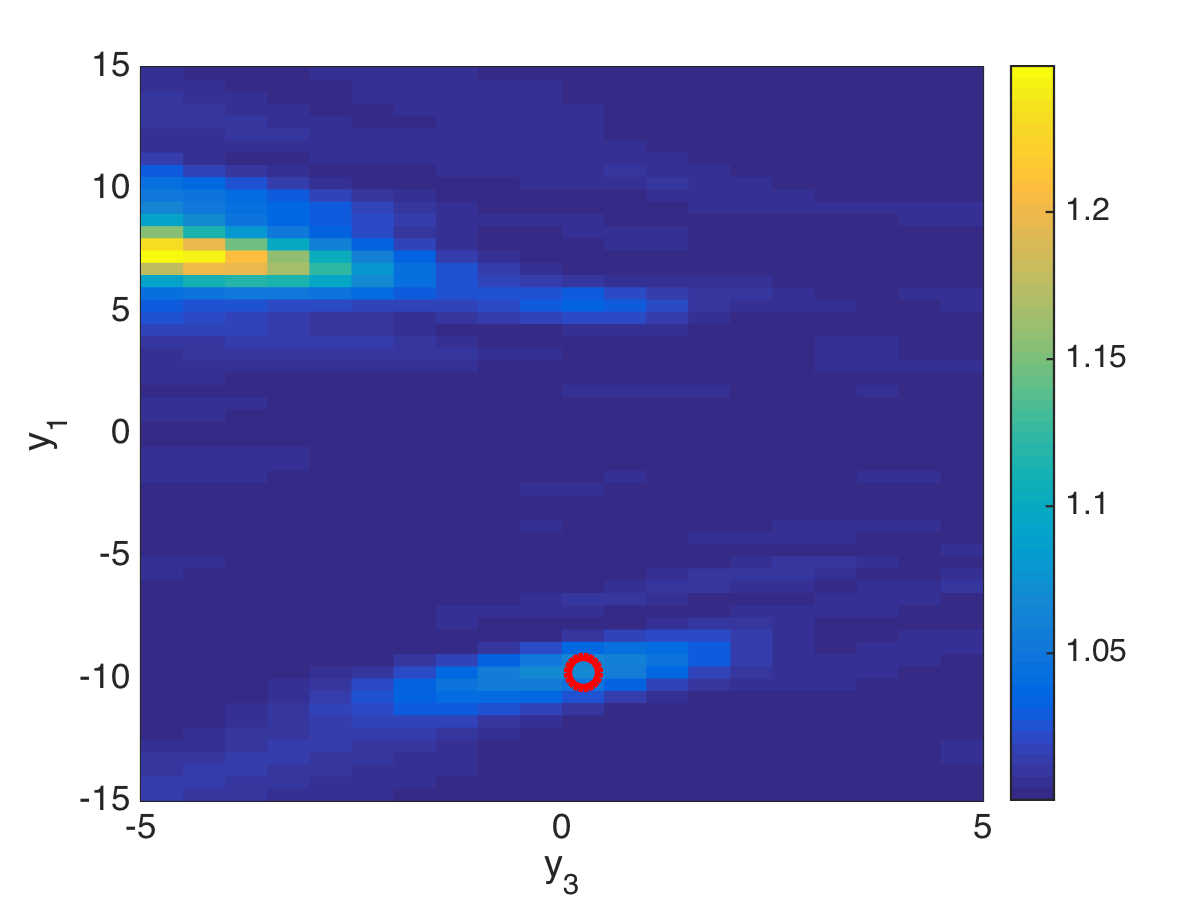}
\hspace{-0.2in}
\includegraphics[width=0.34\textwidth]{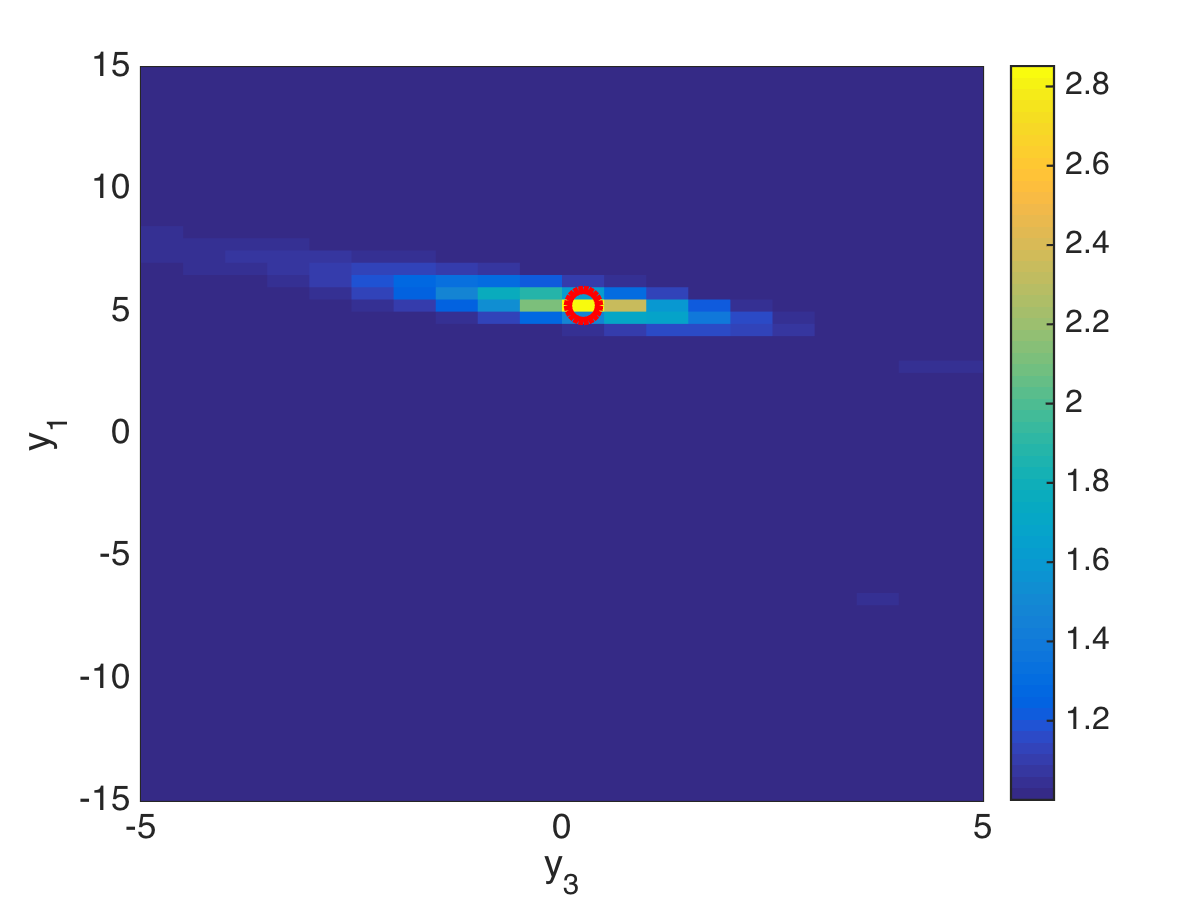}
\end{center}
\vspace{-0.1in}
\caption{Image of three inclusions in the large aperture regime, for complete measurements and $25\%$ noise (top) and 
$75\%$ noise (bottom). The displays are in the planes $y_2 = y_{p,2}$, for $p = 1,2,3$. The axes are in units of the wavelength.}
\label{fig:3T_NEAR_XZ}
\end{figure}

\subsubsection{Large aperture regime} We display in  Figure \ref{fig:3T_NEAR_XY} the imaging function \eqref{eq:IN19} in the plane 
$y_3 = L$. The image on the left is for data with $25\%$ noise,  where the effective rank of $\tDD$ is 
$\widetilde{\mathfrak{R}} = 6$. The right image is for $75\%$ noise, where the effective rank decreases to 
$\widetilde{\mathfrak{R}} = 3$.   The exact locations of the inclusions are  indicated with red circles, and they 
coincide with the peaks of the imaging function. As expected, the peaks are 
more prominent when $\widetilde{\mathfrak{R}}$ is closer to $\mathfrak{R}$ i.e., for the weaker noise. 
Note in particular that the inclusion at $\vy_2$ is barely visible at the $75\%$ noise level,  because the singular vector 
$\bh_1(\vy_2)$ of  $\GG(\vy_2)$ has a small projection on $\mbox{span}\{\bu_1,\bu_2, \bu_3\}$. 
The images in Figure \ref{fig:3T_NEAR_XZ} tell a similar story. They are displayed in the planes $y_2 = y_{p,2}$ for $p = 1, 2,3$, 
and show that the range resolution is worse than that in cross-range, specially at the higher noise level. 
Moreover, consistent with Figure \ref{fig:3T_NEAR_XY}, the inclusion at $\vy_2$ is difficult to localize 
at $75\%$ noise.

We present in Table \ref{tableNF}   the relative errors (in percent) 
of the estimation \eqref{eq:IN23} of the components of the reflectivities $\brho_p$, for $p = 1, 2,3$.   Because the 
inclusions are well localized, the errors of $\brho_p$ are of the form \eqref{eq:IN21}, and are mostly due to the noise. 
The  interaction between the 
inclusions is negligible, as  seen from the very accurate results at zero noise level.  This is expected 
because the inclusions are further apart than $\la L/a = \la$.
The errors increase with the noise, 
but they remain small for some components, mainly on the diagonal of the reflectivity tensors.
\begin{table}
\begin{center}
\begin{tabular}{c c|c | c | c | c | c | c}
 &&$\rho_{11}$ & $\rho_{22}$ & $\rho_{33}$ & $\rho_{12}$ & $\rho_{13}$ & $\rho_{33}$ \\
\hline 
&Inclusion 1 & 0.01 & 0.02 & 0.1 & 0.03 & 0.06 & 0.08 \\
$0\%$ noise & Inclusion 2 & 0.2 & 0.3 & 0.3 & 0.7 & 1.9 & 2.9 \\
& Inclusion 3 & 0.03 & 0.02 & 0.04 & 0.3 & 0.1 & 0.2 \\
\hline
& Inclusion 1 & 0.6 & 0.5 & 1 & 4.9 & 3.8 & 5.0 \\
$25\%$ noise & Inclusion 2 & 3.7 & 5.1 & 2.4 &10.7 & 19.3 & 44.7 \\
& Inclusion 3 & 2.9 & 1.2 & 7.7 & 48.0 & 21.7 & 21.7 \\
\hline
& Inclusion 1 & 3.3 & 1.8 & 5.7 & 7.1 & 31.3 & 9.9 \\
$75\%$ noise & Inclusion 2 & 7.4 & 5.2 & 8.9 & 4.7 & 34.7 & 98.3 \\
& Inclusion 3 &  2.1 & 0.7 & 7.4 & 45.2 & 18.1 & 49.4 \\
\hline 
\end{tabular}
\caption{Table of relative errors (in percent) of the estimated components of the reflectivity tensors in the large aperture regime.}
\end{center}
\label{tableNF}
\end{table}

\begin{figure}[h]
\begin{center}
\includegraphics[width=0.52\textwidth]{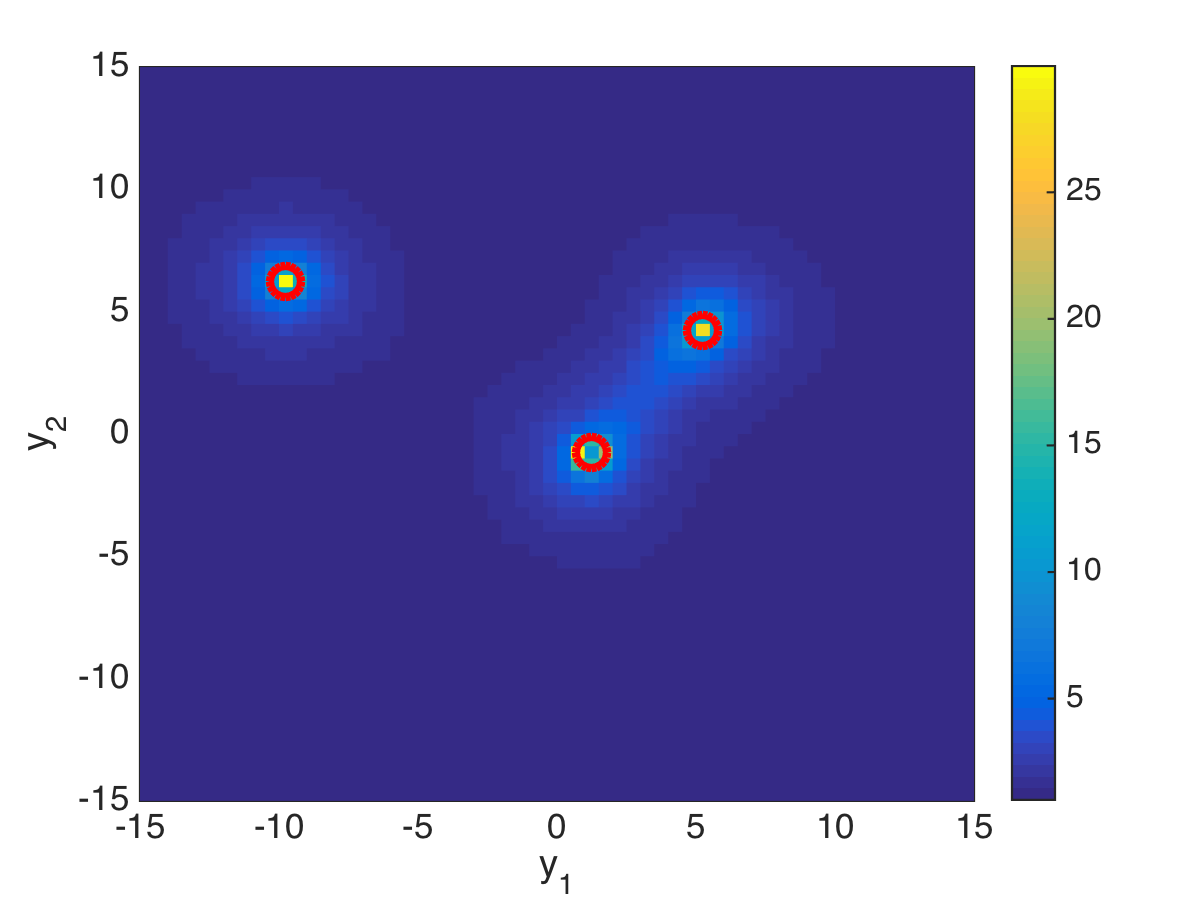}
\hspace{-0.25in}\includegraphics[width=0.52\textwidth]{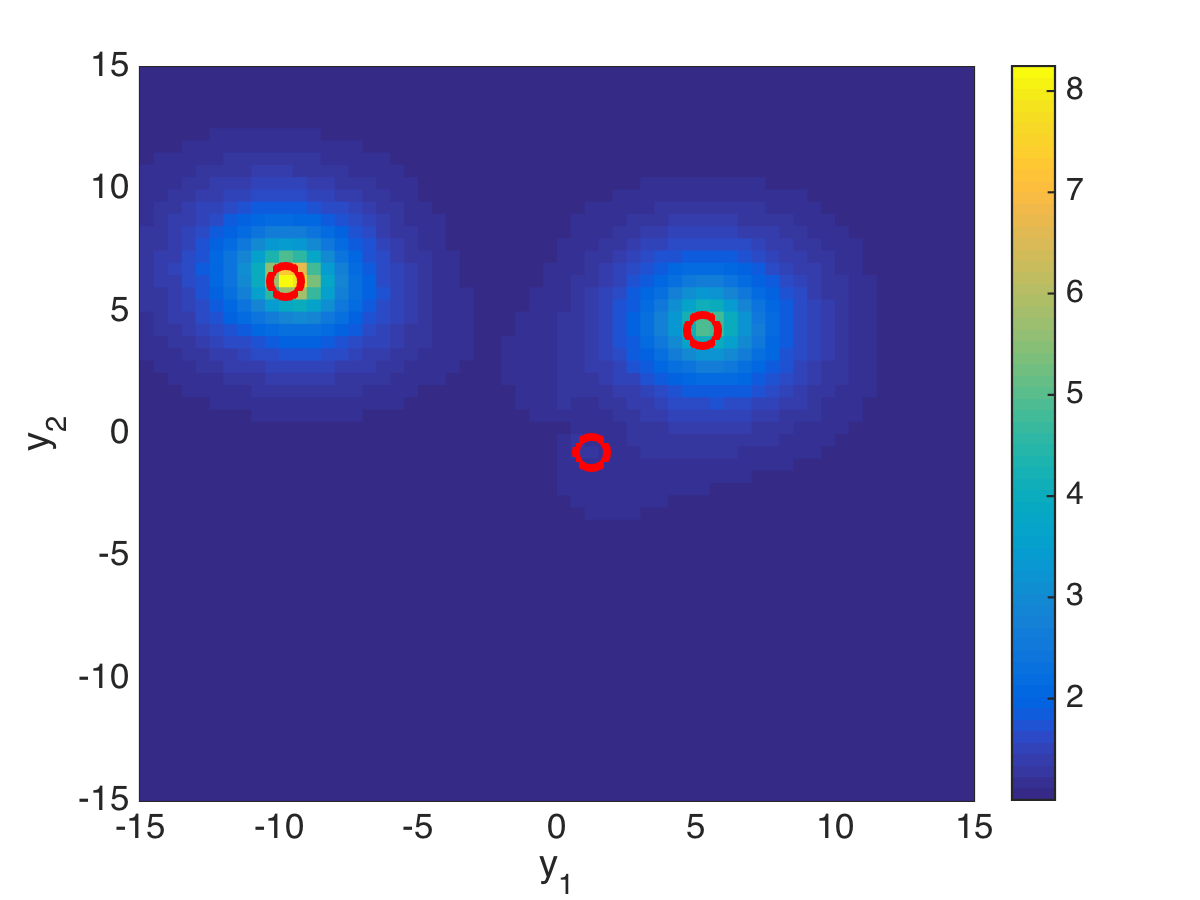}
\end{center}
\vspace{-0.1in}
\caption{Image of three inclusions in the large aperture regime, for complete measurements and $25\%$ noise (left) and 
$50\%$ noise (right). The display is in the plane $y_3 = L$. The axes are in units of the wavelength. The locations of the inclusions
are indicated with red circles.}
\label{fig:3T_FAR_XY}
\end{figure}

\subsubsection{Small aperture regime}
As we saw in the previous section,  imaging is more difficult  in the small aperture regime. Range localization 
is impossible even at small levels of noise. Moreover, the resolution of the left image in Figure \ref{fig:3T_FAR_XY},
obtained at $25\%$ noise, is worse than that in Figure \ref{fig:3T_NEAR_XY}. The inclusions interact in this case,
because they are separated by smaller distances  than $\la L/a = 10\la$, and this adds to the difficulty.
The right image in Figure \ref{fig:3T_FAR_XY} shows that 
one inclusion is obscured at $50\%$ noise.

We do not show the errors in the estimation of the reflectivity tensor  because they are bad even at 
small noise levels, due to the poor conditioning of the matrix 
$\bSigma(\vy)$ in \eqref{eq:defGamma}, and the interaction of the inclusions.

\section{Summary}
\label{sect:sum}

We introduced and analyzed  a robust methodology for detection, localization, and characterization 
of small electromagnetic inclusions from measurements of the time-harmonic electric field  corrupted by additive noise.
The methodology is motivated by the statistical analysis of the noisy data matrix gathered by an array of sensors. 
In particular, it uses random matrix theory results about low rank perturbations of large random matrices.
The detection of the inclusions is carried out by the inspection of the top singular values of the  data matrix.
The localization is done with  an imaging function that uses the predictable angles 
between the singular vectors of the noisy data matrix and of the noiseless matrix. The characterization of the 
inclusions amounts to estimating their reflectivity tensor which 
depends on the shape of their support and their electric permittivity.

The inversion methodology is robust with respect to a significant level of additive noise, much more 
than standard  imaging methods like MUSIC can handle.
We clarify that MUSIC localization fails because the singular vectors of the noisy data matrix are not collinear
to those of the unperturbed matrix. A main result of the paper is that the angles between these singular
vectors can be predicted and used to obtain a robust inclusion localization.
We also explain that the resolution and stability of inversion are dependent on each other,  and quantify 
the stability gain due to  larger array aperture 
and more than one component of the electric field being recorded and processed.

\section*{Acknowledgements}
Liliana Borcea's work was partially supported by AFOSR Grant
FA9550-15-1-0118.


\bibliographystyle{siam} \bibliography{NOISE}

\end{document}